\documentclass[a4paper,11pt,twoside]{article}
\usepackage{mfpaperstuff2}
\usepackage{mfabacus}
\usepackage{mdframed}
\usepackage{ifthen}

\makeatletter
\newsavebox{\@brx}
\newcommand{\llangle}[1][]{\savebox{\@brx}{\(\m@th{#1\langle}\)}%
  \mathopen{\copy\@brx\kern-0.7\wd\@brx\usebox{\@brx}}}
\newcommand{\rrangle}[1][]{\savebox{\@brx}{\(\m@th{#1\rangle}\)}%
  \mathclose{\copy\@brx\kern-0.7\wd\@brx\usebox{\@brx}}}
\makeatother

\newgeometry{margin=0.421in}
\begin{document}
\newcommand\hhh{\calh_}
\newcommand\spe[1]{\mathrm S^{#1}}
\newcommand\jms[1]{\mathrm D^{#1}}
\newcommand\schir[1]{\mathrm L^{#1}}
\newcommand\weyl[1]{\Delta^{#1}}
\newcommand\sym{\mathrm{Sym}}
\newcommand\cm[2]{[#1:#2]}
\newcommand\muset[2]{\{0^{#1},1^{#2}\}}
\newcommand\ev[1]{\operatorname{ev}(#1)}
\newcommand\sprm{spin-removable node\xspace}
\newcommand\sprms{spin-removable nodes\xspace}
\newcommand\spam{spin-addable node\xspace}
\newcommand\spams{spin-addable nodes\xspace}
\newcommand\sdr[1]{\sid{#1\dblreg}}
\newcommand\nextcase[1]{\addtocounter{case}1\subsection*{Case \arabic{case}: #1}}
\newcommand\thr{\scrx_3}
\newcommand\oni{\scrx_1}
\newcommand\thrt{\scry_3}
\newcommand\onit{\scry_1}
\newcommand\frst\scrz
\newcommand\ird\scri
\newcommand\trp{$2$-regular partition\xspace}
\newcommand\trps{$2$-regular partitions\xspace}
\newcommand\fbw{$4$-bar-weight\xspace}
\newcommand\fbc{$4$-bar-core\xspace}
\newcommand\vo[3]{\operatorname{sep}(#1,#2,#3)}
\newcommand\schur[1]{\cals(#1)}
\newcommand\qschur[1]{\cals^\circ(#1)}
\newcommand\odc{\mathring D}
\newcommand\decs[2]{D^{\mathsf{spn}}_{#1#2}}
\renewcommand\v{^{-1}}
\newcommand\sgn[1]{\operatorname{sgn}(#1)}
\newcommand\ke[2]{\spadesuit_{#1,#2}}
\newcommand\ip[2]{\lan#1,#2\ran}
\newcommand\qs{separated\xspace}
\newcommand\rs[2]{#1{\downarrow}_{#2}}
\newcommand\spr{spin residue\xspace}
\newcommand\sprs{spin residues\xspace}
\newcommand\dup[1]{#1\sqcup#1}
\newcommand\chm[2]{\left(#1\!:\!#2\right)}
\newcommand\funny{conjugate-symmetric\xspace}
\renewcommand\deg{\operatorname{deg}}
\renewcommand{\crefrangeconjunction}{--}
\medmuskip=3mu plus 1mu minus 1mu
\renewcommand{\rt}[1]{\rotatebox{90}{$#1$}}
\newcommand\aaa{\mathfrak{A}_}
\newcommand\taaa{\tilde{\mathfrak{A}}_}
\newcommand\tsss{\tilde{\mathfrak{S}}_}
\newcommand\len[1]{l(#1)}
\newcommand\res{{\downarrow}}
\newcommand\ind{{\uparrow}}
\newcommand\reg{^{\operatorname{reg}}}
\newcommand\dbl{^{\operatorname{dbl}}}
\newcommand\dblreg{^{\operatorname{dblreg}}}
\newcommand\ee[1]{\mathrm{e}_{#1}}
\newcommand\ff[1]{\mathrm{f}_{#1}}
\newcommand\eed[2]{\ee{#1}^{(#2)}}
\newcommand\ffd[2]{\ff{#1}^{(#2)}}
\newcommand\emx[1]{\eed{#1}{\operatorname{max}}}
\newcommand\fmx[1]{\ffd{#1}{\operatorname{max}}}
\newcommand\ord[1]{\llbracket#1\rrbracket}
\newcommand\spn[1]{\llangle#1\rrangle}
\newcommand\sid[1]{\varphi(#1)}
\newcommand\prj[1]{\operatorname{prj}(#1)}
\newcommand\prsj[1]{\prj{\sigma+2#1}}
\newcommand\modr[1]{\widebar{#1}}
\newcommand\mspn[1]{\modr{\spn{#1}}}
\newcommand\mord[1]{\modr{\ord{#1}}}
\newcommand\adr[2]{\ifthenelse{\equal{#2}1}{\stackrel{#1}\longrightarrow}{\stackrel{#1^#2}\longrightarrow}}
\newcommand\ads[2]{\ifthenelse{\equal{#2}1}{\stackrel{#1}\Longrightarrow}{\stackrel{#1^#2}\Longrightarrow}}

\Yvcentermath0
\Yboxdim{12pt}

\pdfpagewidth=12.5cm
\pdfpageheight=8.839cm
\setcounter{page}0
{\footnotesize This is the author's version of a work that was accepted for publication in the Proceedings of the London Mathematical Society. Changes resulting from the publishing process, such as peer review, editing, corrections, structural formatting, and other quality control mechanisms may not be reflected in this document. Changes may have been made to this work since it was submitted for publication. A definitive version was subsequently published in\\
\textit{Proc.\ London Math.\ Soc.} \textup{\textbf{116} (2018) 878--928.\\
http://dx.doi.org/10.1112/plms.12087}.\\
%The final publication is available at www.springerlink.com.\normalsize}
}
\newgeometry{margin=2cm,includehead,includefoot}

\title{Irreducible projective representations of the symmetric group which remain irreducible in characteristic $2$}
\runninghead{Irreducible projective representations of the symmetric group in characteristic $2$}

\toptitle

\begin{abstract}
For any finite group $G$ and any prime $p$ one can ask which ordinary irreducible representations remain irreducible in characteristic $p$. We answer this question for $p=2$ when $G$ is a proper double cover of the symmetric group. Our techniques involve constructing part of the decomposition matrix for a Rouquier block of a double cover, restricting to subgroups using the Brundan--Kleshchev modular branching rules and comparing the dimensions of irreducible representations via the bar-length formula.
\end{abstract}

\pdfpagewidth=210mm
\pdfpageheight=297mm

\tableofcontents

\renewcommand\baselinestretch{1.048}

\section{Introduction}

For any finite group $G$ and prime number $p$, it is an interesting question to ask which ordinary irreducible representations of $G$ remain irreducible in characteristic $p$. This question was answered some time ago for the symmetric groups \cite{jm1,jmp2,slred,mfred,mfirred}, and more recently \cite{mfaltred} for the alternating groups. In this paper we address double covers of the symmetric groups; since (in all but finitely many cases) these are Schur covers, we implicitly also address the question of which ordinary irreducible projective representations of the symmetric groups remain irreducible in characteristic $p$.

We denote a chosen double cover of $\sss n$ by $\tsss n$ (though in fact our results apply equally well to the other double cover); thus $\tsss n$ has a central involution $z$ such that $\tsss n/\lan z\ran\cong\sss n$. In any characteristic one can separate the irreducible representations of $\tsss n$ into ``linear'' representations (on which $z$ acts trivially) and ``spin'' representations (on which $z$ acts as $-1$). The irreducible linear representations are also irreducible representations of $\sss n$, and so the answer to our main question for linear representations of $\tsss n$ follows at once from the answer for $\sss n$. So we need only consider spin representations. In fact in this paper we prefer to work in the language of characters: the ordinary irreducible spin characters of $\tsss n$ are labelled by $2$-regular partitions of $n$: there is a single character $\spn\la$ for every $2$-regular partition $\la$ of $n$ with an even number of even parts, and a pair of associate characters $\spn\la_+,\spn\la_-$ for every $\la$ with an odd number of even parts. In the latter case, we write $\spn\la$ to mean ``either $\spn\la_+$ or $\spn\la_-$'', so our main question becomes ``for which $\la$ is the $p$-modular reduction of $\spn\la$ irreducible?''.

There has been a recent resurgence of interest in the representation theory of $\tsss n$ in odd characteristic. We hope to address this case in a future paper, but in the present paper we work instead in characteristic $2$. Apart from the fundamental paper \cite{bo} of Bessenrodt and Olsson in 1997 from which we take some important results, this situation has been relatively neglected, but we hope that the present paper will spark a revival. Note that in characteristic $2$ the central involution $z$ must act as $1$ on an irreducible representation, so there are no irreducible spin representations. So the irreducible Brauer characters in characteristic $2$ are precisely those coming from $\sss n$, namely the characters $\sid\mu$ of the ``James modules'' $\jms\mu$. Thus the $2$-modular reduction of any ordinary irreducible spin character decomposes as a sum of characters $\sid\mu$. This means that the decomposition matrix of $\tsss n$ in characteristic $2$ consists of the decomposition matrix of $\sss n$ extended downwards by a ``spin matrix'' with rows corresponding to ordinary irreducible spin characters. Although this spin matrix can be computed algorithmically given the decomposition matrix of $\sss n$, it remains fairly mysterious in general. A key contribution in the present paper is to describe part of this matrix (specifically, the part with rows labelled by $2$-regular partitions having no even parts greater than $2$) for Rouquier blocks of $\tsss n$.

We now indicate the structure of the paper and of the proof of our main theorem. In \cref{backsec} we set out the background information we shall need; we try to provide enough information to make the paper reasonably self-contained, though we do make some effort at concision. In \cref{mainthmsec} we recall the notion of a \emph{$2$-Carter partition} which enables us to state our main theorem. In \cref{degreesec} we examine character degrees, which are an important ingredient in our proof. The key observation here is that if $\spn\la$ and $\spn\mu$ are ordinary characters having a $2$-modular constituent in common, and if the degree of $\spn\la$ is greater than the degree of $\spn\mu$, then the $2$-modular reduction of $\spn\la$ cannot be irreducible. It turns out that (in contrast to the more famous hook-length formula for the degrees of linear characters) the bar-length formula for spin characters allows us to compare degrees quite easily. In \cref{rouqandsep} we consider spin characters lying in a ``Rouquier'' block of $\tsss n$. While Rouquier blocks of symmetric groups have been studied extensively in the last few years, this appears to be the first time they have been considered for the double covers in characteristic $2$. By constructing various projective characters and exploiting connections with symmetric functions, we describe part of the decomposition matrix for a Rouquier block. We then extend this beyond Rouquier blocks to spin characters labelled by what we call \emph{\qs} partitions, and give a proof of our main theorem for \qs partitions.  In \cref{mainproofsec} we complete the proof of our main theorem, by considering non-\qs partitions. Here we use an inductive proof (similar in spirit to that used for the alternating groups in \cite{mfaltred}), based on Kleshchev's modular branching rules. In \cref{futuresec} we consider possibilities for future work, and in \cref{indexnotnsec} we provide an index of notation.

The author would like to express his gratitude to the referee for a careful and thorough reading of the paper, and some very helpful comments.

The research presented in this paper would not have been possible without extensive computations using GAP \cite{gap}.

\section{Background}\label{backsec}

\emph{Throughout this section $\bbf$ is a fixed field of characteristic $2$.}

\subsection{Symmetric groups, Hecke algebras and Schur algebras}\label{snsec}

We let $\sss n$ denote the symmetric group of degree $n$, and $\bbf\sss n$ its group algebra over $\bbf$. The representation theory of $\bbf\sss n$ is well-studied, and in the last twenty years or so has become increasingly intertwined with the representation theory of a certain (Iwahori--)Hecke algebra. Specifically, we let $\hhh n$ denote the Iwahori--Hecke algebra of $\sss n$ over $\bbc$ with quantum parameter $-1$; this has generators $T_1,\dots,T_{n-1}$, with defining relations $(T_i+1)^2=0$ for each $i$, together with the braid relations $T_iT_{i+1}T_i=T_{i+1}T_iT_{i+1}$ for $1\ls i\ls n-2$ and $T_iT_j=T_jT_i$ for $1\ls i<j-1\ls n-2$. The representation theory of $\hhh n$ is closely related to the representation theory of $\bbf\sss n$: the canonical labelling sets for the Specht modules and simple modules for these algebras are the same, and (as we shall see) there is a close relationship between their decomposition matrices. The main advantage of introducing $\hhh n$ into the subject is that its decomposition numbers are better understood than those of $\bbf\sss n$ (indeed, they can be computed using the \emph{LLT algorithm} \cite{llt}). We shall exploit this later, although the relationship between the symmetric group and the Hecke algebra does not extend well to double covers of $\sss n$.

Another important family of algebras in this area is the family of \emph{Schur algebras} introduced by Green \cite{gbook}, and more generally the \emph{$q$-Schur algebras} defined by Dipper and James \cite[Definition 2.9]{dj}. We just concentrate on two of these algebras, for each $n$: let $\schur n$ denote the classical Schur algebra over $\bbf$ (which is $\cals_\bbf(1,n)$ in the notation of Dipper--James), and let $\qschur n$ denote the $q$-Schur algebra over $\bbc$ with quantum parameter $q=-1$, i.e.\ $\cals_\bbc(-1,n)$ in the Dipper--James notation. $\schur n$ is related to $\bbf\sss n$ via the \emph{Schur functor}, which shows in particular that the decomposition matrix of $\bbf\sss n$ can be obtained from that of $\schur n$ by deleting some of the columns. $\qschur n$ is related to $\hhh n$ is a similar way. As we shall see, the decomposition matrix of $\schur n$ also appears in the adjustment matrices of Rouquier blocks.

\subsection{Double covers of $\sss n$ and projective representations}\label{doublecoversec}

Now we consider double covers. In this paper we work with the group $\tsss n$ with generators $s_1,\dots,s_{n-1},z$, subject to relations
\begin{align*}
z^2&=1,\\
zs_i&=s_iz\tag*{for $1\ls i\ls n-1$,}\\
s_i^2&=1\tag*{for $1\ls i\ls n-1$,}\\
s_is_{i+1}s_i&=s_{i+1}s_is_{i+1}\tag*{for $1\ls i\ls n-2$,}\\
s_is_j&=zs_js_i\tag*{for $1\ls i\ls j-2\ls n-3$.}
\end{align*}
Then $\tsss n$ is a double cover of the symmetric group $\sss n$: the subgroup $\lan z\ran$ is central, and setting $z$ to $1$ in the above presentation yields the Coxeter presentation of $\sss n$.

In fact, $\tsss n$ is a \emph{Schur cover} of $\sss n$ for $n\gs4$. This means that linear representations of $\tsss n$ are equivalent to projective representations of $\sss n$, i.e.\ homomorphisms from $\sss n$ to $\mathrm{PGL}(V)$ for a vector space $V$. Thus our main theorem can be phrased in terms of projective representations. For $n\ls3$, $\sss n$ is its own Schur cover (note that the final relation in the presentation above never occurs when $n\ls3$, so that $\tsss n$ is just the direct product of $\sss n$ and $\lan z\ran$). Nevertheless, it makes sense to work with the group $\tsss n$ for all $n\gs0$, and we do so in this paper to help with our induction proofs.

We remark that for $n\gs4$ there is an alternative Schur cover obtained by replacing the relation $s_i^2=1$ with $s_i^2=z$; this is isomorphic to $\tsss n$ only when $n=6$, but in general is isoclinic to $\tsss n$, and so from our point of view its representation theory is the same.

Now we consider the different types of irreducible representations of $\tsss n$. In an irreducible representation of $\tsss n$ over any field, the central element $z$ must act as $1$ or $-1$. Irreducible representations in which $z$ acts as $1$ correspond exactly to irreducible representations of $\sss n$, and we refer to these as \emph{linear} representations. Studying linear representations of $\tsss n$ is equivalent to studying representations of $\sss n$, so in this paper we shall be chiefly concerned with the representations where $z$ acts as $-1$, which we call \emph{spin} representations.

For an introduction to the representation theory of $\tsss n$, we refer the reader to the book by Hoffman and Humphreys \cite{hh}.

\subsection{Partitions}\label{partnsec}

The representations of the symmetric groups and their double covers are conventionally labelled by partitions. In this paper, a \emph{partition} is an infinite weakly decreasing sequence $\la=(\la_1,\la_2,\dots)$ of non-negative integers with finite sum, which we denote $|\la|$. If $|\la|=n$, then we say that $\la$ is a partition of $n$. The integers $\la_1,\la_2,\dots$ are called the \emph{parts} of $\la$. When writing partitions, we conventionally omit the trailing zeroes and group together equal parts with a superscript (so that we write the partition $(4,4,3,1,1,1,0,0,0,\dots)$ as $(4^2,3,1^3)$), and we write the unique partition of $0$ as $\varnothing$. A partition is \emph{$2$-regular} if it does not have two equal non-zero parts, and \emph{$2$-singular} otherwise. We write $\calp(n)$ for the set of partitions of $n$, and $\cald(n)$ for the set of \trps of $n$, and set $\calp=\bigcup_{n\gs0}\calp(n)$ and $\cald=\bigcup_{n\gs0}\cald(n)$.

If $\la\in\calp$, the \emph{conjugate partition} $\la'$ is defined by setting
\[
\la'_r=\left|\lset{k\gs1}{\la_k\gs r}\right|
\]
for all $r$. The \emph{dominance order} on $\calp(n)$ is defined by saying that $\la$ \emph{dominates} $\mu$ (and writing $\la\dom\mu$) if
\[
\la_1+\dots+\la_r\gs\mu_1+\dots+\mu_r
\]
for all $r\gs1$.

The \emph{Young diagram} of a partition $\la$ is the set
\[
\lset{(r,c)}{r\gs1,\ 1\ls c\ls\la_r}\subset\bbn^2.
\]
We often identify a partition with its Young diagram; for example, we may write $\la\subseteq\mu$ to mean that $\la_i\ls\mu_i$ for all $i$. The elements of the Young diagram of $\la$ are called the \emph{nodes} of $\la$.  More generally, we use the term \emph{node} for any element of $\bbn^2$. We draw Young diagrams as arrays of boxes using the English convention, in which $r$ increases down the page and $c$ increases from left to right.

A node of $\la$ is \emph{removable} if it can be removed from $\la$ to leave the Young diagram of a smaller partition, while a node not in $\la$ is an \emph{addable node} of $\la$ if it can be added to $\la$ to leave the Young diagram of a larger partition. The \emph{residue} of a node $(r,c)$ is the residue of $c-r$ modulo $2$. We refer to a node of residue $i$ as an $i$-node.

Now we introduce some notions relating to spin representations. The \emph{\spr} of a node $(r,c)$ is the residue of $\inp{c/2}$ modulo $2$; that is, $0$ if $c\equiv0$ or $1\ppmod4$, and $1$ otherwise. If $\la$ is a \trp and $i\in\{0,1\}$, we write $\rs\la i$ for the smallest \trp which can be obtained from $\la$ by removing nodes of \spr $i$; the nodes removed are called \emph{$i$-\sprms} of $\la$.

For example, let $\la=(9,6,4,3,1)$. Then $\rs\la0=(7,6,4,3)$ and $\rs\la1=(9,5,4,2,1)$. The following diagram shows the \sprs of the nodes of $\la$, with the \sprms highlighted.
\newcommand\thk{\Ylinethick{1.3pt}}
\newcommand\thn{\Ylinethick{.3pt}}
\[
\young(0110011!\thk00!\thn,01100!\thk1!\thn,0110,01!\thk1,0!\thn)
\]

\smallskip
\begin{mdframed}[innerleftmargin=3pt,innerrightmargin=3pt,innertopmargin=3pt,innerbottommargin=3pt,roundcorner=5pt,innermargin=-3pt,outermargin=-3pt]
\noindent\textbf{Some notation for partitions}

In this paper we shall use some notation for partitions which, although quite natural, is not particularly standard.

Suppose $\la$ and $\mu$ are partitions and $a\in\bbq$ such that $a\la_i\in\bbz_{\gs0}$ for all $i$. Then we write
\begin{itemize}
\item
$a\la$ for the partition $(a\la_1,a\la_2,\dots)$;
\item
$\la+\mu$ for the partition $(\la_1+\mu_1,\la_2+\mu_2,\dots)$;
\item
$\la\sqcup\mu$ for the partition obtained by combining all the parts of $\la$ and $\mu$ and arranging them in decreasing order.
\end{itemize}
Furthermore, we may use these operations simultaneously, and they take precedence in the order they appear above, so that $\la\sqcup a\mu+\nu$ means $\la\sqcup((a\mu)+\nu)$.

For example, if $\la=(11,7,3)$, $\mu=(3,1^2)$ and $\nu=(10,2)$, then we can define the partition
\[
\la+4\mu\sqcup\nu=(23,11,10,7,2).
\]
\end{mdframed}

\subsection{$2$-cores and $2$-quotients}\label{2coresec}

We now set out some of the combinatorics of partitions relevant to the $2$-modular representation theory of the symmetric group. Everything we shall say has analogues for all primes $p$, but we stick to $p=2$ in the interest of brevity.

Suppose $\la\in\calp$. A \emph{rim $2$-hook} of $\la$ is a pair of (horizontally or vertically) adjacent nodes of $\la$ which may be removed to leave a smaller Young diagram; we call the rim hook \emph{horizontal} or \emph{vertical} accordingly. The \emph{$2$-core} of $\la$ is the partition obtained by repeatedly removing rim $2$-hooks until none remain. The $2$-core is well defined, and has the form $(c,c-1,\dots,1)$ for some $c\gs0$. The \emph{$2$-weight} of $\la$ is the number of rim $2$-hooks removed to reach the $2$-core.

This may be visualised using the abacus. Take an abacus with two infinite vertical runners, and label these $0$ and $1$ (with $0$ being the runner on the left). On runner $a$ mark positions $\dots,a-4,a-2,a,a+2,a+4,\dots$, so that position $2i+1$ is directly to the right of position $2i$ for all $i$. Now given a partition $\la$, place a bead on the abacus at position $\la_r-r$ for each $r$. The resulting configuration is called the \emph{abacus display} for $\la$. The key observation motivating the abacus is that removing a rim $2$-hook corresponds to sliding a bead into an empty position immediately above. Hence an abacus display for the $2$-core of $\la$ may be obtained by repeatedly sliding beads up until every bead has a bead immediately above it.

The abacus also allows us to define the $2$-quotient of a partition. Given the abacus display for $\la$ and given $a\in\{0,1\}$, examine runner $a$ in isolation, and let $\la^{(a)}_i$ be the number of empty positions above the $i$th lowest position, for each $i$. Then $\la^{(a)}=(\la^{(a)}_1,\la^{(a)}_2,\dots)$ is a partition, and the pair $(\la^{(0)},\la^{(1)})$ is the \emph{$2$-quotient} of $\la$.

\begin{eg}
The partition $\la=(6,4^2,3)$ has $2$-core $(2,1)$, as we see from the Young diagram (in which we indicate removed rim $2$-hooks) and the abacus displays of these two partitions.
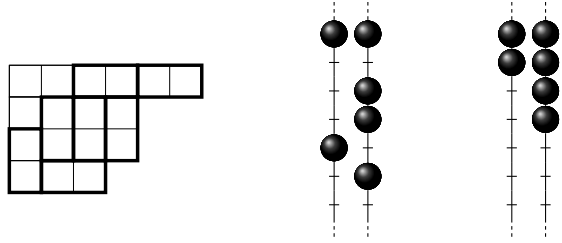
\[
\tikz[baseline=-2cm]{\tyng(0cm,0cm,6,4^2,3)\Yfillopacity{0}\thk\tgyoung(0cm,0cm,::_2_2,:|2|2|2:,|2:,:_2)}\qquad\qquad
\abacus(vv,bb,nn,nb,nb,bn,nb,nn,vv)\qquad\qquad\abacus(vv,bb,bb,nb,nb,nn,nn,nn,vv)
\]
Examining the runners in the display for $\la$, we see that the $2$-quotient of $\la$ is $((3),(2,1^2))$.
\end{eg}

Later we shall need the following simple lemma (a version of this for arbitrary $p$ appears in \cite[\S1.1.1]{mfirred}, among other places).

\begin{lemma}\label{quoconj}
Suppose $\mu\in\calp$, and let $(\mu^{(0)},\mu^{(1)})$ be the $2$-quotient of $\mu$. Then the $2$-quotient of $\mu'$ is $((\mu^{(1)})',(\mu^{(0)})')$.
\end{lemma}

The final thing we need relating to cores and quotients is the \emph{$2$-sign} of a partition. This was originally introduced by Littlewood \cite[p.~338]{litt}, but we use the alternative definition from \cite[p.~229]{j10}. Suppose $\mu$ is a partition, and construct the $2$-core of $\mu$ by repeatedly removing rim $2$-hooks. Let $a$ be the number of vertical hooks removed. $a$ is not well-defined, but its parity is, and so we may safely define the $2$-sign $\epsilon(\mu)=(-1)^a$ of $\mu$.% $\epsilon(\mu)$ is then called the $2$-sign of $\mu$.

We will need the following very simple observation.

\begin{lemma}\label{parityconj}
Suppose $\mu$ is a partition of $2$-weight $w$. Then $\epsilon(\mu')=(-1)^w\epsilon(\mu)$.
\end{lemma}

\subsection{Irreducible modules, irreducible characters and decomposition numbers}\label{irrdecsec}

Having established the necessary background relating to the combinatorics of partitions, we can discuss modules and decomposition numbers for the groups and algebras studied in this paper.

For each $\la\in\calp(n)$, we write $\spe\la$ for the \emph{Specht module} for $\sss n$, as defined (over an arbitrary field) by James \cite{jbook}; in particular, $\spe{(n)}$ is the trivial module. Working over our field $\bbf$ of characteristic $2$, suppose $\la\in\cald(n)$; then $\spe\la$ has a unique irreducible quotient $\jms\la$, which we call the \emph{James module}. The modules $\jms\la$ give all the irreducible $\bbf\sss n$-modules as $\la$ ranges over $\cald(n)$.

A similar situation applies for the Hecke algebra $\hhh n$: for each $\la\in\calp(n)$ there is a Specht module for $\hhh n$ (which we shall also denote $\spe\la$), which has a unique irreducible quotient $\jms\la$ when $\la$ is $2$-regular, and the modules $\jms\la$ are the only irreducible $\hhh n$-modules up to isomorphism.

Now we consider the Schur algebras. For every $\la\in\calp(n)$ there is a \emph{Weyl module} $\weyl\la$ for $\schur n$, which has a unique irreducible quotient $\schir\la$. Note that the labelling we use for these modules is that used in James's paper \cite{j10} (and is the opposite of that in most other places), so that the Specht module $\spe\la$ is the image of $\weyl\la$ (and not $\weyl{\la'}$) under the Schur functor. Similarly for the $(-1)$-Schur algebra $\qschur n$ we have a Weyl module $\weyl\la$ and an irreducible module $\schir\la$ for every $\la\in\calp(n)$.

With these conventions established, we can define decomposition numbers. For any $\la,\mu\in\calp(n)$, we define $D_{\la\mu}$ to be the composition multiplicity $\cm{\weyl\la}{\schir\mu}$ for the Schur algebra $\schur n$. The fact that $\schur n$ is a quasi-hereditary algebra gives the following.

\begin{propn}\label{unitriang}
Suppose $\la,\mu\in\calp(n)$. Then $D_{\la\la}=1$, and $D_{\la\mu}=0$ unless $\mu\dom\la$.
\end{propn}

We emphasise that we use the convention from \cite{j10} for labelling Weyl modules; with the more usual convention the symbol $\dom$ would become $\domby$ in the above \lcnamecref{unitriang}. This convention also means that when $\mu\in\cald(n)$, $D_{\la\mu}$ is also the composition multiplicity $\cm{\spe\la}{\jms\mu}$ for $\sss n$ in characteristic $2$. The matrix with entries $(D_{\la\mu})_{\la,\mu\in\calp(n)}$ is the \emph{decomposition matrix} for $\schur n$, and the matrix with entries $(D_{\la\mu})_{\la\in\calp(n),\mu\in\cald(n)}$ is the decomposition matrix for $\sss n$ in characteristic $2$. We may write either of these matrices simply as $D$, if the context is clear.

A similar situation applies for Hecke algebras and $(-1)$-Schur algebras. For $\la,\mu\in\calp(n)$ we write $\odc_{\la\mu}$ for the composition multiplicity $\cm{\weyl\la}{\schir\mu}$ for the $(-1)$-Schur algebra $\qschur n$; if $\mu\in\cald(n)$, then this also equals the multiplicity $\cm{\spe\la}{\jms\mu}$ for the Hecke algebra $\hhh n$. In addition, \cref{unitriang} holds with $D$ replaced by $\odc$. Furthermore, the decomposition numbers $D_{\la\mu}$ and $\odc_{\la\mu}$ are related by the theory of adjustment matrices: if we fix $n$ and let $D$ denote the decomposition matrix of $\schur n$ and $\odc$ the decomposition matrix of $\qschur n$, then the square matrix $A$ defined by $D=\odc A$ is called the \emph{adjustment matrix} for $\schur n$. $A$ is automatically lower unitriangular (i.e.\ \cref{unitriang} holds with $D$ replaced by $A$), and remarkably $A$ has non-negative integer entries; this statement appears in \cite[Theorem 6.35]{mathbook}, and derives from the theory of decomposition maps \cite{geck2}. By taking only the rows and columns of $A$ labelled by \trps we obtain the adjustment matrix for $\sss n$ in characteristic $2$, which we also call $A$; then if we let $D,\odc$ denote the decomposition matrices for $\sss n$ and $\hhh n$, we again have $D=\odc A$.

Now we come to double covers. Here (and for the rest of the paper) we work with characters, rather than modules. This is because we do not have particularly nice constructions of the modules affording the irreducible characters of double covers, and because statements about the decomposition of projective modules are more readily expressed in terms of projective characters. Since we are only concerned with composition factors and not with module structures, there is no loss in working with characters.

With this in mind, we introduce notation for the characters of the $\sss n$-modules described above: let $\ord\la$ denote the character of the Specht module $\spe\la$ over $\bbc$, so that $\lset{\ord\la}{\la\in\calp(n)}$ is the set of ordinary irreducible characters of $\sss n$. For $\la\in\cald(n)$, let $\sid\la$ denote the $2$-modular Brauer character of the James module $\jms\la$. Since $\sss n$ is a quotient of $\tsss n$, the characters $\ord\la$ and $\sid\la$ naturally become characters of $\tsss n$, and we shall always view them in this way. Note that in characteristic $2$, the central element $z$ lies in the kernel of any irreducible representation of $\tsss n$, so the irreducible representations of $\tsss n$ are precisely those arising from the irreducible representations of $\sss n$. Thus $\lset{\sid\la}{\la\in\cald(n)}$ is the set of irreducible $2$-modular Brauer characters of $\tsss n$.

Now we look at spin characters. For each $\la\in\cald(n)$, let $\ev\la$ denote the number of positive even parts of $\la$, and write
\[
\cald^+(n)=\lset{\la\in\cald(n)}{\ev\la\text{ is even}},\qquad\cald^-(n)=\lset{\la\in\cald(n)}{\ev\la\text{ is odd}}.
\]
Then for each $\la\in\cald^+(n)$, there is an irreducible spin character $\spn\la$ of $\tsss n$, while for each $\la\in\cald^-(n)$, there are two irreducible spin characters $\spn\la_+$ and $\spn\la_-$ of $\tsss n$. The definition of these characters goes back to Schur, who showed that they are all the irreducible spin characters of $\tsss n$. A uniform construction of irreducible representations affording these characters was given much later, by Nazarov \cite{naz}. We do not give the representations or their characters explicitly in this paper; the properties we need will be summarised in later sections.

In fact, for $\la\in\cald^-(n)$ the characters $\spn\la_+$ and $\spn\la_-$ behave very similarly: their values differ only up to sign, and in particular their degrees are the same and their $2$-modular reductions are the same. So to make things simpler in this paper we adopt the convention that we write $\spn\la$ to mean ``either $\spn\la_+$ or $\spn\la_-$'' when $\la\in\cald^-(n)$.

We also follow the usual convention of omitting brackets from a partition when using the notation $\ord\ $ or $\spn\ $; thus we may write $\ord{4,2^2,1}$ instead of $\ord{(4,2^2,1)}$. We remark that historically the spin characters have been denoted $\lan\la\ran$ or $\lan\la\ran_\pm$; we use slightly different notation here because the brackets $\lan\ \ran$ are already used for two other purposes in this paper. We trust that there will be no confusion with the use of $\spn\mu$ in \cite[\S4]{bmo} to denote a projective character of $\tsss n$ in characteristic $3$.

We write $\chm{\ \ }{\ \ }$ for the usual inner product on ordinary characters, with respect to which the irreducible characters are orthonormal. For any character $\chi$ of $\tsss n$, we write $\modr\chi$ for the $2$-modular reduction of $\chi$, i.e.\ the $2$-modular Brauer character obtained by restricting $\chi$ to the $2$-regular conjugacy classes of $\tsss n$. Our main question is then: for which ordinary irreducible characters $\chi$ is $\modr\chi$ irreducible?

We write $\cm{\ \ }{\ \ }$ for the inner product on $2$-modular Brauer characters with respect to which the irreducible characters are orthonormal. Hence the decomposition number $D_{\la\mu}$ defined above equals $\cm{\mord\la}{\sid\mu}$, for $\la\in\calp(n)$ and $\mu\in\cald(n)$. Given $\la,\mu\in\cald(n)$, we write $\decs\la\mu=\cm{\mspn\la}{\sid\mu}$. The matrices of integers $D_{\la\mu}$ and $\decs\la\mu$ together constitute the decomposition matrix of $\tsss n$ in characteristic $2$.

\begin{eg}
The decomposition matrix of $\tsss5$ in characteristic $2$ is given below. The top part of the matrix is the decomposition matrix of $\sss5$, with the $(\la,\mu)$-entry being $D_{\la\mu}$. The lower part gives the decomposition numbers for spin characters, with entries $\decs\la\mu$; note that for $\la\in\cald^-(5)$ there are identical rows corresponding to $\spn{\la}_+$ and $\spn{\la}_-$. (In all matrices given explicitly in this paper, we use a dot to mean $0$.)
\[
\begin{array}{r|ccc|}
&\rt{(5)}&\rt{(4,1)}&\rt{(3,2)}\\\hline
(5)\phantom{_-}&1&\cdot&\cdot\\
(4,1)\phantom{_-}&\cdot&1&\cdot\\
(3,2)\phantom{_-}&1&\cdot&1\\
(3,1^2)\phantom{_-}&2&\cdot&1\\
(2^2,1)\phantom{_-}&1&\cdot&1\\
(2,1^3)\phantom{_-}&\cdot&1&\cdot\\
(1^5)\phantom{_-}&1&\cdot&\cdot\\
\hline
(5)\phantom{_-}&\cdot&\cdot&1\\
(4,1)_+&2&\cdot&1\\
(4,1)_-&2&\cdot&1\\
(3,2)_+&\cdot&1&\cdot\\
(3,2)_-&\cdot&1&\cdot\\\hline
\end{array}
\]
\end{eg}

\begin{rmk}
This example shows an interesting feature of the answer to our main question: We can re-cast our main question and ask ``which irreducible Brauer characters arise as the $2$-modular reductions of ordinary irreducible characters?''.  In the case $n=5$ we see that the Brauer character $\sid{3,2}$ does not arise as $\modr{\ord\la}$ for any $\la$, but does arise as $\mspn5$. So introducing spin characters can provide constructions of irreducible Brauer characters which might otherwise be hard to obtain.
\end{rmk}

A central result in modular representation theory (which we shall use without comment) is Brauer reciprocity, which gives a connection between decomposition numbers and indecomposable projective characters. Given $\mu\in\cald(n)$, the projective cover of the James module $\jms\mu$ may be lifted to an ordinary representation of $\tsss n$, and we write $\prj\mu$ for the character of this representation; this is called an \emph{indecomposable} projective character, and the characters $\prj\mu$ for $\mu\in\cald(n)$ give a basis for the space spanned by all projective characters (i.e.\ characters which vanish on $p$-singular elements of $\tsss n$).

Brauer reciprocity says that $\prj\mu$ is given in terms of irreducible characters by the entries in the column of the decomposition matrix corresponding to $\mu$, i.e.
\[
\prj\mu=\sum_{\la\in\calp(n)}D_{\la\mu}\ord\la\ \ {+}\sum_{\la\in\cald^+(n)}\decs\la\mu\spn\la\ \ {+}\sum_{\la\in\cald^-(n)}\decs\la\mu(\spn\la_++\spn\la_-).
\]
We will use this extensively in \cref{rouqandsep} to derive information on decomposition numbers.

\subsection{The inverse of the decomposition matrix of the Schur algebra}\label{inversesec}

Now we state some results which we shall need later concerning the inverses of the matrices $D$, $\odc$ and $A$. These are taken from James's seminal paper \cite{j10}, and derive ultimately from Steinberg's tensor product theorem.

Given partitions $\alpha,\mu$ with $|\mu|=2|\alpha|$, let $(\mu^{(0)},\mu^{(1)})$ be the $2$-quotient of $\mu$, and define
\[
\kappa(\alpha,\mu)=
\begin{cases}
a^\alpha_{\mu^{(0)}\mu^{(1)}}&(\text{if the $2$-core of $\mu$ is $\varnothing$})\\
0&(\text{otherwise}).
\end{cases}
\]
(This is a special case of \cite[Definition 2.13]{j10}.) Here (and henceforth) $a^\alpha_{\beta\gamma}$ denotes the Littlewood--Richardson coefficient corresponding to partitions $\alpha,\beta,\gamma$ with $|\alpha|=|\beta|+|\gamma|$.

Now we can state two results on the inverse of the decomposition matrix. We begin with a result for even $n$; this is a special case of \cite[Corollary~6.9]{j10}. Recall that $\epsilon(\mu)$ denotes the $2$-sign of a partition $\mu$, defined in \cref{2coresec}.

\begin{propn}\label{steinberg}
Suppose $\la,\mu\in\calp(n)$ and all the columns of $\la$ are of even length, and write $\la=\dup\alpha$. Then
\begin{align*}
\odc_{\la\mu}\v&=(-1)^{n/2}\epsilon(\mu)\kappa(\alpha,\mu)\\
\intertext{and}
D_{\la\mu}\v&=(-1)^{n/2}\epsilon(\mu)\sum_{\beta\in\calp(n/2)}D_{\alpha\beta}\v\kappa(\beta,\mu).
\end{align*}
\end{propn}

The reader following the reference to \cite{j10} may find Corollary~6.9 hard to decipher. Lemma 2.21(iii) in the same paper is also needed to express the term $\upsilon(\beta,\sigma,\mu)$ from Corollary~6.9 in terms of $\kappa(\sigma,\mu)$.

Note also that the first part of \cref{steinberg} follows from the second part of \cite[Corollary 6.9]{j10}; the condition $ep>n$ given there is regarded as automatically true when working over a field of characteristic zero.

As a consequence of \cref{steinberg} we can derive some information on the adjustment matrix~$A$.

\begin{cory}\label{adj}
Suppose $\la,\mu\in\calp(n)$ and all the columns of $\la$ are of even length, and write $\la=\dup\alpha$. Then
\[
A\v_{\la\mu}=
\begin{cases}
D_{\alpha\beta}\v&(\text{if $\mu=\dup\beta$ for $\beta\in\calp(n/2)$})\\
0&(\text{if $\mu$ has a column of odd length}),
\end{cases}
\]
and hence
\[
A_{\la\mu}=
\begin{cases}
D_{\alpha\beta}&(\text{if $\mu=\dup\beta$ for $\beta\in\calp(n/2)$})\\
0&(\text{if $\mu$ has a column of odd length}).
\end{cases}
\]
\end{cory}

\begin{pf}
The two parts of \cref{steinberg} combine to give
\begin{align*}
D_{\la\mu}\v&=\sum_{\beta\in\calp(n/2)}D_{\alpha\beta}\v\odc_{(\dup\beta)\mu}\v\\
\intertext{for all $\mu$. But the definition of the adjustment matrix also gives}
D_{\la\mu}\v&=\sum_{\nu\in\calp(n)}A\v_{\la\nu}\odc_{\nu\mu}\v
\end{align*}
for all $\mu$. Since the rows of $\odc\v$ are linearly independent, this enables us to deduce the given expression for $A\v_{\la\mu}$.

To get the result for $A_{\la\mu}$, consider the row vector $a$ with entries
\[
a_\mu=
\begin{cases}
D_{\alpha\beta}&(\text{if $\mu=\dup\beta$ for $\beta\in\calp(n/2)$})\\
0&(\text{if $\mu$ has a column of odd length}).
\end{cases}
\]
Using the result already proved for $A\v$, we have
\begin{align*}
(aA\v)_\mu&=\sum_{\nu\in\calp(n)} a_\nu A\v_{\nu\mu}\\
&=\sum_{\gamma\in\calp(n/2)}D_{\alpha\gamma}A\v_{(\dup\gamma)\mu}\\
&=
\begin{cases}
\sum_{\gamma\in\calp(n/2)}D_{\alpha\gamma}D\v_{\gamma\beta}&(\text{if $\mu=\dup\beta$ for $\beta\in\calp(n/2)$})\\
0&(\text{if $\mu$ has a column of odd length})
\end{cases}\\
&=\delta_{\la\mu}.
\end{align*}
Hence $a$ is the $\la$-row of the adjustment matrix $A$, as required.
% The result for $A_{\la\mu}$ then follows by inverting.
\end{pf}

Now we give a corresponding result for odd $n$; this is also a special case of \cite[Corollary~6.9]{j10}.

\begin{propn}\label{steinbergod}
Suppose $\la,\mu\in\calp(n)$ and all the columns of $\la$ are of even length except the first column. Write $\la=\dup\alpha\sqcup(1)$, and let $M$ be the set of partitions which can be obtained by removing one node from $\mu$. Then
\begin{align*}
\odc_{\la\mu}\v&=(-1)^{(n-1)/2}\sum_{\nu\in M}\epsilon(\nu)\kappa(\alpha,\nu),\\
\intertext{and}
D_{\la\mu}\v&=(-1)^{(n-1)/2}\sum_{\beta\in\calp((n-1)/2)}D_{\alpha\beta}\v\sum_{\nu\in M}\epsilon(\nu)\kappa(\beta,\nu).
\end{align*}
\end{propn}

Again, we deduce information about the adjustment matrix $A$; this is proved in exactly the same way as \cref{adj}.

\begin{cory}\label{adjod}
Suppose $\la,\mu\in\calp(n)$ and all the columns of $\la$ are of even length except the first, and write $\la=\dup\alpha\sqcup(1)$. Then
\[
A\v_{\la\mu}=
\begin{cases}
D_{\alpha\beta}\v&(\text{if $\mu=\dup\beta\sqcup(1)$ for $\beta\in\calp((n-1)/2)$})\\
0&(\text{if $\mu$ has a column (other than the first column) of odd length}),
\end{cases}
\]
and hence
\[
A_{\la\mu}=
\begin{cases}
D_{\alpha\beta}&(\text{if $\mu=\dup\beta\sqcup(1)$ for $\beta\in\calp((n-1)/2)$})\\
0&(\text{if $\mu$ has a column (other than the first column) of odd length}).
\end{cases}
\]
\end{cory}

\subsection{The degree of an irreducible spin character}\label{degsec}

Let $\deg(\chi)$ denote the degree of an irreducible character, i.e.\ the value $\chi(1)$. In this section we recall the ``bar-length formula'' which gives the degrees of the irreducible spin characters. This goes back to Schur \cite{schu}. Recall that for $\la\in\cald^-$ we write $\spn\la$ to mean either $\spn\la_+$ or $\spn\la_-$.

\begin{thm}[\textbf{(The bar-length formula)}]\label{barlength}
Suppose $\la\in\cald(n)$ has length $m$. Then
\[
\deg\spn\la=2^{\lfloor\frac12(n-m)\rfloor}\frac{n!}{\prod_{1\ls i\ls m}\la_i!}\prod_{1\ls i<j\ls m}\frac{\la_i-\la_j}{\la_i+\la_j}.
\]
\end{thm}

We shall use this formula extensively in \cref{degreesec}.

\subsection{Regularisation and doubling}\label{regnsec}

One of the most useful results in the modular representation theory of $\sss n$ is James's regularisation theorem, which gives an explicit composition factor (occurring with multiplicity $1$) in the $p$-modular reduction of $\ord\la$, for each $\la$. To state this result for $p=2$, we need to introduce the $2$-regularisation of a partition. For $l\gs0$, define the $l$th \emph{ladder} in $\bbn^2$ to be the set of nodes $(r,c)$ for which $r+c=l+2$. The intersection of this ladder with the Young diagram of a partition $\la$ is referred to as the $l$th ladder of $\la$. Now define the \emph{$2$-regularisation} $\la\reg$ of $\la$ to be the \trp whose Young diagram is obtained by moving the nodes of $\la$ as far up their ladders as possible. For example, if $\la=(4^2,3^4,1)$, then $\la\reg=(8,6,5,2)$, as we see from the following diagrams, in which we label nodes according to the ladders in which they lie.
\[
\young(0123,1234,234,345,456,567,6)\qquad\qquad\young(01234567,123456,23456,34)
\]
Now we can state a simple form of James's theorem.

\begin{thm}[\xcite{j1}{Theorem A}]\label{jreg}
Suppose $\la\in\calp$. Then $D_{\la(\la\reg)}=1$.
\end{thm}

Of course, this theorem is extremely useful in determining whether $\mord\la$ is irreducible, since it tells us that if $\mord\la$ is irreducible, then $\mord\la=\sid{\la\reg}$. For the main result in the present paper it will be useful to have an analogue of \cref{jreg} for spin characters. Fortunately there is such a result; this is due to Bessenrodt and Olsson, though a special case was proved earlier by Benson \cite[Theorem 1.2]{ben}. To state this, we need another definition. Given $\la\in\cald$, define the \emph{double} of $\la$ to be the partition
\[
\la\dbl=\big(\lceil\la_1/2\rceil,\lfloor\la_1/2\rfloor,\lceil\la_2/2\rceil,\lfloor\la_2/2\rfloor,\lceil\la_3/2\rceil,\lfloor\la_3/2\rfloor,\dots\big).
\]
In other words, $\la\dbl$ is obtained from $\la$ by replacing each part $\la_i$ with two parts which are equal (if $\la_i$ is even) or differ by $1$ (if $\la_i$ is odd). The fact that $\la$ is $2$-regular guarantees that $\la\dbl$ is a partition. We write $\la\dblreg$ to mean $(\la\dbl)\reg$.

Now we can give the ``spin-regularisation theorem''; recall that $\ev\la$ denotes the number of positive even parts of a partition $\la$.

\begin{thm}[\xcite{bo}{Theorem 5.2}]\label{spinreg}
Suppose $\la\in\cald$. Then $\decs\la{(\la\dblreg)}=2^{\lfloor\ev\la/2\rfloor}$.
\end{thm}

This result is extremely useful for our main problem: it tells us that $\mspn\la$ can only be irreducible if $\la$ has at most one non-zero even part, and that if $\mspn\la$ is irreducible, then $\mspn\la=\sdr\la$.

In order to exploit \cref{spinreg}, we will often want to show that $\la\dblreg=\mu\dblreg$ for certain $\la,\mu\in\cald(n)$ without actually calculating $\la\dblreg$ or $\mu\dblreg$. We do this using the combinatorics of \emph{slopes}. For $l\gs0$, define the $l$th slope in $\bbn^2$ to be the set of all nodes $(r,c)$ satisfying $2r+\inp{c/2}=l+2$. Define the $l$th slope of a partition $\la$ to be the intersection of this slope with the Young diagram of $\la$. (n.b.\ the slopes are essentially the ``ladders in the $\bar4$-residue diagram'' considered in \cite[\S3]{bo}, but we prefer not to over-tax the word ``ladder''. The reader should also note that unlike in \cite{bo} we do not consider shifted Young diagrams.) Now we have the following result, which is implicit in \cite{bo}.

\begin{lemma}\label{slopelad}
Suppose $\la$ is a \trp and $l\gs0$. Then the number of nodes in the $l$th slope of $\la$ equals the number of nodes in the $l$th ladder of $\la\dbl$. Hence if $\mu$ is another \trp, then $\la\dblreg=\mu\dblreg$ \iff $\la$ and $\mu$ have the same number of nodes in slope $l$ for each $l$.
\end{lemma}

For example, in the following diagram we label the nodes of $(14,8,7,1)$ with the numbers of the slopes containing them, and the nodes of $(14,8,7,1)\dbl=(7^2,4^3,3,1)$ with the numbers of the ladders containing them.
\[
\young(01122334455667,23344556,4556677,6)\qquad\young(0123456,1234567,2345,3456,4567,567,6)
\]

\begin{pf}
At the level of Young diagrams, replacing $\la$ with $\la\dbl$ involves replacing the nodes $(i,1),\dots,(i,\la_i)$ with
\[
(2i-1,1),\dots,(2i-1,\lceil\la_i/2\rceil),\quad(2i,1),\dots,(2i,\inp{\la_i/2})
\]
for each $i$. The former nodes lie in slopes
\[
2i-2,2i-1,2i-1,2i,2i,\dots,2i-3+\lceil\la_i/2\rceil,2i-2+\inp{\la_i/2}
\]
respectively, while the latter nodes lie in ladders
\[
2i-2,2i-1,\dots,2i-3+\lceil\la_i/2\rceil,\quad2i-1,2i,\dots,2i-2+\inp{\la_i/2}
\]
respectively. The result follows.
\end{pf}

\subsection{Blocks}

We now recall the $2$-block classification for $\tsss n$, which will be very useful in this paper. We begin by looking at the $2$-blocks of $\sss n$. The following result is a special case of the Brauer--Robinson Theorem \cite{brauer,gdbr}, first conjectured by Nakayama.

\begin{thm}\label{brarob}
Suppose $\la,\mu\in\calp(n)$. Then $\ord\la$ and $\ord\mu$ lie in the same $2$-block of $\sss n$ \iff $\la$ and $\mu$ have the same $2$-core.
\end{thm}

The same statement applies for $2$-blocks of $\tsss n$, and a corresponding statement \cite[Theorem 4.29]{jm1} holds for blocks of $\hhh n$. Since we only consider characteristic $2$ in this paper, we will henceforth say ``block'' to mean ``$2$-block''. Given \cref{brarob}, we may speak of the \emph{core} of a $2$-block $B$ of $\sss n$, meaning the common $2$-core of the partitions labelling the linear irreducible characters in $B$. These partitions necessarily have the same $2$-weight as well, and we call this the \emph{weight} of the block. We can immediately deduce the distribution of irreducible Brauer characters into blocks: if $\mu\in\cald(n)$, then since $\sid\mu$ occurs as a composition factor of $\mord\mu$, it lies in the same block as $\ord\mu$.

The block classification may alternatively be expressed in terms of residues. Recall that the \emph{residue} of a node $(r,c)$ is the residue of $c-r$ modulo $2$. Define the \emph{$2$-content} of a partition $\la$ to be the multiset of $0$s and $1$s comprising the residues of all the nodes of $\la$. We write a multiset of $0$s and $1$s in the form $\muset ab$. For example, the $2$-content of the partition $(7,4,1^3)$ is $\muset86$, as we see from the following ``$2$-residue diagram''.
\[
\yngres(2,7,4,1^3)
\]
It is an easy combinatorial exercise to show that two partitions have the same $2$-content \iff they have the same $2$-core and $2$-weight, so the block classification may alternatively be stated by saying that $\ord\la$ and $\ord\mu$ lie in the same block \iff $\la$ and $\mu$ have the same $2$-content. Accordingly, we can define the \emph{content} of a block $B$ to be the common $2$-content of the partitions labelling the linear irreducible characters in $B$.

Now we consider spin characters. Since the only irreducible $2$-modular characters are the characters $\sid\la$, the spin characters fit into the $2$-blocks described above, i.e.\ those containing linear characters. The distribution of the spin characters among these blocks is given by the following theorem, which was originally conjectured by Kn\"orr and Olsson \cite[p.~246]{ol}.

\begin{thm}[\xcite{bo}{Theorem 4.1}]\label{spinblockclass}
Suppose $\la\in\cald(n)$, and let $\sigma$ be the $2$-core of $\la\dbl$. Then $\spn\la$ lies in the $2$-block of $\tsss n$ with $2$-core $\sigma$.
\end{thm}

It will be helpful to have alternative descriptions of this block classification, for which we need different notions of core and content. If $\la\in\cald$, we define the \emph{\fbc} of $\la$ to be the \trp obtained by repeatedly applying the following operations to $\la$:
\begin{itemize}
\item
removing all even parts;
\item
removing any two parts whose sum is a multiple of $4$;
\item
replacing any odd part $\la_i\gs5$ with $\la_i-4$, if $\la_i-4$ is not already a part of $\la$.
\end{itemize}
The \fbc of $\la$ is easily seen to be well defined, and equals either $(4l-1,4l-5,\dots,3)$ or $(4l-3,4l-7,\dots,1)$ for some $l\gs0$. Note that these are precisely the partitions whose double is a $2$-core. Moreover, the double of the \fbc of $\la\in\cald$ coincides with the $2$-core of $\la\dbl$ \cite[Lemma 3.6]{bo}. So \cref{spinblockclass} may alternatively be stated by saying that if $\tau$ is the \fbc of $\la\in\cald(n)$, then $\spn\la$ lies in the block of $\tsss n$ with $2$-core $\tau\dbl$. In particular, two spin characters $\spn\la$ and $\spn\mu$ lie in the same block \iff $\la$ and $\mu$ have the same \fbc. So we can define the \emph{\fbc} of a block $B$ to be the common \fbc of the \trps labelling spin characters in $B$.

The \fbc of $\la$ is a partition of $n-2w$ for some $w$, which we call the \emph{\fbw} of $\la$. By the comments above, the \fbw of $\la$ equals the weight of the block containing $\spn\la$.

\begin{eg}
Consider the block $B$ of $\tsss{10}$ with $2$-core $(3,2,1)$. The partitions of $10$ with $2$-core $(3,2,1)$ are $(7,2,1)$, $(5,4,1)$, $(5,2,1^3)$, $(3,2^3,1)$ and $(3,2,1^5)$. $(3,2,1)$ is the double of the \fbc $(5,1)$, so $(5,1)$ is the \fbc of $B$. The partitions in $\cald(10)$ with \fbc $(5,1)$ are $(9,1)$ and $(5,4,1)$. So the ordinary irreducible characters in $B$ are
\[
\ord{7,2,1},\ord{5,4,1},\ord{5,2,1^3},\ord{3,2^3,1},\ord{3,2,1^5},
\]
\[
\spn{9,1},\spn{5,4,1}_+,\spn{5,4,1}_-,
\]
and the irreducible $2$-modular characters in $B$ are $\sid{7,2,1}$ and $\sid{5,4,1}$.
\end{eg}

\subsection{Branching rules}\label{brnchsec}

If $n>0$, then $\tsss{n-1}$ is naturally embedded in $\tsss n$, and a lot of information can be obtained by inducing and restricting characters between these two groups. In the modular representation theory of the symmetric groups, more delicate information can be gleaned by using the so-called \emph{$i$-induction} and \emph{$i$-restriction} functors. We summarise the key points here, specialising to characteristic $2$ (where the results automatically extend to $\tsss n$). The definition of the $i$-induction and $i$-restriction functors goes back to Robinson \cite{gdbrbook}, though our main reference will be the survey by Brundan and Kleshchev \cite{bk}.

Given a character $\chi$ of $\tsss n$, we write $\chi\res_{\tsss{n-1}}$ for its restriction to $\tsss{n-1}$, and $\chi\ind^{\tsss{n+1}}$ for the corresponding induced character for $\tsss{n+1}$. Now suppose $\chi$ lies in a single block $B$, with content $\muset ab$. Then we write $\ee0\chi$ for the component of $\chi\res_{\tsss{n-1}}$ lying in the block with content $\muset{a-1}b$ if there is such a block, and set $\ee0\chi=0$ otherwise. Similarly, we write $\ee1\chi$ for the component of $\chi\res_{\tsss{n-1}}$ lying in the block with content $\muset a{b-1}$ if there is such a block, and set $\ee1\chi=0$ otherwise. We extend the functions $\ee0,\ee1$ linearly. These functions $\ee0,\ee1$ can be applied to either ordinary characters or $2$-modular Brauer characters. In any case, it turns out that for any character $\chi$ we have $\chi\res_{\tsss{n-1}}=\ee0\chi+\ee1\chi$.

The effect of these functions on ordinary irreducible characters is well understood, via the following ``branching theorems''. The first of these dates back to Young, but the second is much more recent, due to Dehuai and Wybourne \cite[\S8]{dw}. To state these, we introduce some notation. Suppose $\la$ and $\mu$ are partitions, and $i\in\{0,1\}$. We write $\mu\adr ir\la$ to mean that $\la$ is obtained from $\mu$ by adding $r$ addable $i$-nodes (omitting the superscript $r$ when it equals $1$). Similarly, we write $\mu\ads ir\la$ to mean that $\la$ is obtained from $\mu$ by adding $r$ $i$-\spams.

\begin{thm}[\textbf{(The branching rule)}]\label{classbrnch}
Suppose $\la\in\calp(n)$ and $i\in\{0,1\}$, and let
\[
\La=\lset{\mu\in\calp(n-1)}{\mu\adr i1\la}.
\]
Then
\[
\ee i\ord\la=\sum_{\mu\in\La}\ord\mu.
\]
\end{thm}

\begin{thm}[\textbf{(The spin branching rule)}]\label{spinbrnch}
Suppose $\la\in\cald(n)$ and $i\in\{0,1\}$, and let
\[
\La=\lset{\mu\in\cald(n-1)}{\mu\ads i1\la}.
\]
\begin{enumerate}
\item
If $\la\in\cald^+(n)$, then
\begin{align*}
\ee i\spn\la&=\sum_{\mu\in\La\cap\cald^+(n-1)}\spn\mu\ +\ \sum_{\mu\in\La\cap\cald^-(n-1)}(\spn\mu_++\spn\mu_-).
\\
\intertext{\item
If $\la\in\cald^-(n)$, then}
\ee i\spn\la_\pm&=\sum_{\mu\in\La\cap\cald^+(n-1)}\spn\mu\ +\ \sum_{\mu\in\La\cap\cald^-(n-1)}\spn\mu_\pm.
\end{align*}
\end{enumerate}
\end{thm}

To state further results, we need to consider powers: $\ee0$ and $\ee1$ are defined for any $n$, so we can define powers $\ee i^r$ and (if we allow $\bbq$-linear combinations of characters) divided powers $\eed ir=\ee i^r/r!$. For any non-zero character $\chi$ and $i\in\{0,1\}$, we can then define
\[
\epsilon_i\chi=\max\lset{r\gs0}{\ee i^r\chi\neq0},\qquad\emx i\chi=\eed i{\epsilon_i\chi}\chi.
\]

The following result comes immediately from the classical branching rule.

\begin{propn}\label{ordmaxbranch}
Suppose $\la\in\calp$ and $i\in\{0,1\}$, and let $\la^-$ be the partition obtained by removing all the removable $i$-nodes from $\la$. Then $\emx i\ord\la=\ord{\la^-}$.
\end{propn}

For the spin branching rule, a result of the same form holds, but it is more difficult to keep track of multiplicities. We begin with the following result, which will also be useful in later sections.

\begin{propn}\label{spinbranchpower}
Suppose $\la\in\cald(n)$, $\mu\in\cald(n-r)$, $i\in\{0,1\}$, and $\mu\ads ir\la$. Let
\[
c=\big|\lset{x\gs1}{\text{there are nodes of $\la\setminus\mu$ in columns $x$ and $x+1$}}\big|.
\]
\begin{itemize}
\item
If $\la\in\cald^+(n)$ and $\mu\in\cald^+(n-r)$, then
\begin{align*}
\chm{\eed ir\spn\la}{\spn\mu}&=2^{\lfloor r/2\rfloor-c}.
\\
\intertext{\item
If $\la\in\cald^+(n)$ and $\mu\in\cald^-(n-r)$, then}
\chm{\eed ir\spn\la}{\spn\mu_++\spn\mu_-}&=2^{\lfloor(r+1)/2\rfloor-c}.
\\
\intertext{\item
If $\la\in\cald^-(n)$ and $\mu\in\cald^+(n-r)$, then}
\chm{\eed ir\spn\la_\pm}{\spn\mu}&=2^{\lfloor(r-1)/2\rfloor-c}.
\\
\intertext{\item
If $\la\in\cald^-(n)$ and $\mu\in\cald^-(n-r)$, then}
\chm{\eed ir\spn\la_\pm}{\spn\mu_++\spn\mu_-}&=2^{\lfloor r/2\rfloor-c}.
\end{align*}
\end{itemize}
\end{propn}

\begin{pf}
Recall that we write $\ev\nu$ for the number of positive even parts of $\nu\in\cald(n)$, and that by definition $\nu\in\cald^+(n)$ \iff $\ev\nu$ is even. When we remove a \sprm from $\nu$, we change the parity of $\ev\nu$, unless the removed node is the only node in its row. Applying this observation $r$ times, we find the following, when $\mu\ads ir\la$:
\begin{itemize}
\item
if $\ev\la+\ev\mu+r$ is even, then none of the nodes of $\la\setminus\mu$ lies in the first column (i.e.\ $\la$ and $\mu$ have the same length);
\item
if $\ev\la+\ev\mu+r$ is odd, then one of the nodes of $\la\setminus\mu$ does lie in column $1$ (and so in particular $i=0$).
\end{itemize}

Now we use induction on $r$, with the case $r=0$ being trivial. For the inductive step, we consider only the case where $\la\in\cald^+(n)$, $\mu\in\cald^+(n-r)$ and $r$ is odd; the other cases are similar (and often simpler). Let $C$ be the set of column labels of the nodes of $\la\setminus\mu$.

By the above remarks, we must have $i=0$ in this case, and $1\in C$. Apply the spin branching rule, and consider the partitions $\nu\in\cald(n-1)$ such that $\mu\subseteq\nu\subset\la$. These come in three types.
\begin{enumerate}
\item
The partition $\nu$ obtained from $\la$ by removing the node at the bottom of column $1$. Then $\nu\in\cald^+(n-1)$, $\chm{\ee0\spn\la}{\spn\nu}=1$ by the branching rule, and by induction $\chm{\eed0{r-1}\spn\nu}{\spn\mu}=2^{(r-1)/2-c}$.
\item
There are $r-1-2c$ different partitions $\nu$ obtained by removing a node in column $x>1$, where neither $x-1$ nor $x+1$ lies in $C$. In these cases $\nu\in\cald^-(n-1)$, $\chm{\ee0\spn\la}{\spn\nu_+}=\chm{\ee0\spn\la}{\spn\nu_-}=1$ by the spin branching rule, and by induction $\chm{\eed0{r-1}\spn\nu_\pm}{\spn\mu}=2^{(r-3)/2-c}$.
\item
For each pair $(x,x+1)$ of consecutive integers in $C$, there is one \trp $\nu$ obtained by removing a node in column $x$ or $x+1$: if the nodes of $\la\setminus\mu$ lying in columns $x$ and $x+1$ are of the form $(m,x),(m,x+1)$, then $\nu$ is obtained by removing the node $(m,x+1)$, while if they have the form $(m,x),(m-1,x+1)$, then $\nu$ is obtained by removing $(m,x)$. In any case $\nu\in\cald^-(n-1)$, $\chm{\ee0\spn\la}{\spn\nu_+}=\chm{\ee0\spn\la}{\spn\nu_-}=1$ by the spin branching rule, and by induction $\chm{\eed0{r-1}\spn\nu_\pm}{\spn\mu}=2^{(r-3)/2-(c-1)}$.
\end{enumerate}
So in total we find that
\begin{align*}
\chm{\eed0{r-1}\ee0\spn\la}{\spn\mu}&=2^{(r-1)/2-c}+(r-2c-1)\times2\times2^{(r-3)/2-c}+c\times2\times2^{(r-3)/2-(c-1)}\\
&=r\times2^{(r-1)/2-c},
\end{align*}
which is what we need.
\end{pf}

\begin{eg}
We give an example which illustrates a different case of \cref{spinbranchpower}. Take $\la=(11,9,7,5,4,1)$ and $r=2$. The $0$-\sprms of $\la$ are highlighted in the following diagram.
\[
\young(01100110011,01100110!\thk0!\thn,0110011,0110!\thk0!\thn,011!\thk0,0)
\]
The spin branching rule gives
\begin{align*}
\ee0\spn\la_\pm=\ &\spn{11,9,7,5,4}_\pm+\spn{11,9,7,5,3,1)}+\spn{11,8,7,5,4,1)}.
\\
\intertext{The case $r=1$ of the \lcnamecref{spinbranchpower} gives}
\ee0\spn{11,9,7,5,4}_\pm&=\spn{11,9,7,5,3}+\spn{11,8,7,5,4},\\
\ee0\spn{11,9,7,5,3,1}&=\spn{11,9,7,5,3}+\spn{11,9,7,4,3,1}_++\spn{11,9,7,4,3,1}_-\\
&\phantom=\ +\spn{11,8,7,5,3,1}_++\spn{11,8,7,5,3,1}_-,\\
\ee0\spn{11,8,7,5,4,1}&=\spn{11,8,7,5,4}+\spn{11,8,7,5,3,1}_++\spn{11,8,7,5,3,1}_-,
\\
\intertext{and we obtain}
\eed02\spn\la_\pm&=\spn{11,9,7,5,3}+\spn{11,8,7,5,4}+\tfrac12(\spn{11,9,7,4,3,1}_++\spn{11,9,7,4,3,1}_-)\\
&\phantom=\ +\spn{11,8,7,5,3,1}_++\spn{11,8,7,5,3,1}_-
\end{align*}
as predicted by \cref{spinbranchpower}.
\end{eg}

Now we can deduce the following about $\emx i\spn\la$ for $\la\in\cald$; this will be central to the proof of our main theorem.

\begin{cory}\label{spinmaxbranch}
Suppose $\la\in\cald$ and $i\in\{0,1\}$. Suppose $\la$ has $r$ $i$-\sprms, and let $\la^-$ be the partition obtained by removing all these nodes. Then:
\begin{enumerate}
\item
if $\la^-\in\cald^+(n-r)$, then $\emx i\spn\la=a\spn{\la^-}$ for some $a\in\bbn$;
\item
if $\la^-\in\cald^-(n-r)$, then $\emx i\spn\la=a\spn{\la^-}_++b\spn{\la^-}_-$ with $a+b\in\bbn$.
\end{enumerate}
\end{cory}

Analogous results apply to induced modules. For a character $\chi$ of $\tsss n$, we have $\chi\ind^{\tsss{n+1}}=\ff0\chi+\ff1\chi$, where (if $\chi$ lies in the block with content $\muset ab$) $\ff0\chi$ is the component of $\chi\ind^{\tsss{n+1}}$ lying in the block with content $\muset{a+1}b$, and $\ff1\chi$ is defined similarly. For a non-zero character $\chi$ we define $\varphi_i\chi=\max\lset{r\gs0}{\ff i^r\chi\neq0}$, and we set $\fmx i\chi=\ffd i{\varphi_i\chi}\chi$. Then we have analogues of \cref{classbrnch,spinbrnch,ordmaxbranch,spinbranchpower,spinmaxbranch} (which follow from these results by Frobenius reciprocity), in which we add nodes rather than removing nodes. We will not state these results explicitly, but instead refer to the ``induction versions'' of the restriction results above.

Next we recall some of Kleshchev's ``modular branching rules'', to describe what happens to the irreducible Brauer characters $\sid\mu$ under the $i$-restriction operators. This involves the combinatorics of \emph{normal nodes}.

Suppose $\mu\in\cald$ and $i\in\{0,1\}$, and construct a sequence of $+$ and $-$ signs by reading along the edge of the Young diagram of $\mu$ from top to bottom and writing a $+$ for each addable $i$-node and a $-$ for each removable $i$-node. This sequence is called the \emph{$i$-signature} of $\mu$. Now successively delete all adjacent pairs $+-$ from this sequence until none remain. The resulting sequence is called the \emph{reduced $i$-signature} of $\mu$. The removable nodes corresponding to the $-$ signs in the reduced $i$-signature are called the \emph{normal $i$-nodes} of $\mu$, and the addable nodes corresponding to the $+$ signs are the \emph{conormal $i$-nodes} of $\mu$. Now we have the following result. See \cite{bk} (in particular, the discussion following Lemma 2.12) for this and more general modular branching results.

\begin{thm}\label{modbranch}
Suppose $\mu\in\cald$ and $i\in\{0,1\}$. Let $\mu^-$ be the \trp obtained by removing all the normal $i$-nodes from $\mu$, and let $\mu^+$ be the partition obtained by adding all the conormal $i$-nodes to $\mu$. Then $\mu^-,\mu^+$ are $2$-regular, and
\[
\emx i\sid\mu=\sid{\mu^-},\qquad\fmx i\sid\mu=\sid{\mu^+}.
\]
In particular, if $\chi$ is an irreducible Brauer character of $\tsss n$, then $\emx i\chi$ and $\fmx i\chi$ are irreducible Brauer characters.
\end{thm}

The last statement in this theorem is extremely useful for us, and will form the basis of our induction proof of our main theorem.

\begin{eg}
Let $\mu=(15,11,8,6,5,2)$ and $i=0$. The Young diagram of $\mu$ with the residues of addable and removable nodes indicated, is as follows.
\Yinternals0\Yaddables1
\[
\yngres(2,15,11,8,6,5,2)
\]
So the $0$-signature of $\mu$ is $-++---+$, and hence the reduced $0$-signature is $--+$, with the normal nodes being $(1,15)$ and $(6,2)$, and the conormal node $(7,1)$. Hence
\[
\emx0\sid\mu=\sid{14,11,8,6,5,1},\qquad\fmx0\sid\mu=\sid{15,11,8,6,5,2,1}.
\]
\end{eg}

\subsection{Symmetric functions}\label{symfnsec}

The combinatorics of partitions involved in representation theory is closely connected with the combinatorics involved in symmetric functions, and in this paper we shall exploit this connection. In this section we set out the notation and background results we shall need. We give just the bare essentials, since symmetric function theory is covered in detail elsewhere.

We consider the space $\sym$ of symmetric functions over $\bbc$ in infinitely many variables. For $\la\in\calp$, let $s_\la$ denote the corresponding \emph{Schur function}. Let $\ip\,\,$ be the inner product on $\sym$ for which the Schur functions are orthonormal.

For each $r$, we let $h_r$ denote the complete homogeneous symmetric function of degree $r$ (which coincides with the Schur function $s_{(r)}$). For any partition $\mu$, we write $h_\mu$ for the product $h_{\mu_1}h_{\mu_2}\dots$. The \emph{Pieri rule} says that for any $r$ and for $\la\in\calp$ we have $h_rs_\la=\sum s_\nu$, summing over all partitions $\nu$ which can be obtained from $\la$ by adding $r$ nodes in distinct columns. Inductively, this enables us to express any $h_\mu$ as a sum of Schur functions; the coefficients obtained are the \emph{Kostka numbers}. Two immediate consequences of this are that $\ip{h_\mu}{s_\mu}=1$ for each $\mu$, and that $\ip{h_\mu}{s_\la}$ is non-zero only if $\la\dom\mu$.

We write $e_r$ for the $r$th elementary symmetric function (which coincides with $s_{(1^r)}$), and for any $\mu\in\calp$ we write $e_\mu=e_{\mu_1}e_{\mu_2}\dots$. The \emph{dual Pieri rule} says that for any $r$ and any $\la$ we have $e_rs_\la=\sum_\nu s_\nu$, summing over all partitions $\nu$ which can be obtained from $\la$ by adding $r$ nodes in distinct rows. Inductively, this enables us to express any $e_\mu$ as a sum of Schur functions (again, in terms of Kostka numbers). This has the consequences that $\ip{e_\mu}{s_{\mu'}}=1$ for any $\mu$ and that $\ip{e_\mu}{s_\la}$ is non-zero only if $\mu'\dom\la$. The combination of the Pieri and dual Pieri rules yields the identity $\ip{e_\mu}{s_\la}=\ip{h_\mu}{s_{\la'}}$.

The sets $\lset{s_\la}{\la\in\calp}$, $\lset{h_\la}{\la\in\calp}$ and $\lset{e_\la}{\la\in\calp}$ are all bases for $\sym$. We shall need to consider the transition coefficients which allow us to express a complete homogeneous function $h_\la$ in terms of the elementary functions $e_\mu$. So define $\ke\la\mu$ for all $\la,\mu\in\calp$ by $h_\la=\sum_\mu\ke\la\mu e_\mu$. In the case where $\la=(k)$, these coefficients are given by a special case of the \emph{second Jacobi--Trudi formula}, which says that $h_k$ equals the determinant of the $k\times k$ matrix
\[
\left\lvert\begin{matrix}e_1&e_2&e_3&\cdots&e_k\\1&e_1&e_2&\ddots&\vdots\\0&1&e_1&\ddots&e_3\\\vdots&\ddots&\ddots&\ddots&e_2\\0&\cdots&0&1&e_1\end{matrix}\right\rvert.
\]

We now derive some consequences of this for partitions with only even parts.

\begin{lemma}\label{ekdouble}
If $\mu\in\calp(k)$, then $\ke{(2k)}{2\mu}=(-1)^k\ke{(k)}\mu$ for any partition $\mu$.
\end{lemma}

\begin{pf}
Using the second Jacobi--Trudi formula and expanding the determinant in terms of permutations, we find that $\ke{(k)}\mu$ is sum of the signs of all the permutations $\pi\in\sss k$ with the property that $\pi(i)\gs i-1$ for each $i$ and the integers $\pi(1),\pi(2)-1,\dots,\pi(k)-k+1$ equal $\mu_1,\dots,\mu_k$ in some order. Let $\Pi_\mu$ be the set of such permutations. Then it suffices to show that there is a bijection from $\Pi_\mu$ to $\Pi_{2\mu}$ which preserves the sign of every permutation if $k$ is even, or changes it if $k$ is odd.

To construct the required map from $\Pi_\mu$ to $\Pi_{2\mu}$, we take $\pi\in\Pi_\mu$, and define $\hat\pi\in\Pi_{2\mu}$ by
\[
\hat\pi(i)=
\begin{cases}
2\pi((i+1)/2)&(i\text{ odd})\\
i-1&(i\text{ even}).
\end{cases}
\]
To construct the inverse map, we take $\sigma\in\Pi_{2\mu}$, and observe that because $\sigma(i)-i$ is odd and $\sigma(i)\gs i-1$ for every $i$, we must have $\sigma(i)=i-1$ for all even $i$. Now we can define $\check\sigma\in\Pi_\mu$ by $\check\sigma(i)=\sigma(2i-1)/2$ for each $i$. It is clear that these two maps are mutually inverse, so it remains to show that $\sgn{\hat\pi}=(-1)^k\sgn\pi$ for every $\pi\in\Pi_\mu$.

Regarding $\pi$ as an element of $\sss{2k}$ by the usual embedding of $\sss k$ in $\sss{2k}$, we can write $\hat\pi=\kappa\circ\pi\circ\iota$, where $\iota,\kappa\in\sss{2k}$ are given by
\[
\iota(i)=
\begin{cases}
(i+1)/2&(i\text{ odd})\\
i/2+k&(i\text{ even}),
\end{cases}
\qquad
\kappa(i)=
\begin{cases}
2i&(i\ls k)\\
2(i-k)-1&(i>k).
\end{cases}
\]
It is easy to see that the signs of $\iota$ and $\kappa$ are $(-1)^{\binom k2}$ and $(-1)^{\binom{k+1}2}$ respectively, which gives the result.
\end{pf}

By writing $h_\la$ as the product $h_{\la_1}h_{\la_2}\dots$ and applying \cref{ekdouble}, we obtain the following result. Note that for partitions $\nu^1,\nu^2,\dots$ we have $e_{\nu^1}e_{\nu^2}\dots=e_{2\mu}$ \iff each $\nu^i$ has only even parts and $e_{\frac12\nu^1}e_{\frac12\nu^2}\dots=e_\mu$, i.e.\ $\frac12\nu^1\sqcup\frac12\nu^2\sqcup\dots=\mu$.

\begin{propn}\label{ekdoublep}
Suppose $\la,\mu\in\calp(k)$. Then $\ke{2\la}{2\mu}=(-1)^k\ke\la\mu$.
\end{propn}

Finally we note the following simple result.

\begin{lemma}\label{simpleep}
Suppose $\mu\in\calp(k)$ and $\la\in\calp(2k)$ with $\la_i$ odd for at least one value of $i$. Then $\ke\la{2\mu}=0$.
\end{lemma}

\begin{pf}
$h_\la$ includes a factor $h_{\la_i}$; when this is expressed in terms of elementary symmetric polynomials $e_\nu$, each $\nu$ satisfies $|\nu|=\la_i$, and so must have at least one odd part. Hence each $e_\nu$ appearing when $h_\la$ is expressed in terms of elementary symmetric polynomials must contain a factor $e_l$ for some odd $l$.
\end{pf}

\section{$2$-Carter partitions and the main theorem}\label{mainthmsec}

Having set out the background results we need, we are now able to state our main theorem. This relies on another combinatorial definition: say that $\la\in\calp$ is \emph{$2$-Carter} if for every $r\gs1$, $\la_r-\la_{r+1}+1$ is divisible by a power of $2$ greater than~$\la_{r+1}-\la_{r+2}$.

\begin{eg}
The $2$-Carter partitions $\la$ with $\la_1\ls5$ are
\[
\varnothing,(1),(2),(2,1),(3),(3,2,1),(4),(4,1),(4,3,2,1),(5),(5,2),(5,2,1),(5,4,3,2,1).
\]
\end{eg}

The notion of a $2$-Carter partition was introduced for the classification of irreducible Specht and Weyl modules. We shall need the following result later (recall that $\schur n$ is the Schur algebra over a field $\bbf$ of characteristic $2$).

\begin{thm}[\xcite{jm1}{Theorem 4.5}]\label{irredweyl}
Suppose $\la\in\calp(n)$. Then the Weyl module $\weyl\la$ for $\schur n$ is irreducible \iff $\la$ is a $2$-Carter partition.
\end{thm}

(We remark that \cite[Theorem 4.5]{jm1} is phrased in terms of hook lengths in the Young diagram of $\la$. The fact that this formulation is equivalent to the one we have given is explained by James in \cite[Lemma~3.14]{j0}.)

Now we come to the classification of ordinary irreducible characters for $\tsss n$ that remain irreducible in characteristic $2$. For the characters $\ord\la$ of the Specht modules, the answer is the same as for $\sss n$, and is given by the following theorem.

\begin{thm}[\xcite{jmp2}{Main Theorem}]\label{jmp2main}
Suppose $\la\in\calp$. Then $\modr{\ord\la}$ is irreducible \iff either $\la$ or $\la'$ is a $2$-Carter partition, or $\la=(2^2)$.
\end{thm}

So it remains to consider spin characters. Recall that when $\la\in\cald^-(n)$ we write $\spn\la$ to mean either $\spn\la_+$ or $\spn\la_-$, and $\mspn\la$ is then unambiguously defined, since $\modr{\spn\la_+}=\modr{\spn\la_-}$. So we can phrase the main question for spin characters simply as ``for which $\la\in\cald$ is $\mspn\la$ irreducible?''. Now we can state our main theorem.

\begin{thm}\label{main}
Suppose $\la\in\cald$. Then $\mspn\la$ is irreducible \iff one of the following occurs.
\begin{enumerate}
\item
$\la=\tau+4\alpha$, where $\tau=(4l-1,4l-5,\dots,3)$ for some $l\gs0$ and $\alpha$ is a $2$-Carter partition with $\len\alpha\ls l$.
\item
$\la=\tau+4\alpha$, where $\tau=(4l-3,4l-7,\dots,1)$ for some $l\gs1$ and $\alpha$ is a $2$-Carter partition with $\len\alpha\ls l$.
\item
$\la=\tau+4\alpha\sqcup(2)$, where $\tau=(4l-1,4l-5,\dots,3)$ for some $l\gs0$ and $\alpha$ is a $2$-Carter partition with $\len\alpha\ls l$.
\item
$\la=\tau+4\alpha\sqcup(2)$, where $\tau=(4l-3,4l-7,\dots,1)$ for some $l\gs1$ and $\alpha$ is a $2$-Carter partition with $\len\alpha\ls l-1$.
\item
$\la$ equals $(2b)$ or $(4b-2,1)$ for some $b\gs2$.
\item
$\la=(3,2,1)$.
\end{enumerate}
\end{thm}

\cref{degreesec,rouqandsep,mainproofsec} are devoted to the proof of this theorem.

\section{Degrees of irreducible characters}\label{degreesec}

In this section we use the bar-length formula to compare the degrees of certain characters. The motivating observation here is the following.

\begin{lemma}\label{samerdoub}
Suppose $\la,\mu\in\cald(n)$, with $\la\dblreg=\mu\dblreg$ and $\deg\spn\la>\deg\spn\mu$. Then $\mspn\la$ is reducible.
\end{lemma}

\begin{pf}
By \cref{spinreg}, $\sdr\la$ occurs as a constituent of $\mspn\la$, so if $\mspn\la$ were irreducible, we would have $\mspn\la=\sdr\la$, and in particular $\deg(\sdr\la)=\deg\spn\la$. But $\sdr\la=\sdr\mu$ also occurs as a constituent of $\mspn\mu$, so $\deg(\sdr\la)\ls\deg\spn\mu<\deg\spn\la$.
\end{pf}

We now construct various families of pairs of partitions $\la,\mu$ satisfying the hypotheses of \cref{samerdoub}. Our first result is as follows.

\begin{lemma}\label{dimen1}
Suppose $m\gs2$, and let
\begin{align*}
\la^m&=(4m,4m-3,4m-7,\dots,9,5),\\
\mu^m&=(4m+1,4m-3,4m-7,\dots,9,4).
\end{align*}
Then $(\la^m)\dblreg=(\mu^m)\dblreg$, and $\deg\spn{\la^m}>\deg\spn{\mu^m}$.
\end{lemma}

\begin{pf}
$\mu^m$ is obtained from $\la^m$ by replacing the node $(m,5)$ with the node $(1,4m+1)$. Since these nodes both lie in the same slope, we have $(\la^m)\dblreg=(\mu^m)\dblreg$ by \cref{slopelad}.

Now we consider the degrees. We use induction on $m$, with the case $m=2$ an easy check. For the inductive step, it suffices to prove that when $m\gs3$
\[
\frac{\deg\spn{\la^m}\deg\spn{\mu^{m-1}}}{\deg\spn{\mu^m}\deg\spn{\la^{m-1}}}>1.
\]
Directly from the bar-length formula, this ratio equals
\[
\frac{(4m-5)(4m-1)(4m+1)^2(4m+5)}{(2m+1)(2m+3)(4m-3)(8m-7)(8m-3)}.
\]
To show that this is always greater than $1$, we show that the numerator exceeds the denominator for all $m\gs2$. The difference between the numerator and the denominator is $f(m)=256m^4+368m^3-1220m^2-368m+214$. We have $f(2)$ and $f(0)$ positive, while $f(1)$ and $f(-1)$ are negative, so that all the roots of $f$ are less than $2$.
\end{pf}

Now we give a similar dimension argument which will help when we consider Rouquier blocks.

\begin{propn}\label{firstrdim}
Given $a>0$ and $m\gs0$, define
\begin{align*}
\la^{a,m}&=(4a+4m-3,4a+4m-7,\dots,1)\sqcup(4a),\\
\mu^{a,m}&=(4a+4m+1,4a+4m-3,\dots,4m+5,4m-3,4m-7,\dots,1).
\end{align*}
Then $(\la^{a,m})\dblreg=(\mu^{a,m})\dblreg$, and $\deg\spn{\la^{a,m}}>\deg\spn{\mu^{a,m}}$.
\end{propn}

\begin{pf}
To see that $(\la^{a,m})\dblreg=(\mu^{a,m})\dblreg$ we use \cref{slopelad}, showing that $\mu^{a,m}$ can be obtained from $\la^{a,m}$ by moving some nodes, but keeping each node in the same slope. There are two cases.
\begin{enumerate}
\item
Suppose $a>m$. Let $A$ be the set of nodes of $\la^{a,m}$ comprising
\begin{itemize}
\item
the last four nodes in each of rows $a+1,\dots,a+m$, and
\item
the unique node in row $a+m+1$.
\end{itemize}
Now observe that $\mu^{a,m}$ can be obtained from $\la^{a,m}$ by replacing each node $(r,c)\in A$ with $(r-a,c+4a)$. Since $(r,c)$ and $(r-a,c+4a)$ lie in the same slope, the result follows.
\item
Now suppose $a\ls m$. Let $A$ be the set of nodes of $\la^{a,m}$ comprising
\begin{itemize}
\item
the last three nodes in row $m+1$, and
\item
the last four nodes in each of rows $m+2,\dots,a+m$.
\end{itemize}
Now $\mu^{a,m}$ can be obtained from $\la^{a,m}$ by replacing each node $(r,c)\in A$ with $(r-m,c+4m)$, and also replacing the node $(a+m+1,1)$ with $(1,4a+4m+1)$. Again, each moved node remains in the same slope, so the result follows.
\end{enumerate}

Now we consider degrees. Direct from the bar-length formula we get
\begin{align*}
\frac{\deg\spn{\la^{a,m}}}{\deg\spn{\mu^{a,m}}}=&\frac{(4a+4m+1)!\big((8a+8m-2)(8a+8m-6)\dots(4a+8m+6)\big)\big(4m(4m-4)\dots4\big)}
{(4a)!(4m)!\big((4a+4m-2)(4a+4m-6)\dots(4m+2)\big)}\\
&\times\frac{\big((4m-3)(4m-7)\dots1\big)\big((4a-1)(4a-5)\dots3\big)}{\big((8a+4m-3)(8a+4m-7)\dots(4a+1)\big)\big((4a+4m)(4a+4m-4)\dots(4a+4)\big)}.\tag*{(\textasteriskcentered)}
\end{align*}
Now we proceed by induction on $a$. In the case $a=1$, (\textasteriskcentered) becomes $(4m+3)/2$, which is greater than $1$. For the inductive step, it suffices to show that
\[
\frac{\deg\spn{\la^{a+1,m}}}{\deg\spn{\mu^{a+1,m}}}\frac{\deg\spn{\mu^{a,m}}}{\deg\spn{\la^{a,m}}}>1.
\]
From (\textasteriskcentered), this ratio equals
\[
\frac{(4a+4m+5)(4a+4m+3)^2(4a+4m+1)}{(8a+4m+5)(8a+4m+1)(2a+4m+3)(2a+1)}.
\]
The difference between the numerator and denominator in this fraction is
\begin{align*}
&256m^4+896am^3+960a^2m^2+256a^3m+704m^3+1728am^2+1056a^2m+64a^3\\
&+656m^2+984am+204a^2+244m+152a+30
\end{align*}
which is obviously positive.
\end{pf}

The next three results are proved in exactly the same way.

\begin{propn}\label{secondrdim}
Given $a>0$ and $m\gs0$, define
\begin{align*}
\la^{a,m}&=(4a+4m+1,4a+4m-3,\dots,1)\sqcup(4a+2),\\
\mu^{a,m}&=(4a+4m+5,4a+4m+1,\dots,4m+9,4m+1,4m-3,\dots,5,2,1).
\end{align*}
Then $(\la^{a,m})\dblreg=(\mu^{a,m})\dblreg$, and $\deg\spn{\la^{a,m}}>\deg\spn{\mu^{a,m}}$.
\end{propn}

\begin{propn}\label{thirdrdim}
Given $a>0$ and $m\gs0$, define
\begin{align*}
\la^{a,m}&=(4a+4m-1,4a+4m-5,\dots,3)\sqcup(4a),\\
\mu^{a,m}&=(4a+4m+3,4a+4m-1,\dots,4m+7,4m-1,4m-5,\dots,3).
\end{align*}
Then $(\la^{a,m})\dblreg=(\mu^{a,m})\dblreg$, and $\deg\spn{\la^{a,m}}>\deg\spn{\mu^{a,m}}$.
\end{propn}

\begin{propn}\label{fourthrdim}
Given $a>0$ and $m\gs0$, define
\begin{align*}
\la^{a,m}&=(4a+4m-1,4a+4m-5,\dots,3)\sqcup(4a+2),\\
\mu^{a,m}&=(4a+4m+3,4a+4m-1,\dots,4m+7,4m-1,4m-5,\dots,3,2).
\end{align*}
Then $(\la^{a,m})\dblreg=(\mu^{a,m})\dblreg$, and $\deg\spn{\la^{a,m}}>\deg\spn{\mu^{a,m}}$.
\end{propn}

To make these results more general, we now consider row removal. We start with the following lemma.

\begin{lemma}\label{domdim}
Suppose $\la,\mu\in\calp(n)$ with $\mu\doms\la$, and that $l$ is an integer greater than $\mu_1$. Then
\[
\prod_{i\gs1}\frac{l+\mu_i}{l-\mu_i}>\prod_{i\gs1}\frac{l+\la_i}{l-\la_i}.
\]
\end{lemma}

\begin{pf}
We may assume that $\mu$ covers $\la$ in the dominance order, in which case $\la$ is obtained from $\mu$ by moving a single node down to a lower row. So suppose that for some $j<k$ we have $\mu_j=a+1$, $\la_j=a$, $\mu_k=b$, $\la_k=b+1$, and that $\la_i=\mu_i$ for all $i\neq j,k$. Then the ratio of the left-hand expression to the right-hand expression is
\[
\frac{(l-a)(l-b-1)(l+a+1)(l+b)}{(l+a)(l+b+1)(l-a-1)(l-b)}.
\]
The difference between the numerator and denominator here is $2l(a-b)(a+b+1)$, which is positive, so the ratio is greater than $1$.
\end{pf}

\begin{lemma}\label{dimrowrem}
Suppose $\la,\mu\in\cald(n)$ with $\mu\doms\la$, and $l$ is an odd integer with $l>\mu_1$. Define
\begin{align*}
\la^+&=(l,\la_1,\la_2,\dots),\\
\mu^+&=(l,\mu_1,\mu_2,\dots).
\end{align*}
Then
\[
\frac{\deg\spn{\la^+}}{\deg\spn{\mu^+}}>\frac{\deg\spn\la}{\deg\spn\mu}.
\]
Furthermore, if $\la\dblreg=\mu\dblreg$, then $(\la^+)\dblreg=(\mu^+)\dblreg$.
\end{lemma}

\begin{pf}
From the bar-length formula, the ratio of the left-hand side to the right-hand side is
\[
\prod_{i\gs1}\frac{l+\mu_i}{l-\mu_i}\frac{l-\la_i}{l+\la_i},
\]
and by \cref{domdim} this is greater than $1$. The second statement follows by considering the slopes containing the nodes of $\la$ and the corresponding nodes of $\la^+$.
\end{pf}

\section{Rouquier blocks and \qs partitions}\label{rouqandsep}

Rouquier blocks are a certain class of particularly well-behaved blocks of symmetric groups (and more generally of Iwahori--Hecke algebras and $q$-Schur algebras). In this section we summarise some of the important properties of Rouquier blocks of symmetric groups and their double covers in characteristic $2$, and then examine the decomposition numbers for Rouquier blocks of $\tsss n$. We compute the rows of the spin part of the decomposition matrix labelled by \trps with only odd parts, showing that some of the corresponding spin characters are irreducible in characteristic $2$. We extend these results to what we call ``\qs'' \trps, and prove our main theorem for spin characters labelled by \qs partitions.

\subsection{Rouquier blocks}

Suppose $B$ is a block of $\tsss n$, with $2$-core $\sigma=(c,c-1,\dots,1)$ and weight $w$. We say that $B$ is \emph{Rouquier} if $w\ls c+1$.

The first thing that makes Rouquier blocks easy to understand is a simple description of the \trps labelling characters in a Rouquier block.

\begin{lemma}\label{trpsinrouq}
Suppose $\sigma=(c,c-1,\dots,1)$ and that $\la$ is a partition with $2$-core $\sigma$ and $2$-weight $w\ls c+1$. \Tfae.
\begin{enumerate}
\item
$\la$ is $2$-regular.
\item
The length of $\la$ is at most $c+1$.
\item
$\la$ has the form $\sigma+2\alpha$, where $\alpha\in\calp(w)$.
\end{enumerate}
\end{lemma}

\needspace{3em}
\begin{pf}
\indent
\begin{description}
\vspace{-\topsep}
\item[{\rm(1$\Rightarrow$2)}]
Suppose $\la$ is $2$-regular, and suppose for a contradiction that the length of $\la$ is greater than $c+1$. Then $(c+2,1)$ is a node of $\la$, so (since $\la$ is $2$-regular) $(c+1,2),(c,3),(c-1,4),\dots,(1,c+2)$ are all nodes of $\la$. But then $|\la|\gs\frac12(c+2)(c+3)>|\sigma|+2w$, a contradiction.
\item[{\rm(2$\Rightarrow$3)}]
We use induction on $w$, with the case $w=0$ being trivial. Assuming $w>0$, let $r\ls c+1$ be maximal such that $\la_r\gs c+2-r$ (there must be such an $r$, since $\la\supset\sigma$). Then we claim that $\la_r\gs\la_{r+1}+2$. If $r\ls c$ this is immediate from the choice of $r$, so suppose $r=c+1$. The node $(c+1,1)$ is a node of $\la$ but not of $\sigma$, so must be part of one of the rim $2$-hooks added to obtain $\la$ from $\sigma$. The other node in this rim $2$-hook must be $(c+1,2)$, since $(c,1)$ is a node of $\sigma$, and (by assumption) $(c+2,1)$ is not a node of $\la$. Hence $\la_{c+1}\gs2=\la_{c+2}+2$ as claimed.

Hence we can remove the rim $2$-hook $\{(r,\la_r-1),(r,\la_r)\}$ from $\la$ and leave a partition $\mu$ which satisfies the hypotheses of part (2). By induction $\mu=\sigma+2\beta$ for some $\beta\in\calp(w-1)$, and hence $\la=\sigma+2\alpha$, where $\alpha$ is the partition obtained from $\beta$ by adding a node at the end of row $r$.
%We use induction on $w$, with the case $w=0$ being trivial. Assuming $w>0$, let $r\ls c+1$ be maximal such that $\la_r\gs c+2-r$ (there must be such an $r$, since $\la\supset\sigma$). Then we claim that $\la_r\gs c+3-r$: note that $(r,c+2-r)$ is a node of $\la$ but not of $\sigma$, so must be part of one of the rim $2$-hooks added to obtain $\la$ from $\sigma$. The other node in this rim $2$-hook must be $(r,c+3-r)$, since $(r-1,c+2-r)$ and $(r,c+1-r)$ are nodes of $\sigma$, and (by the choice of $r$) $(r+1,c+2-r)$ is not a node of $\la$. Hence $\la_r\gs c+3-r$. By the choice of $r$ either $\la_{r+1}\ls c-r$ or $\la_{r+1}=0$, which means that we can remove the rim $2$-hook $\{(r,\la_r-1),(r,\la_r)\}$ from $\la$ and leave a partition $\mu$ which is still $2$-regular. By induction $\mu=\sigma+2\beta$ for some $\beta\in\calp(w-1)$, and hence $\la=\sigma+2\alpha$, where $\alpha$ is the partition obtained from $\beta$ by adding a node at the end of row $r$.
\item[{\rm(3$\Rightarrow$1)}]
Suppose $\la=\sigma+2\alpha$ with $\alpha\in\calp(w)$. Then the length of $\alpha$ is at most $w\ls c+1$. So $\la_r=0$ for $r>c+1$, and $\la_r-\la_{r+1}=1+2(\alpha_r-\alpha_{r+1})>0$ for $r\ls c$, and hence $\la$ is $2$-regular.\qedhere
\end{description}
\end{pf}

We can also describe the \trps labelling spin characters in Rouquier blocks. This requires some preliminary work.

\begin{lemma}\label{rouqlem3}
Suppose $\tau=(4l-1,4l-5,\dots,3)$ for $l\gs0$, and $\la$ is a \trp with \fbc $\tau$ and \fbw $w\ls2l+1$. Then $\la$ has no parts congruent to $1$ modulo $4$, and $\la$ includes all the integers $3,7,11,\dots,4t-1$, where $t=\lceil l-\frac12w\rceil$.
\end{lemma}

\begin{pf}
We use induction on $w$. In the case $w=0$ we have $\la=\tau$ and the result holds. Now take $w>0$ and suppose the result holds for smaller $w$. By the definition of \fbc, there is a \trp $\mu$ satisfying one of the following.
\begin{itemize}
\item
$\mu$ is obtained from $\la$ by reducing some $\la_i$ by $4$ (and re-ordering). By the inductive hypothesis $\mu$ has no parts congruent to $1$ modulo $4$, and hence neither does $\la$. $\mu$ has \fbw $w-2$, so by hypothesis $\mu$ includes all the integers $3,7,\dots,4t+3$. Hence $\la$ includes the integers $3,7,\dots,4t-1$.
\item
$\la=\mu\sqcup(2)$. In this case, the assumption that $\mu$ satisfies the given conditions means that $\la$ does too.
\item
$\la=\mu\sqcup(3,1)$. This means that $3\notin\mu$; but since $\mu$ has \fbw $w-2$ the inductive hypothesis says that $3,7,\dots,4t+3\in\mu$, so that $t<0$, i.e.\ $w>2l+1$, contrary to assumption.
\qedhere
\end{itemize}
\end{pf}

\begin{cory}\label{rouqspinform3}
Suppose $B$ is a Rouquier block of $\tsss n$ with \fbc $\tau=(4l-1,4l-5,\dots,3)$ and \fbw $w$, and that $\la\in\cald(n)$ such that $\spn\la$ lies in $B$. Then $\la$ can be written in the form $\tau+4\alpha\sqcup2\beta$, where $\alpha\in\calp$, $\beta\in\cald$ and $2|\alpha|+|\beta|=w$. Furthermore, if $\la_r$ is any even part of $\la$, then $\la$ includes all positive integers less than $\la_r$ which are congruent to $3$ modulo $4$.
\end{cory}

\begin{pf}
By sorting the parts of $\la$ according to parity, we can write $\la=\nu\sqcup2\beta$ where all the positive parts of $\nu$ are odd. Then the \fbc of $\nu$ is also $\tau$, and the \fbw of $\nu$ is at most $w$. Since $B$ is a Rouquier block we have $w\ls 2l+1$, so by \cref{rouqlem3} all the positive parts of $\nu$ are congruent to $3$ modulo $4$. So if we let $m=\len\nu$ and let $\upsilon$ denote the \fbc $(4m-1,4m-5,\dots,3)$, then we have $\nu=\upsilon+4\alpha$ for some partition $\alpha$. But then the \fbc of $\nu$ is clearly $\upsilon$, and so we have $\upsilon=\tau$. So we can write $\la=\tau+4\alpha\sqcup2\beta$, and it follows that $2|\alpha|+|\beta|=w$.

The \fbw of $\nu$ is $2|\alpha|$, so applying the second statement of \cref{rouqlem3} with $\nu$ in place of $\la$, we find that $\nu$ (and hence $\la$) contains all the integers $3,7,\dots,4(l-|\alpha|)-1$. %The fact that $B$ is a Rouquier block means that $w\ls2l+1$. 
Now any even part of $\la$ equals $2\beta_s$ for some $s$, and
\[
2\beta_s\ls2|\beta|=2w-4|\alpha|\ls4l-4|\alpha|+2.
\]
Hence $\la$ contains all positive integers less than $2\beta_s$ which are congruent to $3$ modulo $4$.
\end{pf}

In a very similar way, we obtain the following.

\begin{lemma}\label{rouqlem1}
Suppose $\tau=(4l-3,4l-7,\dots,1)$ for $l\gs1$, and $\la$ is a \trp with \fbc $\tau$ and \fbw $w\ls2l$. Then $\la$ has no parts congruent to $3$ modulo $4$, and $\la$ includes all the integers $1,5,9,\dots,4t-3$, where $t=\lceil l-\frac12w\rceil$.
\end{lemma}

\begin{cory}\label{rouqspinform1}
Suppose $B$ is a Rouquier block of $\tsss n$ with \fbc $\tau=(4l-3,4l-7,\dots,1)$ and \fbw $w$, and that $\la\in\cald(n)$ such that $\spn\la$ lies in $B$. Then $\la$ can be written in the form $\tau+4\alpha\sqcup2\beta$, where $\alpha\in\calp$, $\beta\in\cald$ and $2|\alpha|+|\beta|=w$. Furthermore, if $\la_r$ is any even part of $\la$, then $\la$ includes all positive integers less than $\la_r$ which are congruent to $1$ modulo $4$.
\end{cory}

Note in particular that by \cref{rouqspinform3,rouqspinform1}, if $\la\in\cald$ has no even parts and $\spn\la$ lies in a Rouquier block, then the \fbw of $\la$ must be even.

\subsection{Decomposition numbers for Rouquier blocks}

An advantage of working with Rouquier blocks is that their decomposition numbers are relatively well understood. We summarise this situation, beginning with Hecke algebras, where the decomposition numbers for Rouquier blocks are known explicitly in terms of Littlewood--Richardson coefficients. This result is due to James and Mathas \cite[Corollary 2.6]{jmq-1}; here we only give two special cases.

\begin{thm}\label{heckerouqdec}
Suppose $\sigma=(c,c-1,\dots,1)$ and $w\ls c+1$.
\begin{enumerate}
\item
If $\la,\mu\in\calp(w)$, then
\begin{align*}
\odc_{(\sigma+2\la)(\sigma+2\mu)}&=\delta_{\la\mu}.
\\
\intertext{
\item
If $\la\in\calp(w-1)$ and $\mu\in\calp(w)$, then}
\odc_{(\sigma+2\la\sqcup(1^2))(\sigma+2\mu)}&=
\begin{cases}
1&\text{if $\la\subset\mu$}\\
0&\text{if $\la\not\subset\mu$}.
\end{cases}
\end{align*}
\end{enumerate}
\end{thm}

Now we consider Rouquier blocks of symmetric groups. The following result is due to Turner.

\begin{thm}[\xcite{turn}{Theorem 132}]\label{snrouqdec}
Suppose $B$ is a Rouquier block of $\sss n$ with $2$-core $\sigma=(c,c-1,\dots,1)$ and weight $w\ls c+1$. Then
\[
D_{(\sigma+2\la)(\sigma+2\mu)}=D_{\la\mu}
\]
for all $\la,\mu$ in $\calp(w)$.

Hence the adjustment matrix $A$ for the block $B$ is just the decomposition matrix of $\schur w$, i.e.
\[
A_{(\sigma+2\la)(\sigma+2\mu)}=D_{\la\mu}
\]
for $\la,\mu$ in $\calp(w)$.
\end{thm}

\subsection{Some virtual projective characters in a Rouquier block}\label{projcharsec}

We now examine the decomposition numbers in a Rouquier block of $\tsss n$ by considering projective characters.  We fix some notation.

\smallskip
\begin{mdframed}[innerleftmargin=3pt,innerrightmargin=3pt,innertopmargin=3pt,innerbottommargin=3pt,roundcorner=5pt,innermargin=-3pt,outermargin=-3pt]
\noindent\textbf{Assumptions and notation in force for \cref{projcharsec,somedecevensec}:}
$B$ is a Rouquier block of $\tsss n$ with $2$-core $\sigma$ and weight $w$. $\tau$ is the \fbc of $B$, i.e.\ the common \fbc of the \trps labelling spin characters in $B$.
\end{mdframed}

Recall that a projective character is one which vanishes on $2$-singular elements; a \emph{virtual} projective character is a $\bbz$-linear combination of projective characters. It is well known that induction and restriction send (virtual) projective characters to (virtual) projective characters, and hence so do the functors $\ee i$ and $\ff i$.

Recall that $\chm{\ \ }{\ \ }$ is the usual inner product on ordinary characters, and that $\ip\,\,$ is the standard inner product on the space of symmetric functions.

\begin{propn}\label{induceproj}
Suppose $\mu\in\calp(w)$. There is a projective character $\psi^\mu$ of $\tsss n$ with the following properties.
\begin{enumerate}
\item
$\chm{\psi^\mu}{\ord{\sigma+2\nu}}=\ip{e_{\mu'}}{s_\nu}$ for any $\nu\in\calp(w)$. In particular, $\chm{\psi^\mu}{\ord{\sigma+2\mu}}=1$, while $\chm{\psi^\mu}{\ord{\sigma+2\nu}}=0$ if $\mu\ndom\nu$.
\item
If the column lengths of $\mu$ are all even, say $\mu=\dup\alpha$, then $\chm{\psi^\mu}{\spn{\tau+4\beta}}=\ip{e_{\alpha'}}{s_\beta}$ for every $\beta\in\calp(w/2)$. In particular, $\chm{\psi^\mu}{\spn{\tau+4\alpha}}=1$ while $\chm{\psi^\mu}{\spn{\tau+4\beta}}=0$ if $\alpha\ndom\beta$.
\item
If $w$ is even but $\mu$ has at least one column of odd length, then $\chm{\psi^\mu}{\spn{\tau+4\beta}}=0$ for every $\beta\in\calp(w/2)$.
\end{enumerate}
\end{propn}

We prove \cref{induceproj} by applying induction functors. Note that all the addable nodes of $\sigma$ have the same residue; we assume for the rest of this section that they have residue $0$. (The proof in the opposite case is identical, but with $0$ and $1$ swapped throughout.)

Given $r\in\bbn$, consider the induction functor $\ffd1r\ffd0r$. The classical branching rule gives the following.

\begin{lemma}\label{induceord}
Suppose $r\ls v\ls w$, and $\xi$ is a partition with $2$-core $\sigma$ and $2$-weight $v-r$. Suppose $\nu\in\calp(v)$.
\begin{enumerate}
\item
If $\xi$ is $2$-regular, say $\xi=\sigma+2\kappa$, then $\chm{\ffd1r\ffd0r\ord\xi}{\ord{\sigma+2\nu}}=\ip{e_rs_\kappa}{s_\nu}$.
\item
If $\xi$ is $2$-singular, then $\chm{\ffd1r\ffd0r\ord\xi}{\ord{\sigma+2\nu}}=0$.
\end{enumerate}
\end{lemma}

\begin{pf}
This result is effectively a special case of \cite[Lemma 3.1]{ct}, but exploiting this requires a lot of translation of notation and transfer of results from one context to another, so we give a full proof here.

First note that $\ffd1r\ffd0r\ord\xi$ lies in the block with $2$-core $\sigma$ and weight $v$, and the condition $v\ls w$ guarantees that this is a Rouquier block. Now the branching rule says that $\ffd1r\ffd0r\ord\xi$ is the sum of $\ord\rho$ over all pairs of partitions $(\pi,\rho)$ such that $\xi\adr0r\pi\adr1r\rho$. This gives (2) straight away, because if $\xi$ is $2$-singular then by \cref{trpsinrouq} the length of $\xi$ is at least $\len\sigma+2$, so $\sigma+2\nu\nsupseteq\xi$.

So we are left with (1). Suppose $\xi$ is $2$-regular and write $\xi=\sigma+2\kappa$. If $\nu$ is obtained from $\kappa$ by adding $r$ nodes in distinct rows, then let $\pi$ be the partition obtained from $\xi$ by adding one node at the end of each of these rows. Then (since all the addable nodes of $\xi$ have residue $0$) we have $\xi\adr0r\pi\adr1r\sigma+2\nu$, so $\ord{\sigma+2\nu}$ occurs in $\ffd1r\ffd0r\ord\xi$, and clearly occurs once only. Conversely, if we have $\xi\adr0r\pi\adr1r{\sigma+2\nu}$, then $\xi\subset\sigma+2\nu$, so $\kappa\subset\nu$. Moreover, the nodes added to $\xi$ to obtain $\pi$ all have the same residue, so lie in different rows, and so $\sigma+2\nu$ differs from $\xi$ in at least $r$ rows. So $\nu$ differs from $\kappa$ in at least $r$ rows, so $\nu$ must be obtained from $\kappa$ by adding $r$ nodes in distinct rows.
\end{pf}

Now we prove a corresponding result for spin characters. Recall that we write $\spn\xi$ to mean ``either $\spn\xi_+$ or $\spn\xi_-$'' when $\xi\in\cald^-(n)$.

\begin{propn}\label{inducespin}
Suppose $r\ls v\ls w$ with $v$ even, and $\xi$ is a \trp with \fbc $\tau$ and \fbw $v-r$. Suppose $\beta\in\calp(v/2)$.
\begin{enumerate}
\item
Suppose all the parts of $\xi$ are odd, say $\xi=\tau+4\alpha$ for $\alpha\in\calp((v-r)/2)$. Then $\chm{\ffd0r\ffd1r\spn\xi}{\spn{\tau+4\beta}}=\ip{e_{r/2}s_\alpha}{s_\beta}$.
\item
If $\xi$ has at least one even part (and in particular if $r$ is odd), then $\chm{\ffd0r\ffd1r\spn\xi}{\spn{\tau+4\beta}}=0$.
\end{enumerate}
\end{propn}

\begin{pf}
Let $\rho=\tau+4\beta$. By the spin branching rule, $\spn\rho$ occurs in $\ffd0r\ffd1r\spn\xi$ if and only if there is $\pi$ such that $\xi\ads0r\pi\ads1r\rho$, i.e.\ $\xi$ can be obtained from $\rho$ by removing $r$ nodes of \spr $1$ followed by $r$ nodes of \spr $0$. The form of $\rho$ means that there are two nodes of \spr $1$ at the end of each non-empty row of $\rho$; removing both of these nodes from a row then exposes two nodes of \spr $0$ at the end of that row (or one, if the row initially only has length $3$). So the only way $\pi$ can have $r$ nodes of \spr $0$ that can be removed is if $r$ is even, and the nodes removed from $\rho$ to obtain $\pi$ occur in pairs at the end of $r/2$ different rows; then the only possible way to remove $r$ nodes of \spr $0$ from $\pi$ is to remove them in pairs from these same $r/2$ rows. Thus if $\xi\ads0r\pi\ads1r\rho$, then we must have $\rho_i=\xi_i+4$ for $r/2$ different values of $i$, with $\rho_i=\xi_i$ for all other values of $i$. Hence $\xi=\tau+4\alpha$, where $\alpha$ is obtained from $\beta$ by removing $r/2$ nodes from distinct rows.

This is sufficient to prove the \lcnamecref{inducespin}, except that in the case where $\beta$ is obtained from $\alpha$ by adding nodes in distinct rows we must show that the coefficient of $\spn\rho$ in $\ffd0r\ffd1r\spn{\tau+4\alpha}$ is $1$. But this follows from the induction version of \cref{spinbranchpower}.
\end{pf}

Now we can prove \cref{induceproj}.

\begin{pf}[Proof of \cref{induceproj}]
We use induction on $w$. In the case $w=0$, we can define $\psi^\varnothing$ to be the unique indecomposable projective character in the block with $2$-core $\sigma$ and weight $0$, i.e.\ $\prj\sigma$. By \cref{jreg,spinreg}, this equals $\ord\sigma+\spn\tau$.

Now suppose $w\gs1$, and let $\rho$ be the partition obtained by removing the last non-empty column (of length $r$, say) from $\mu$, so that $e_{\mu'}=e_re_{\rho'}$. By the inductive hypothesis the \lcnamecref{induceproj} holds with $\rho$ in place of $\mu$, and we define $\psi^\mu=\ffd1r\ffd0r\psi^\rho$. We must check that the conditions given in the \lcnamecref{induceproj} hold for $\psi^\mu$, given that they hold for $\psi^\rho$.

Suppose $\nu\in\calp(w)$. By \cref{induceord} we have
\begin{align*}
\chm{\psi^\mu}{\ord{\sigma+2\nu}}&=\chm{\ffd1r\ffd0r\psi^{\rho}}{\ord{\sigma+2\nu}}\\
&=\sum_{\kappa\in\calp(w-r)}\chm{\psi^{\rho}}{\ord{\sigma+2\kappa}}\chm{\ffd1r\ffd0r\ord{\sigma+2\kappa}}{\ord{\sigma+2\nu}}\\
&=\sum_{\kappa\in\calp(w-r)}\ip{e_{\rho'}}{s_\kappa}\ip{e_rs_\kappa}{s_\nu}\\
&=\ip{e_re_{\rho'}}{s_\nu}\\
&=\ip{e_{\mu'}}{s_\nu}.
\end{align*}

Now suppose $w$ is even and $\beta\in\calp(w/2)$. If $r$ is odd, then by \cref{inducespin}(2) we have
\[
\chm{\psi^\mu}{\spn{\tau+4\beta}}=\chm{\ffd1r\ffd0r\psi^{\rho}}{\spn{\tau+4\beta}}=0.
\]
So suppose instead that $r$ is even. If any of the column lengths of $\rho$ are odd, then by the inductive hypothesis $\chm{\psi^\rho}{\spn{\tau+4\gamma}}=0$ for every $\gamma\in\calp((w-r)/2)$, so by \cref{inducespin}(2) $\chm{\psi^\mu}{\spn{\tau+4\beta}}=0$. On the other hand, if the column lengths of $\mu$ are all even, then write $\mu=\dup\alpha$, so that $\rho=\dup\delta$ where $\delta$ is obtained by removing the last column from $\alpha$. Then
\begin{align*}
\chm{\psi^\mu}{\spn{\tau+4\beta}}&=\chm{\ffd1r\ffd0r\psi^{\rho}}{\spn{\tau+4\beta}}\\
&=\sum_{\gamma\in\calp((w-r)/2)}\chm{\psi^{\rho}}{\spn{\tau+4\gamma}}\chm{\ffd1r\ffd0r\spn{\tau+4\gamma}}{\spn{\tau+4\beta}}\\
&=\sum_{\gamma\in\calp((w-r)/2)}\ip{e_{\delta'}}{s_\gamma}\ip{e_{r/2}s_\gamma}{s_\beta}\\
&=\ip{e_{\alpha'}}{s_\beta}.\qedhere
\end{align*}
\end{pf}

\begin{egs}
We take $w=4$ and $\sigma=(3,2,1)$, so that $\tau=(5,1)$. We will show how to construct the characters $\psi^{(2^2)}$ and $\psi^{(2,1^2)}$. In both cases we start from $\psi^\varnothing=\ord{3,2,1}+\spn{5,1}$.

Applying $\ffd02\ffd12$ to $\psi^\varnothing$, we obtain
\begin{align*}
\psi^{(1^2)}&=\ord{5,4,1}+\ord{5,1^4}+\ord{3,2^3,1}\\
&\phantom=\ +\spn{9,1}+\spn{5,4,1}_++\spn{5,4,1}_-.\\
\intertext{
Applying $\ffd02\ffd12$ again, we obtain
}
\psi^{(2^2)}&=\ord{7,6,1}+\ord{7,4,3}+\ord{5,4,3,2}+\chi_1\\
&\phantom=\ +\spn{13,1}+\spn{9,5}+\chi_2,
\end{align*}
where $\chi_1$ is a sum of characters of the form $\ord\nu$ with $\nu$ $2$-singular, and $\chi_2$ is a sum of characters $\spn\xi$ for $\xi$ not of the form $\tau+4\beta$ for any $\beta$. Note that the characters $\ord{5,1^4}$ and $\ord{3,2^3,1}$ do not contribute any terms $\ord\nu$ with $\nu$ $2$-regular, and the characters $\spn{5,4,1}_\pm$ do not contribute any terms $\spn{\tau+4\beta}$.

Now we look at $\mu=(2,1^2)$. Starting with $\psi^\varnothing$ and applying $\ffd03\ffd13$, we obtain
\begin{align*}
\psi^{(1^3)}&=\ord{5,4,3}+\ord{5,4,1^3}+\ord{5,2^3,1}+\ord{3^3,2,1}\\
&\phantom=\ +\spn{9,2,1}_++\spn{9,2,1}_-+\spn{6,5,1}_++\spn{6,5,1}_-.\\
\intertext{
Applying $\ff0\ff1$, we obtain
}
\psi^{(2,1^2)}&=\ord{7,4,3}+\ord{5,4,3,2}+\chi,
\end{align*}
where $\chi$ is a sum of characters of the form $\ord\nu$ with $\nu$ $2$-singular, or $\spn\xi$ with $\xi$ not of the form $\tau+4\beta$.
\end{egs}

Next we want to prove a dual result to \cref{induceproj}, using virtual projective characters.

\begin{propn}\label{induceprojdual}
Suppose $\la\in\calp(w)$. There is a virtual projective character $\upsilon^\la$ with the following properties.
\begin{enumerate}
\item
$\chm{\upsilon^\la}{\ord{\sigma+2\nu}}=\ip{h_\la}{s_\nu}$ for any $\nu\in\calp(w)$.
\item
If $\la_i$ is even for every $i$, say $\la=2\gamma$, then $\chm{\upsilon^\la}{\spn{\tau+4\beta}}=(-1)^{w/2}\ip{h_\gamma}{s_\beta}$ for every $\beta\in\calp(w/2)$.
\item
If $w$ is even but $\la_i$ is odd for some $i$, then $\chm{\upsilon^\la}{\spn{\tau+4\beta}}=0$ for every $\beta\in\calp(w/2)$.
\end{enumerate}
\end{propn}

\begin{pf}
Recall the coefficients $\ke\la\mu$ defined by $h_\la=\sum_\mu\ke\la\mu e_\mu$, and define
\[
\upsilon^\la=\sum_{\mu\in\calp(w)}\ke\la\mu\psi^{\mu'}.
\]
Now for $\nu\in\calp(w)$ we have
\begin{align*}
\chm{\upsilon^\la}{\ord{\sigma+2\nu}}&=\ip{\textstyle\sum_\mu\ke\la\mu e_\mu}{s_\nu}\\
&=\ip{h_\la}{s_\nu}.
\end{align*}

Now suppose $w$ is even. Then for $\beta\in\calp(w/2)$ we have
\[
\chm{\upsilon^\la}{\spn{\tau+4\beta}}=\sum_{\alpha\in\calp(w/2)}\ke\la{2\alpha}\ip{e_\alpha}{s_\beta}
\]
by \cref{induceproj}. If $\la_i$ is odd for some $i$, then by \cref{simpleep} $\ke\la{2\alpha}=0$ for every $\alpha$, so that $\chm{\upsilon^\la}{\spn{\tau+4\beta}}=0$. On the other hand, if $\la=2\gamma$, then $\ke\la{2\alpha}=(-1)^{w/2}\ke\gamma{\alpha}$ by \cref{ekdoublep}, so that
\begin{align*}
\chm{\upsilon^\la}{\spn{\tau+4\beta}}&=(-1)^{w/2}\ip{\textstyle\sum_\alpha\ke\gamma\alpha e_\alpha}{s_\beta}\\
&=(-1)^{w/2}\ip{h_\gamma}{s_\beta}.\qedhere
\end{align*}
\end{pf}

We complete this subsection by examining a third set of virtual projective characters with a nice symmetry property.

By \cref{trpsinrouq}, the indecomposable projective characters are the characters $\prsj\mu$ for $\mu\in\calp(w)$; by Brauer reciprocity these satisfy $\chm{\prsj\mu}{\ord{\sigma+2\la}}=D_{(\sigma+2\la)(\sigma+2\mu)}$ for all $\la$, and hence by the first statement in \cref{snrouqdec} $\chm{\prsj\mu}{\ord{\sigma+2\la}}=D_{\la\mu}$ for all $\la,\mu\in\calp(w)$. The invertibility of the decomposition matrix of the Schur algebra then implies that for each $\mu$ there is a unique virtual projective character $\omega^\mu$ in $B$ such that $\chm{\omega^\mu}{\ord{\sigma+2\la}}=\delta_{\la\mu}$ for each $\la$. In fact, we can write
\[
\omega^\mu=\sum_{\nu\in\calp(w)}D_{\nu\mu}\v\prsj\nu,
\]
where (as set out in \cref{irrdecsec}) $D$ is the decomposition matrix of the Schur algebra in characteristic $2$, i.e.\ $D_{\la\mu}=\cm{\Delta(\la)}{L(\mu)}$. Our aim is to study the multiplicities of the spin characters $\spn{\tau+4\beta}$ in the $\omega^\mu$ when $w$ is even, which will enable us to deduce information about the decomposition numbers.

The triangularity of the projective characters $\psi^\mu$ (i.e.\ the second statement in \cref{induceproj}(1)) means that the set $\lset{\psi^\mu}{\mu\in\calp(w)}$ is a basis for the space of virtual projective characters in $B$. Hence each $\omega^\la$ can be written uniquely as a linear combination of the $\psi^\mu$. So define coefficients $a_{\la\mu}$ by $\omega^\la=\sum_{\mu\in\calp(w)} a_{\la\mu}\psi^\mu$. The next result shows that these coefficients can also be used to write $\omega^\la$ in terms of the virtual characters~$\upsilon^\mu$.

\begin{lemma}\label{omegaups}
For any $\la\in\calp(w)$,
\[
\omega^\la=\sum_{\mu\in\calp(w)}a_{\la'\mu'}\upsilon^\mu.
\]
\end{lemma}

\begin{pf}
For any $\nu\in\calp(w)$ we have
\begin{align*}
\chm{\sum_{\mu\in\calp(w)}a_{\la'\mu'}\upsilon^\mu}{\ord{\sigma+2\nu}}&=\sum_{\mu\in\calp(w)}a_{\la'\mu'}\ip{h_\mu}{s_\nu}\\
&=\sum_{\mu\in\calp(w)}a_{\la'\mu'}\ip{e_\mu}{s_{\nu'}}\\
&=\sum_{\mu\in\calp(w)}a_{\la'\mu}\ip{e_{\mu'}}{s_{\nu'}}\\
&=\sum_{\mu\in\calp(w)}a_{\la'\mu}\chm{\psi^\mu}{\ord{\sigma+2\nu'}}\\
&=\chm{\omega^{\la'}}{\ord{\sigma+2\nu'}}\\
&=\delta_{\la'\nu'}\\
&=\delta_{\la\nu}\\
&=\chm{\omega^\la}{\ord{\sigma+2\nu}}.
\end{align*}
Now the uniqueness property defining $\omega^\la$ gives the result.
\end{pf}

Now we can deduce the following symmetry property of the coefficients $\chm{\omega^\mu}{\spn{\tau+4\beta}}$.

\begin{propn}\label{virtualsym}
Suppose $w$ is even, $\la\in\calp(w)$ and $\beta\in\calp(w/2)$. Then
\[
\chm{\omega^\la}{\spn{\tau+4\beta}}=(-1)^{w/2}\chm{\omega^{\la'}}{\spn{\tau+4\beta'}}.
\]
\end{propn}

\begin{pf}
On the one hand, we have
\begin{align*}
\chm{\omega^\la}{\spn{\tau+4\beta}}&=\sum_{\mu\in\calp(w)}a_{\la\mu}\chm{\psi^\mu}{\spn{\tau+4\beta}}\\
&=\sum_{\alpha\in\calp(w/2)}a_{\la(\dup\alpha)}\ip{e_{\alpha'}}{s_\beta}
\\
\intertext{by \cref{induceproj}(2,3). On the other hand,}
\chm{\omega^{\la'}}{\spn{\tau+4\beta'}}&=\sum_{\mu\in\calp(w)}a_{\la\mu'}\chm{\upsilon^\mu}{\spn{\tau+4\beta'}}\tag*{by \cref{omegaups}}\\
&=\sum_{\gamma\in\calp(w/2)}a_{\la(2\gamma)'}(-1)^{w/2}\ip{h_\gamma}{s_{\beta'}}\tag*{by \cref{induceprojdual}}\\
&=(-1)^{w/2}\sum_{\gamma\in\calp(w/2)}a_{\la(\dup{\gamma'})}\ip{e_\gamma}{s_\beta}.
\end{align*}
Replacing $\gamma$ with $\alpha'$ gives the result.
\end{pf}

\subsection{Some decomposition numbers for Rouquier blocks of even weight}\label{somedecevensec}

We now use the results of the previous subsection to find some explicit decomposition numbers for spin characters in Rouquier blocks. We retain the assumptions set out at the start of \cref{projcharsec}, and \emph{we assume throughout \cref{somedecevensec} that $w$ is even}. Our aim is to study the decomposition numbers $\decs{(\tau+4\alpha)}{(\sigma+2\mu)}$, for $\alpha\in\calp(w/2)$ and $\mu\in\calp(w)$, and to express them in terms of the adjustment matrix of the Schur algebra.

We start by defining several matrices. Let $E$ be the matrix with rows indexed by $\calp(w/2)$ and columns by $\calp(w)$, with
\begin{align*}
E_{\alpha\mu}&=\decs{(\tau+4\alpha)}{(\sigma+2\mu)}.
\\
\intertext{Now define a matrix $J$ with the same indexing as $E$, with}
J_{\alpha\mu}&=
\begin{cases}
1&(\mu=\dup\alpha)\\
0&(\text{otherwise}).
\end{cases}
\end{align*}
Let $D$ be the decomposition matrix of the Schur algebra $\schur w$ in characteristic $2$, i.e.\ $D_{\la\mu}=[\Delta(\la):L(\mu)]$. Then (as explained in \cref{irrdecsec}) $D=\odc A$, where $\odc$ is the decomposition matrix of the $(-1)$-Schur algebra $\qschur w$ and $A$ is its adjustment matrix. Now we can state the main result of this section.

\begin{thm}\label{allodd}
Suppose $B$ is a Rouquier block of even weight $w$, and define $E,J,A$ as above. Then $E=JA$.
\end{thm}

For example, if $w=4$ then we have the following matrices, and we see that \cref{allodd} is true in this case.
\begin{alignat*}2
E&=
\begin{array}{@{}r@{\,}|@{\,}c@{\,}c@{\,}c@{\,}c@{\,}c}
\phantom{(2,1^2)}&\rt{(4)}&\rt{(3,1)}&\rt{(2^2)}&\rt{(2,1^2)}&\rt{(1^4)}\\\hline
(2)&\cdot&\cdot&1&\cdot&\cdot\\
(1^2)&\cdot&\cdot&1&\cdot&1
\end{array}
&\qquad
J&=
\begin{array}{@{}r@{\,}|@{\,}c@{\,}c@{\,}c@{\,}c@{\,}c}
\phantom{(2,1^2)}&\rt{(4)}&\rt{(3,1)}&\rt{(2^2)}&\rt{(2,1^2)}&\rt{(1^4)}\\\hline
(2)&\cdot&\cdot&1&\cdot&\cdot\\
(1^2)&\cdot&\cdot&\cdot&\cdot&1
\end{array}
\\[6pt]
D&=\begin{array}{@{}r@{\,}|@{\,}c@{\,}c@{\,}c@{\,}c@{\,}c}
&\rt{(4)}&\rt{(3,1)}&\rt{(2^2)}&\rt{(2,1^2)}&\rt{(1^4)}\\\hline
(4)&1&\cdot&\cdot&\cdot&\cdot\\
(3,1)&1&1&\cdot&\cdot&\cdot\\
(2^2)&\cdot&1&1&\cdot&\cdot\\
(2,1^2)&1&1&1&1&\cdot\\
(1^4)&1&\cdot&1&1&1
\end{array}&\qquad
\odc&=\begin{array}{@{}r@{\,}|@{\,}c@{\,}c@{\,}c@{\,}c@{\,}c}
&\rt{(4)}&\rt{(3,1)}&\rt{(2^2)}&\rt{(2,1^2)}&\rt{(1^4)}\\\hline
(4)&1&\cdot&\cdot&\cdot&\cdot\\
(3,1)&1&1&\cdot&\cdot&\cdot\\
(2^2)&\cdot&1&1&\cdot&\cdot\\
(2,1^2)&1&1&1&1&\cdot\\
(1^4)&1&\cdot&\cdot&1&1
\end{array}
\\[6pt]
A&=\begin{array}{@{}r@{\,}|@{\,}c@{\,}c@{\,}c@{\,}c@{\,}c}
&\rt{(4)}&\rt{(3,1)}&\rt{(2^2)}&\rt{(2,1^2)}&\rt{(1^4)}\\\hline
(4)&1&\cdot&\cdot&\cdot&\cdot\\
(3,1)&\cdot&1&\cdot&\cdot&\cdot\\
(2^2)&\cdot&\cdot&1&\cdot&\cdot\\
(2,1^2)&\cdot&\cdot&\cdot&1&\cdot\\
(1^4)&\cdot&\cdot&1&\cdot&1
\end{array}&&
\end{alignat*}

\cref{induceproj} gives us the following information about the matrix $E$, which shows that \cref{allodd} is true up to a triangular adjustment.

\begin{cory}\label{triang}
Suppose $\beta\in\calp(w/2)$ and $\mu\in\calp(w)$.
\begin{enumerate}
\item
If $\mu$ has at least one column of odd length, then $E_{\beta\mu}=0$.
\item
If the columns of $\mu$ all have even length, say $\mu=\dup\alpha$, then $E_{\alpha\mu}=1$ while $E_{\beta\mu}=0$ unless $\alpha\dom\beta$.
\end{enumerate}
\end{cory}

\begin{pf}
$E_{\beta\mu}$ is the coefficient $\chm{\prsj\mu}{\spn{\tau+4\beta}}$.  So we examine the projective characters $\prsj\la$, for $\la\in\calp(w)$.

By the triangularity of the decomposition matrix of the symmetric group we have
\[
\chm{\prsj\mu}{\ord{\sigma+2\mu}}=1,\qquad\chm{\prsj\mu}{\ord{\sigma+2\nu}}=0\text{ for }\mu\ndom\nu.
\]
By \cref{induceproj} the same property holds if we replace $\prsj\mu$ with the projective character $\psi^\mu$ introduced in that \lcnamecref{induceproj}. Since $\psi^\mu$ is a linear combination of the characters $\prsj\la$ with non-negative coefficients, we can obtain $\prsj\mu$ from $\psi^\mu$ by subtracting multiples of the characters $\prsj\nu$ for $\nu\domsby\mu$. That is,
\[
\prsj\mu=\psi^\mu-\sum_{\nu\domsby\mu}x_\nu\prsj\nu
\]
with each $x_\nu$ a non-negative integer. In particular, for any $\beta\in\calp(w/2)$
\[
\chm{\prsj\mu}{\spn{\tau+4\beta}}\ls\chm{\psi^\mu}{\spn{\tau+4\beta}},
\]
so (by \cref{induceproj}) $\chm{\prsj\mu}{\spn{\tau+4\beta}}=0$ if $\mu$ has at least one column of odd length or $\mu=\dup\alpha$ with $\alpha\ndom\beta$. Hence if $\mu=\dup\alpha$, the coefficient $\chm{\prsj\nu}{\spn{\tau+4\alpha}}$ is zero for any $\nu\domsby\mu$, and so
\[
\chm{\prsj\mu}{\spn{\tau+4\alpha}}=\chm{\psi^\mu}{\spn{\tau+4\alpha}}=1.\qedhere
\]
\end{pf}

Now we note that the same property holds for the matrix $EA\v$.

\begin{cory}\label{ehv}
Suppose $\beta\in\calp(w/2)$ and $\mu\in\calp(w)$.
\begin{enumerate}
\item
If $\mu$ has at least one column of odd length, then $(EA\v)_{\beta\mu}=0$.
\item
If the columns of $\mu$ all have even length, say $\mu=\dup\alpha$, then $(EA\v)_{\alpha\mu}=1$ while $(EA\v)_{\beta\mu}=0$ unless $\alpha\dom\beta$.
\end{enumerate}
\end{cory}

\begin{pf}
We have $(EA\v)_{\beta\mu}=\sum_{\la\in\calp(w)}E_{\beta\la}A\v_{\la\mu}$, and by \cref{triang}(1) we may restrict the range of summation to those $\la$ with all columns of even length; that is,
\[
(EA\v)_{\beta\mu}=\sum_{\rho\in\calp(w/2)}E_{\beta(\dup\rho)}A\v_{(\dup\rho)\mu}.
\]
If $\mu$ has a column of odd length, then by \cref{adj} $A\v_{(\dup\rho)\mu}=0$ for every $\rho$, which gives (1). If $\mu=\dup\alpha$, then by \cref{adj}
\[
(EA\v)_{\beta\mu}=\sum_{\rho\in\calp(w/2)}E_{\beta(\dup\rho)}A\v_{(\dup\rho)(\dup\alpha)}=\sum_{\rho\in\calp(w/2)}E_{\beta(\dup\rho)}K\v_{\rho\alpha},
\]
where $K$ is the decomposition matrix of the Schur algebra $\schur{w/2}$. The term $E_{\beta(\dup\rho)}$ is non-zero only if $\rho\dom\beta$ (and equals $1$ when $\rho=\beta$), while the term $K\v_{\rho\alpha}$ is non-zero only if $\alpha\dom\rho$ (and is $1$ if $\alpha=\rho$). The result follows.
\end{pf}

\cref{ehv} may be alternatively phrased as follows: there is a square matrix $T$ with rows and columns indexed by $\calp(w/2)$, such that $EA\v=TJ$. Furthermore, $T$ is lower unitriangular in the sense that $T_{\alpha\alpha}=1$ for each $\alpha$, while $T_{\alpha\beta}=0$ for $\beta\ndom\alpha$.

Our aim is to prove that $T$ is the identity matrix. To do this, we consider the matrices $ED\v$ and $J\odc\v$, and show that they both satisfy a certain symmetry property. Given a matrix $B$ with rows indexed by $\calp(w/2)$ and columns by $\calp(w)$, we say that $B$ is \emph{\funny} if
\[
B_{\alpha'\la'}=(-1)^{w/2}B_{\alpha\la}
\]
for every $\alpha\in\calp(w/2)$ and $\la\in\calp(w)$.

\begin{propn}\label{jfunny}
The matrix $J\odc\v$ is \funny.
\end{propn}

\begin{pf}
The definition of $J$ means that $(J\odc\v)_{\alpha\mu}=\odc_{(\dup\alpha)\mu}\v$, and this is given explicitly by \cref{steinberg} as $(-1)^{w/2}\epsilon(\mu)\kappa(\alpha,\mu)$, where $\kappa(\alpha,\mu)$ is defined in \cref{inversesec}. Now $\mu$ and $\mu'$ have the same $2$-core, and if this $2$-core if not $\varnothing$ then $\kappa(\alpha,\mu)=\kappa(\alpha',\mu')=0$, so that $(J\odc\v)_{\alpha\mu}=(J\odc\v)_{\alpha'\mu'}=0$. So assume that $\mu$ and $\mu'$ have empty $2$-core, and therefore have $2$-weight $w/2$. Then by \cref{parityconj} we have $\epsilon(\mu)\epsilon(\mu')=(-1)^{w/2}$, so we just need to show that $\kappa(\alpha,\mu)=\kappa(\alpha',\mu')$ for any $\alpha$. It is a standard property of Littlewood--Richardson coefficients that $a^\alpha_{\beta\gamma}=a^{\alpha'}_{\gamma'\beta'}$ for any $\alpha,\beta,\gamma$, so the result follows from \cref{quoconj}.
\end{pf}

\begin{propn}\label{efunny}
The matrix $ED\v$ is \funny.
\end{propn}

\begin{pf}
Recall that $E_{\alpha\mu}$ is the multiplicity $\chm{\prsj\mu}{\spn{\tau+4\alpha}}$. Recall also that the virtual projective characters $\omega^\mu$ from \S\ref{projcharsec} are given by
\[
\omega^\mu=\sum_{\nu\in\calp(w)}D_{\nu\mu}\v\prsj\nu.
\]
This means that $(ED\v)_{\alpha\mu}=\chm{\omega^\mu}{\spn{\tau+4\alpha}}$. Now \cref{virtualsym} gives the required result.
\end{pf}

For example, in the case $w=4$ the matrix $J\odc\v=ED\v$ is the following.
\[
\begin{array}{r|ccccc}
&\rt{(4)}&\rt{(3,1)}&\rt{(2^2)}&\rt{(2,1^2)}&\rt{(1^4)}\\\hline
(2)&1&-1&1&\cdot&\cdot\\
(1^2)&\cdot&\cdot&1&-1&1
\end{array}
\]

Now we can complete the proof of \cref{allodd}.

\begin{pf}[Proof of \cref{allodd}]
We have seen that there is a matrix $T$ such that $EA\v=TJ$, and that $T$ is unitriangular in the sense that $T_{\alpha\alpha}=1$ and $T_{\alpha\beta}=0$ if $\alpha\ndomby\beta$. So to prove that $T$ is the identity matrix, we just need to prove that $T_{\alpha\beta}=0$ for all $\alpha\domsby\beta$, and we do this by induction on $\beta$ using the dominance order. So assume that we have proved that $T_{\alpha\gamma}=0$ whenever $\gamma\doms\alpha,\beta$.

Multiplying the equation $EA\v=TJ$ by $\odc\v$, we obtain $ED\v=TJ\odc\v$. Consider the column of $J\odc\v$ labelled by the partition $\dup{\beta'}$. For any $\gamma$ we have $(J\odc\v)_{\gamma(\dup{\beta'})}=\odc_{(\dup\gamma)(\dup{\beta'})}\v$, and the triangularity of $\odc$ means that $(J\odc\v)_{\gamma(\dup{\beta'})}=1$ if $\gamma=\beta'$, and $0$ if $\gamma\ndomby\beta'$. The triangularity of $T$ then means that the same is true with $J\odc\v$ replaced by $TJ\odc\v=ED\v$.

Using the fact that both $J\odc\v$ and $ED\v$ are \funny (and the fact that the dominance order is reversed by conjugating partitions) we get
\[
(J\odc\v)_{\alpha(2\beta)}=
\begin{cases}
(-1)^{w/2}&\text{if }\alpha=\beta\\
0&\text{if }\alpha\ndom\beta
\end{cases}
\]
and the same with $ED\v$ in place of $J\odc\v$. Hence for $\alpha\domsby\beta$ we have
\begin{align*}
0&=(ED\v)_{\alpha(2\beta)}\\
&=(TJ\odc\v)_{\alpha(2\beta)}\\
&=\sum_\gamma T_{\alpha\gamma}(J\odc\v)_{\gamma(2\beta)}.
\end{align*}
The term $(J\odc\v)_{\gamma(2\beta)}$ is zero unless $\gamma\dom\beta$. If $\gamma\doms\beta$ then by our induction hypothesis $T_{\alpha\gamma}=0$, so we only need to consider the term with $\gamma=\beta$, and we obtain $0=(-1)^{w/2}T_{\alpha\beta}$.

Hence $T$ is the identity matrix, and therefore $E=JA$.
\end{pf}

\subsection{Some decomposition numbers for Rouquier blocks of odd weight}\label{somedecoddsec}

We now prove a similar result to \cref{allodd} for the case where $w$ is odd. We retain the assumptions set out at the start of \cref{projcharsec}, and \emph{we assume throughout \cref{somedecoddsec} that $w$ is odd}.

As before, we let $D$ be the decomposition matrix of the Schur algebra $\schur w$, and factorise $D$ as $\odc A$, where $\odc$ is the decomposition matrix of the corresponding $(-1)$-Schur algebra and $A$ is the adjustment matrix. We now define $E$ to be the matrix with rows indexed by $\calp((w-1)/2)$ and columns by $\calp(w)$, and
\begin{align*}
E_{\alpha\mu}&=\decs{(\tau+4\alpha\sqcup(2))}{(\sigma+2\mu)}.
\\
\intertext{$J$ is defined to have the same indexing as $E$, with}
J_{\alpha\mu}&=
\begin{cases}
1&(\mu=\dup\alpha\sqcup(1))\\
0&(\text{otherwise}).
\end{cases}
\end{align*}

Now we have the following.

\begin{thm}\label{allodd2}
Suppose $B$ is a Rouquier block of odd weight $w$, and define $E,J,A$ as above. Then $E=JA$.
\end{thm}

For example, taking $w=5$, we have the following matrices.

{
\allowdisplaybreaks
\begin{alignat*}2
E&=
\begin{array}{@{}r@{\,}|@{\,}c@{\,}c@{\,}c@{\,}c@{\,}c@{\,}c@{\,}c}
\phantom{(3,1^2)}&\rt{(5)}&\rt{(4,1)}&\rt{(3,2)}&\rt{(3,1^2)}&\rt{(2^2,1)}&\rt{(2,1^3)}&\rt{(1^5)}\\\hline
(2)&\cdot&\cdot&\cdot&\cdot&1&\cdot&\cdot\\
(1^2)&\cdot&\cdot&\cdot&\cdot&1&\cdot&1
\end{array}
&\qquad
J&=
\begin{array}{@{}r@{\,}|@{\,}c@{\,}c@{\,}c@{\,}c@{\,}c@{\,}c@{\,}c}
\phantom{(3,1^2)}&\rt{(5)}&\rt{(4,1)}&\rt{(3,2)}&\rt{(3,1^2)}&\rt{(2^2,1)}&\rt{(2,1^3)}&\rt{(1^5)}\\\hline
(2)&\cdot&\cdot&\cdot&\cdot&1&\cdot&\cdot\\
(1^2)&\cdot&\cdot&\cdot&\cdot&\cdot&\cdot&1
\end{array}
\\[6pt]
D&=\begin{array}{@{}r@{\,}|@{\,}c@{\,}c@{\,}c@{\,}c@{\,}c@{\,}c@{\,}c}
\phantom{(3,1^2)}&\rt{(5)}&\rt{(4,1)}&\rt{(3,2)}&\rt{(3,1^2)}&\rt{(2^2,1)}&\rt{(2,1^3)}&\rt{(1^5)}\\\hline
(5)&1&\cdot&\cdot&\cdot&\cdot&\cdot&\cdot\\
(4,1)&\cdot&1&\cdot&\cdot&\cdot&\cdot&\cdot\\
(3,2)&1&\cdot&1&\cdot&\cdot&\cdot&\cdot\\
(3,1^2)&2&\cdot&1&1&\cdot&\cdot&\cdot\\
(2^2,1)&1&\cdot&1&1&1&\cdot&\cdot\\
(2,1^3)&\cdot&1&\cdot&\cdot&\cdot&1&\cdot\\
(1^5)&1&\cdot&\cdot&1&1&\cdot&1
\end{array}&\qquad
\odc&=\begin{array}{@{}r@{\,}|@{\,}c@{\,}c@{\,}c@{\,}c@{\,}c@{\,}c@{\,}c}
\phantom{(3,1^2)}&\rt{(5)}&\rt{(4,1)}&\rt{(3,2)}&\rt{(3,1^2)}&\rt{(2^2,1)}&\rt{(2,1^3)}&\rt{(1^5)}\\\hline
(5)&1&\cdot&\cdot&\cdot&\cdot&\cdot&\cdot\\
(4,1)&\cdot&1&\cdot&\cdot&\cdot&\cdot&\cdot\\
(3,2)&\cdot&\cdot&1&\cdot&\cdot&\cdot&\cdot\\
(3,1^2)&1&\cdot&1&1&\cdot&\cdot&\cdot\\
(2^2,1)&\cdot&\cdot&1&1&1&\cdot&\cdot\\
(2,1^3)&\cdot&1&\cdot&\cdot&\cdot&1&\cdot\\
(1^5)&1&\cdot&\cdot&1&\cdot&\cdot&1
\end{array}\\[6pt]
A&=\begin{array}{@{}r@{\,}|@{\,}c@{\,}c@{\,}c@{\,}c@{\,}c@{\,}c@{\,}c}
\phantom{(3,1^2)}&\rt{(5)}&\rt{(4,1)}&\rt{(3,2)}&\rt{(3,1^2)}&\rt{(2^2,1)}&\rt{(2,1^3)}&\rt{(1^5)}\\\hline
(5)&1&\cdot&\cdot&\cdot&\cdot&\cdot&\cdot\\
(4,1)&\cdot&1&\cdot&\cdot&\cdot&\cdot&\cdot\\
(3,2)&1&\cdot&1&\cdot&\cdot&\cdot&\cdot\\
(3,1^2)&\cdot&\cdot&\cdot&1&\cdot&\cdot&\cdot\\
(2^2,1)&\cdot&\cdot&\cdot&\cdot&1&\cdot&\cdot\\
(2,1^3)&\cdot&\cdot&\cdot&\cdot&\cdot&1&\cdot\\
(1^5)&\cdot&\cdot&\cdot&\cdot&1&\cdot&1
\end{array}&&
\end{alignat*}
}

To prove \cref{allodd2}, we use restriction, comparing $B$ with the block $B_-$ of weight $w-1$ and $2$-core $\sigma$ and exploiting \cref{allodd}. We define $\ee\bullet$ to be: 
\begin{itemize}
\item
$\ee0\ee1$ if $\sigma$ is one of $\varnothing,(2,1),(4,3,2,1),\dots$ (and hence $\tau$ is one of $\varnothing,(3),(7,3),\dots$);
\item
$\ee1\ee0$ if $\sigma$ is one of $(1),(3,2,1),(5,4,3,2,1),\dots$ (and hence $\tau$ is one of $(1),(5,1),(9,5,1),\dots$).
\end{itemize}

As in the even-weight case, for each $\la\in\calp(w)$ we define $\omega^\la$ to be the unique virtual projective character in $B$ for which $\chm{\omega^\la}{\ord{\sigma+2\mu}}=\delta_{\la\mu}$.

Now consider applying $\ee\bullet$ to a character in $B$, and examining the coefficients of $\ord{\sigma+2\la}$ and $\spn{\tau+4\alpha}$ in this restriction.

\begin{lemma}\label{ebullord}
Suppose $\la\in\calp(w-1)$ and $\xi\in\calp(n)$ with $\ord\xi$ in $B$. Then $\chm{\ee\bullet\ord\xi}{\ord{\sigma+2\la}}$ equals:
\[
\begin{cases}
1&\text{if $\xi=\sigma+2\mu$ where $\mu$ is obtained by adding a node to $\la$;}\\
1&\text{if $\xi=\sigma+2\la\sqcup(1^2)$;}\\
0&\text{otherwise.}
\end{cases}
\]
\end{lemma}

\begin{pf}
This is a simple application of the branching rule. Assume that $\ee\bullet=\ee0\ee1$ (the other case is very similar). Then by the branching rule $\chm{\ee\bullet\ord\xi}{\ord{\sigma+2\la}}$ equals $1$ if there is a partition $\pi$ with $\sigma+2\la\adr01\pi\adr11\xi$, and $0$ otherwise. Since $\len\sigma\gs w-1\gs\len\la$, all the addable nodes of $\sigma+2\la$ have residue $0$, so if $\pi$ exists then the two nodes added to $\sigma+2\la$ to obtain $\xi$ must be adjacent. If these nodes are both in the same row, say row $i$, then $\xi=\sigma+2\mu$, where $\mu$ is obtained from $\la$ by adding a node in row $i$. On the other hand, if the added nodes are both in the same column, then (since $\sigma+2\la$ is $2$-regular) this must be column $1$, and hence $\xi=\sigma+2\la\sqcup(1^2)$.
\end{pf}

\begin{cory}\label{ebullomega}
Suppose $\mu\in\calp(w)$ and $\la\in\calp(w-1)$. Then $\chm{\ee\bullet\omega^\mu}{\ord{\sigma+2\la}}$ equals $2$ if $\la\subset\mu$, and $0$ otherwise.  Hence $\ee\bullet\omega^\mu=2\sum_\la\omega^\la$, summing over those partitions $\la$ obtained by removing a node from~$\mu$.
\end{cory}

\begin{pf}
To apply \cref{ebullord}, we need to know the multiplicities $\chm{\omega^\mu}{\ord{\sigma+2\la\sqcup(1^2)}}$. We have
\begin{align*}
\chm{\omega^\mu}{\ord{\sigma+2\la\sqcup(1^2)}}&=\sum_{\nu\in\calp(w)}D\v_{\nu\mu}\chm{\prsj\nu}{\ord{\sigma+2\la\sqcup(1^2)}}\\
&=\sum_{\nu\in\calp(w)}A\v_{(\sigma+2\nu)(\sigma+2\mu)}D_{(\sigma+2\la\sqcup(1^2))(\sigma+2\nu)}\tag*{by \cref{snrouqdec}}\\
&=(DA\v)_{(\sigma+2\la\sqcup(1^2))(\sigma+2\mu)}\\
&=\odc_{(\sigma+2\la\sqcup(1^2))(\sigma+2\mu)},
\end{align*}
and by \cref{heckerouqdec}(2) this equals $1$ if $\la\subset\mu$ and $0$ otherwise.

Now we can apply \cref{ebullord}:
\begin{align*}
\chm{\ee\bullet\omega^\mu}{\ord{\sigma+2\la}}&=\sum_\xi\chm{\omega^\mu}{\ord\xi}\chm{\ee\bullet\ord\xi}{\ord{\sigma+2\la}}\\
&=\chm{\omega^\mu}{\ord{\sigma+2\la\sqcup(1^2)}}+\sum_{\substack{\nu\in\calp(w)\\\nu\supset\la}}\chm{\omega^\mu}{\ord{\sigma+2\nu}}\tag*{by \cref{ebullord}.}
\end{align*}
The first term is $1$ if $\la\subset\mu$ and $0$ otherwise from above, while the second term is $1$ if $\la\subset\mu$ and $0$ otherwise by the definition of $\omega^\mu$.

Now the final statement follows from the fact that the set $\lset{\omega^\la}{\la\in\calp(w-1)}$ is a basis for the space of virtual projective characters in $B_-$, with $\chm{\omega^\la}{\ord{\sigma+2\nu}}=\delta_{\la\nu}$ for all $\la,\nu\in\calp(w-1)$.
\end{pf}

Now we consider applying $\ee\bullet$ to spin characters.

\begin{lemma}\label{ebullspin}
Suppose $\alpha\in\calp((w-1)/2)$ and $\pi\in\cald(n)$ with $\spn\pi$ in $B$. Then $\chm{\ee\bullet\spn\pi}{\spn{\tau+4\alpha}}$ equals $1$ if $\pi=\tau+4\alpha\sqcup(2)$, and $0$ otherwise. Hence if $\omega$ is a virtual projective character in $B$, then $\chm{\omega}{\spn{\tau+4\alpha\sqcup(2)}}=\frac12\chm{\ee\bullet\omega}{\spn{\tau+4\alpha}}$.
\end{lemma}

\begin{pf}
We apply the spin branching rule. Assume first that $\ee\bullet=\ee0\ee1$. In order to have $\chm{\ee\bullet\spn\pi}{\spn{\tau+4\alpha}}>0$, we need $\tau+4\alpha\ads01\rho\ads11\pi$ for some $\rho\in\cald$. We have $\len\tau\gs(w-1)/2\gs\len\alpha$, so every addable node of $\tau+4\alpha$ has \spr $0$, and every addable node also has a node of \spr $0$ to its immediate right, with the exception of the addable node in the first column. So the only possible way of adding two nodes of \sprs $0$ and then $1$ is to add the addable node at the bottom of the first column, and then the addable node at the bottom of the second column, which gives the partition $\pi=\tau+4\alpha\sqcup(2)$. The spin branching rule then gives the coefficient $\chm{\ee\bullet\spn\pi}{\spn{\tau+4\alpha}}$ as $1$.

The case where $\ee\bullet=\ee1\ee0$ is very similar: in this case $\len\tau>\len\alpha$, so $\tau+4\alpha$ has addable nodes in columns $1$ and $2$, and the only way to add a node of \spr $1$ followed by a node of \spr $0$ and end up with a \trp is to add the node in column $2$ followed by the node in column $1$.

The second statement now follows easily, except that we must take the two characters $\spn{\tau+4\alpha\sqcup(2)}_\pm$ into account. Because $\omega$ is a virtual projective character, we have $\chm{\omega}{\spn{\tau+4\alpha\sqcup(2)}_+}=\chm{\omega}{\spn{\tau+4\alpha\sqcup(2)}_-}$ (so it makes sense to write $\chm{\omega}{\spn{\tau+4\alpha\sqcup(2)}}$ for either of these multiplicities), and the characters $\spn{\tau+4\alpha\sqcup(2)}_\pm$ both contribute terms $\spn{\tau+4\alpha}$ when we apply $\ee\bullet$. Hence we have the factor $\frac12$ in the final expression.
\end{pf}

\begin{pf}[Proof of \cref{allodd2}]
Take $\mu\in\calp(w)$ and $\alpha\in\calp((w-1)/2)$, and let $M=\lset{\la\in\calp(w-1)}{\la\subset\mu}$. Let $E_-$, $J_-$, $D_-$, $\odc_-$ be the matrices defined in \cref{somedecevensec} for the block with $2$-core $\sigma$ and weight $w-1$ (so $D_-$ is the decomposition matrix of $\schur{w-1}$, and so on). Then
\allowdisplaybreaks
\begin{align*}
(ED\v)_{\alpha\mu}&=\sum_{\nu\in\calp(w)}D_{\nu\mu}\v\chm{\prsj\nu}{\spn{\tau+4\alpha\sqcup(2)}}\\
&=\chm{\omega^\mu}{\spn{\tau+4\alpha\sqcup(2)}}\\
&=\tfrac12\chm{\ee\bullet\omega^\mu}{\spn{\tau+4\alpha}}\tag*{by \cref{ebullspin}}\\
&=\sum_{\la\in M}\chm{\omega^\la}{\spn{\tau+4\alpha}}\tag*{by \cref{ebullomega}}\\
&=\sum_{\la\in M}(E_-D_-\v)_{\alpha\la}\\
&=\sum_{\la\in M}(J_-\odc_-\v)_{\alpha\la}\tag*{by \cref{allodd}}\\
&=\sum_{\la\in M}(\odc_-\v)_{(\dup\alpha)\la}.
\end{align*}
Now comparing \cref{steinberg,steinbergod} we see that this equals $\odc\v_{(\dup\alpha\sqcup(1))\mu}=(J\odc\v)_{\alpha\mu}$. Hence $ED\v=J\odc\v$, and so $E=JA$.
\end{pf}

\subsection{Irreducibility of spin characters in Rouquier blocks}\label{irrspinrouq}

From the results of the last two subsections we can finally deduce some information about irreducibility of spin characters modulo $2$.

We continue to assume that $B$ is a Rouquier block with weight $w$, $2$-core $\sigma$ and \fbc $\tau$. We no longer make any assumption on the parity of $w$. In this section we will classify the spin characters in $B$ that remain irreducible in characteristic $2$.

\begin{propn}\label{rouqirred}
Suppose $\alpha\in\calp(\lfloor w/2\rfloor)$, and let
\[
\la=
\begin{cases}
\tau+4\alpha&(\text{if $w$ is even})\\
\tau+4\alpha\sqcup(2)&(\text{if $w$ is odd}).
\end{cases}
\]
Then $\mspn\la$ is irreducible \iff $\alpha$ is a $2$-Carter partition.
\end{propn}

\begin{pf}
By \cref{allodd} or \cref{allodd2}, the row of the decomposition matrix corresponding to $\la$ equals the row of the adjustment matrix of $\schur w$ corresponding to $\dup\alpha$ or $\dup\alpha\sqcup(1)$. By \cref{adj} or \cref{adjod}, the sum of the entries of this row equals the sum of the entries of the row of the decomposition matrix of $\schur{\lfloor w/2\rfloor}$ labelled by $\alpha$. Thus $\mspn\la$ is irreducible \iff the Weyl module $\Delta(\alpha)$ is irreducible, and by \cref{irredweyl} this happens \iff $\alpha$ is a $2$-Carter partition.
\end{pf}

We now want to show that if $\la$ is not of the form $\tau+4\alpha$ or $\tau+4\alpha\sqcup(2)$ then $\mspn\la$ is reducible.

\begin{propn}\label{sameregdoub}
Suppose $\spn\la$ lies in $B$, and that $\la$ has the form $\tau+4\alpha\sqcup(4a)$ for some $a\gs1$. Then there is $\mu\in\cald(n)$ such that $\la\dblreg=\mu\dblreg$ and $\deg\spn\la>\deg\spn\mu$.
\end{propn}

\begin{pf}
We assume first that $\tau=(4l-1,4l-5,\dots,3)$ for $l\gs0$. Let $k$ be the length of $\alpha$. By \cref{rouqspinform3}, $\la$ includes all the integers $3,7,\dots,4a-1$, which is the same as saying that $k\ls l-a$. So if we let $m=l-a-k$, then the partition $(\la_{k+1},\la_{k+2},\dots)$ equals $(4a+4m-1,4a+4m-3,\dots,3)\sqcup(4a)$, i.e.\ the partition $\la^{a,m}$ from \cref{thirdrdim}. Furthermore, if $k>0$ then $\la_k>4a+4m+3$, so we can define a partition
\begin{align*}
\mu=(&4l-1+4\alpha_1,\dots,4l-4k+3+4\alpha_k,\\
&4a+4m+3,4a+4m-1,\dots,4m+7,\\
&4m-1,4m-5,\dots,3),
\end{align*}
and we then have $\la_i=\mu_i$ for $i=1,\dots,k$, while the partition $(\mu_{k+1},\mu_{k+2},\dots)$ coincides with the partition $\mu^{a,m}$ from \cref{thirdrdim}. Now by combining \cref{thirdrdim,dimrowrem} we obtain $\la\dblreg=\mu\dblreg$ and $\deg\spn\la>\deg\spn\mu$.

The case where $\tau=(4l-3,4l-7,\dots,1)$ is the same, using \cref{firstrdim,rouqspinform1} instead of \cref{thirdrdim,rouqspinform3}.
\end{pf}

In the same way (using \cref{secondrdim,fourthrdim} rather than \cref{firstrdim,thirdrdim}) we obtain the following.

\begin{propn}\label{sameregdoub2}
Suppose $\spn\la$ lies in $B$, and $\la$ has the form $\tau+4\alpha\sqcup(4a+2)$ for some $a\gs1$. Then there is $\mu\in\cald(n)$ such that $\la\dblreg=\mu\dblreg$ and $\deg\spn\la>\deg\spn\la$.
\end{propn}

Now we can give the main result for Rouquier blocks.

\begin{cory}\label{mainirrrouq}
Suppose $\la\in\cald$ with \fbc $\tau$, and that $\spn\la$ lies in a Rouquier block. Then $\mspn\la$ is irreducible \iff $\la$ has the form $\tau+4\alpha$ or $\tau+4\alpha\sqcup(2)$ for $\alpha$ a $2$-Carter partition.
\end{cory}

\begin{pf}
By \cref{spinreg}, we may assume $\la$ has at most one non-zero even part. Hence by \cref{rouqspinform3} or \cref{rouqspinform1} we can write $\la=\tau+4\alpha$ or $\tau+4\alpha\sqcup(2b)$ for some natural number $b$. In the case where $\la=\tau+4\alpha$ or $\tau+4\alpha\sqcup(2)$, the result follows from \cref{rouqirred}. If $\la=\tau+4\alpha\sqcup(2b)$ with $b\gs2$, then by \cref{sameregdoub} or \cref{sameregdoub2} we can find another \trp $\mu$ with the same regularised double as $\la$ and such that $\spn\mu$ has smaller degree than $\spn\la$. Hence by \cref{samerdoub} $\mspn\la$ is reducible.
\end{pf}

\subsection{Separated partitions}

Our aim in the remainder of this section is to extend the results of \cref{irrspinrouq} beyond Rouquier blocks. So now we drop the assumption that we are working in a Rouquier block, but we make a definition motivated by \cref{rouqspinform3,rouqspinform1}: given $\la\in\cald$ and $i=1$ or $3$, we say that $\la$ is \emph{$i$-\qs} if
\begin{itemize}
\item
$\la$ has at most one non-zero even part,
\item
all the odd parts of $\la$ are congruent to $i$ modulo $4$, and
\item
if $\la$ includes an even part $2b$, then $\la$ includes all positive integers less than $2b$ which are congruent to $i$ modulo $4$.
\end{itemize}
Say that $\la$ is \emph{\qs} if it is either $1$- or $3$-\qs. If we let $\tau$ be the \fbc of $\la$, then $\la$ is \qs if $\la$ has one of the following two forms:
\begin{itemize}
\item
$\tau+4\alpha$, where $\alpha\in\calp$ with $\len\alpha\ls\len\tau$;
\item
$\tau+4\alpha\sqcup(2b)$, where $b>0$, $\alpha\in\calp$ and $\len\alpha$ is at most the number of parts of $\tau$ which are greater than $2b$.
\end{itemize}

So let's write $\vo\tau\alpha0=\tau+4\alpha$ and $\vo\tau\alpha b=\tau+4\alpha\sqcup(2b)$ when $\tau$ is a \fbc and $\alpha$ a partition.

\begin{rmk}
The idea of the definition of \qs \trps is that the corresponding characters $\spn\la$ are well-behaved, inheriting properties of corresponding spin characters in Rouquier blocks. The terminology ``\qs'' is intended to reference the term ``$2$-quotient separated'' introduced by James and Mathas \cite{jmq-1} for a class of partitions labelling characters which behave in a similar way for symmetric groups and Hecke algebras.
\end{rmk}

Our aim is to prove the following, which generalises \cref{rouqirred}.

\begin{propn}\label{voddirr}
Suppose $\la$ is a \qs \trp, and write $\la=\vo\tau\alpha b$ with $b\gs0$. Then $\mspn\la$ is irreducible \iff $b\ls1$ and $\alpha$ is a $2$-Carter partition.
\end{propn}

To prove \cref{voddirr}, we use a downwards induction. The inductive step is achieved via the following \lcnamecref{voindstep}s.

\begin{lemma}\label{voindstep}
Given $l\gs1$, let $\tau=(4l-3,4l-7,\dots,1)$ and let $\upsilon=(4l-1,4l-5,\dots,3)$. Suppose $\la=\vo\tau\alpha b$ is \qs, and let $\mu=\vo\upsilon\alpha b$.
\begin{enumerate}
\item
If $b=0$, then
\[
\fmx1\spn\la=\spn\mu,\qquad\emx1\spn\mu=\spn\la.
\]
\item
If $b\gs1$, then
\[
\fmx1\spn\la_\pm=\tfrac12(\spn\mu_++\spn\mu_-),\qquad\emx1\spn\mu_\pm=\tfrac12(\spn\la_++\spn\la_-).
\]
\end{enumerate}
Hence $\mspn\la$ is irreducible \iff $\mspn\mu$ is.
\end{lemma}

\begin{pf}\indent
\begin{enumerate}
\vspace{-\topsep}
\item
$\la$ has two $1$-\spams in each of rows $1,\dots,l$. So by the induction version of \cref{spinbranchpower}, $\fmx1\spn\la=\ffd1{2l}\spn\la=\spn\mu$. Similarly, since $\mu$ has two $1$-\sprms in each row, $\emx1\spn\mu=\eed1{2l}\spn\mu=\spn\la$.
\item
We let $k$ be such that $\la_k=2b$, and consider two cases according to the parity of $b$. If $b$ is even, then $\la$ has two $1$-\spams in each of rows $1,\dots,l+1$ except row $k$.  Adding these nodes yields $\mu$. $\mu$ has no other $1$-\sprms, so removing all $1$-\sprms from $\mu$ leaves $\la$. Now the expressions for $\fmx1\spn\la$ and $\emx1\spn\mu$ follow from \cref{spinbranchpower} and its induction inversion.

If $b$ is odd, then each of rows $1,\dots,l+1$ contains two $1$-\spams of $\la$ except rows $k$ and $k+1$, which each contain one (note that the \qs condition means that $\la_{k+1}=2b-1$, so that $(k+1,\la_{k+1}+2)$ is not a \spam of $\la$). Adding all these nodes yields $\mu$, and $\mu$ has no other $1$-\sprms. Again, the expressions for $\fmx1\spn\la$ and $\emx1\spn\mu$ follow from \cref{spinbranchpower} and its induction inversion.
\end{enumerate}
In either case, we obtain
\[
\fmx1\mspn\la=\mspn\mu,\qquad\emx1\mspn\mu=\mspn\la.
\]
so the final statement follows from \cref{modbranch}.
\end{pf}

In a similar way we obtain the following.

\begin{lemma}\label{voindstep2}
Given $l\gs0$, let $\tau=(4l-1,4l-5,\dots,3)$ and $\upsilon=(4l+1,4l-3,\dots,1)$. Suppose $\la=\vo\tau\alpha b$ is \qs, and let $\mu=\vo\upsilon\alpha b$.
\begin{enumerate}
\item
If $b=0$, then
\[
\fmx0\spn\la=\spn\mu,\qquad\emx0\spn\mu=\spn\la.
\]
\item
If $b\gs1$, then
\[
\fmx0\spn\la_\pm=\tfrac12(\spn\mu_++\spn\mu_-),\qquad\emx0\spn\mu_\pm=\tfrac12(\spn\la_++\spn\la_-).
\]
\end{enumerate}
Hence $\mspn\la$ is irreducible \iff $\mspn\mu$ is.
\end{lemma}

\begin{pf}[Proof of \cref{voddirr}]
We use downwards induction on the size of the \fbc of $\la$. Note that by \cref{mainirrrouq} the result holds if $\spn\la$ lies in a Rouquier block. In general if we assume that $\la$ is \qs and write $\la=\vo\tau\alpha b$, then for any \fbc $\upsilon$ larger than $\tau$, the partition $\vo\upsilon\alpha b$ is also \qs and has the same \fbw as $\la$, i.e.\ $2|\alpha|+b$. So if we fix $\alpha$ and $b$ and allow $\tau$ to vary, then for sufficiently large $\tau$ the character $\spn{\vo\tau\alpha b}$ lies in a Rouquier block. So to prove \cref{voddirr} for a specific $\la=\vo\tau\alpha b$, it suffices to show that $\mspn\la$ is irreducible \iff $\mspn{\vo\upsilon\alpha b}$ is irreducible, for some \fbc $\upsilon$ larger than $\tau$. If $\la$ is $1$-\qs this follows from \cref{voindstep}, while if $\la$ is $3$-\qs then we use \cref{voindstep2} instead.
\end{pf}

Now we come to the conclusion of this section, which is that \cref{main} is true for spin characters labelled by \qs partitions.

\begin{pf}[Proof of \cref{main} when $\la$ is \qs{}]
Looking at the six cases listed in \cref{main}, we see that $\la$ is \qs in (and only in) cases 1--4: in cases 1 and 2 we have $\la=\vo\tau\alpha0$ with $\alpha$ a $2$-Carter partition, and cases 3 and 4 we have $\la=\vo\tau\alpha1$ with $\alpha$ a $2$-Carter partition. So when $\la$ is \qs, \cref{main} is equivalent to \cref{voddirr}.
\end{pf}

\section{Proof of \cref{main} for non-\qs partitions}\label{mainproofsec}

Now we prove \cref{main} in the case where $\la$ is not \qs. To do this, we introduce some notation for the various classes of partitions appearing in \cref{main}:

\begin{itemize}
\item
let $\thr$ denote the set of \trps $\tau+4\alpha$, where $\tau=(4l-1,4l-5,\dots,3)$ with $l\gs0$ and $\alpha$ a $2$-Carter partition with $\len\alpha\ls l$;
\item
let $\oni$ denote the set of \trps $\tau+4\alpha$, where $\tau=(4l-3,4l-7,\dots,1)$ with $l\gs1$ and $\alpha$ a $2$-Carter partition with $\len\alpha\ls l$;
\item
let $\thrt$ denote the set of \trps $\tau+4\alpha\sqcup(2)$, where $\tau=(4l-1,4l-5,\dots,3)$ with $l\gs0$ and $\alpha$ a $2$-Carter partition with $\len\alpha\ls l$;
\item
let $\onit$ denote the set of \trps $\tau+4\alpha\sqcup(2)$, where $\tau=(4l-3,4l-7,\dots,1)$ with $l\gs1$ and $\alpha$ a $2$-Carter partition with $\len\alpha\ls l-1$;
\item
let $\frst$ denote the set of \trps $(2b)$ or $(4b-2,1)$ for $b\gs2$.
\end{itemize}

Now let $\ird=\thr\cup\oni\cup\thrt\cup\onit\cup\frst\cup\{(3,2,1)\}$.

To prove the ``if'' part of \cref{main} in the case where $\la$ is not \qs, we need the following theorem of Wales.

\begin{thm}[\xcite{wal}{Theorem 7.7}]\label{wales}
Suppose $p$ is a prime. Then the $p$-modular reduction of $\spn n$ is irreducible unless $n$ is odd and $p$ divides $n$.
\end{thm}

In particular, we have $\mspn{2b}$ irreducible for every $b$. Next consider the case $\la=(4b-2,1)$ and let $\mu=(4b-2)$. Then (as we have just noted) $\mspn\mu$ is irreducible, and by examining \sprs we see that $\mspn\la=\fmx0\mspn\mu$, so by \cref{modbranch} $\mspn\la$ is irreducible too. The final case to consider to complete the ``if'' part of the proof is $\la=(3,2,1)$, for which one can just check the readily-available decomposition matrix (for example, at the Modular Atlas homepage www.math.rwth-aachen.de/homes/MOC/).

Now we proceed with the ``only if'' part of \cref{main}, which is more complex. We will prove by induction on $|\la|$ that if $\la$ is a \trp not in $\ird$ then $\mspn\la$ is a reducible $2$-modular character. The key result which underlies the combinatorics of the proof is the following (recall that we write $\rs\la i$ for the partition obtained by removing all the $i$-\sprms from $\la\in\cald$).

\begin{propn}\label{rsprop}
Suppose $\la\in\cald$ and $i\in\{0,1\}$. If $\mspn\la$ is irreducible, then so is $\mspn{\rs\la i}$.
\end{propn}

\begin{pf}
If $\mspn\la$ is irreducible, then by \cref{modbranch} so is $\emx i\mspn\la$. But by \cref{spinmaxbranch}, $\emx i\mspn\la$ equals $a\mspn{\rs\la i}$ for some $a\in\bbn$, so $\mspn{\rs\la i}$ must be irreducible (and in fact $a$ must equal $1$).
\end{pf}

So using our inductive hypothesis we may assume that for $i=0,1$ either $\rs\la i\in\ird$ or $\rs\la i=\la$. We also assume that $\la$ is not \qs, and by \cref{spinreg} we can assume that $\la$ has at most one non-zero even part.

We summarise our assumptions for ease of reference.

\smallskip
\begin{mdframed}[innerleftmargin=3pt,innerrightmargin=3pt,innertopmargin=3pt,innerbottommargin=3pt,roundcorner=5pt,innermargin=-3pt,outermargin=-3pt]
\noindent\textbf{Assumptions and notation in force for the rest of \cref{mainproofsec}:}

$\la$ is a \trp which has at most one non-zero even part, is not \qs and does not lie in $\ird$. For $i=0,1$, either $\rs\la i\in\ird$ or $\rs\la i=\la$. We let $\mu=\la\dblreg$.
\end{mdframed}

Recall that by \cref{spinreg}, if $\mspn\la$ is irreducible, then $\mspn\la=\sid\mu$. We will often exploit this by finding an induction or restriction functor (or a composition of such functors) which kills $\sid\mu$ but not $\mspn\la$, so that $\sid\mu$ and $\mspn\la$ cannot be equal.

The assumptions above have several immediate consequences.
\begin{itemize}
\item
$\rs\la1$ cannot lie in $\thr$ or $\thrt$, since partitions in $\thr\cup\thrt$ all have $1$-\sprms (with the exception of $\varnothing$, but in this case $\la$ would also equal $\varnothing$, contrary to assumption).
\item
Similarly, $\rs\la0$ cannot lie in $\oni$ or $\onit$.
\item
Neither $\rs\la0$ nor $\rs\la1$ can lie in $\frst$, since then $\la$ would lie in $\frst$ as well.
\item
Neither $\rs\la0$ nor $\rs\la1$ can equal $(3,2,1)$.
\end{itemize}

So our assumptions give $\rs\la0\in\thr\cup\thrt\cup\{\la\}$ and $\rs\la1\in\oni\cup\onit\cup\{\la\}$. We now consider several cases.

\newcounter{case}
\setcounter{case}0

\nextcase{$\rs\la0=\la$, $\rs\la1\in\oni$}

In this case, write
\[
\rs\la1=(4l-3,4l-7,\dots,1)+4\nu
\]
with $\nu$ a $2$-Carter partition of length $m\ls l$. We can obtain $\la$ from $\rs\la1$ by adding $1$-\spams, and there are two of these nodes which can be added in each row from $1$ to $l$. We must add at least one node in each of these rows, since otherwise $\la$ would have a $0$-\sprm, so that $\rs\la0\neq\la$. On the other hand, $\la$ has at most one even part, so there can be at most one row where we add exactly one node. If we add two nodes in each row, then
\[
\la=(4l-1,4l-5,\dots,3)+4\nu\in\thr,
\]
contrary to assumption, so instead there must be some $k$ such that we add one node in row $k$ and two nodes in each other row. That is,
\[
\la=(4l-1,4l-5,\dots,4l-4k+7,\quad 4l-4k+2,\quad 4l-4k-1,\dots,7,3)+4\nu.
\]
Hence
\[
\la\dbl=(2l,2l-1,\dots,2l-2k+3,\quad (2l-2k+1)^2,\quad 2l-2k,2l-2k-1,\dots,1)+2(\dup\nu).
\]
We now consider several subcases; recall that $m$ is the length of $\nu$.
\begin{enumerate}[label=(\arabic{case}\alph{enumi})]
\item
If $k>m$, then $\la$ is \qs, contrary to assumption.
\item
If $k=m=l=1$, then $\la\in\frst$, again contrary to assumption.
\item\label{subef}
If $k=m=l>1$, then
\[
\mu=(2l,2l-1,\dots,2)+2(\dup\nu).
\]
We apply the modular branching rules to compute $\ffd0{2l-1}\sid\mu$. Note that $\mu$ has an addable $0$-node at the end of each row from $1$ to $2l-1$, and (since $\nu_l>0$) a removable $0$-node at the end of row $2l$ and an addable $0$-node at the start of row $2l+1$. Hence all the addable $0$-nodes except the one in row $2l-1$ are conormal, so $\ffd0{2l-1}\sid\mu=\sid\xi$, where
\[
\xi=(2l+1,2l,\dots,5,4,2,0,1)+2(\dup\nu).
\]
Applying the modular branching rules again we obtain $\ee1\ffd0{2l-1}\sid\mu=0$. On the other hand, the spin branching rules give $\ee1\ffd0{2l-1}\mspn\la\neq0$, so we cannot possibly have $\mspn\la=\sid\mu$.
\item
If $l>k\ls m$, then $\nu_k>\nu_{k+1}$ (since $\nu$ is a $2$-Carter partition and hence $2$-regular), which means that
\[
\mu=(2l,2l-1,\dots,2l-2k+3,\quad 2l-2k+2,(2l-2k)^2,\quad 2l-2k-1,2l-2k-2,\dots,1)+2(\dup\nu).
\]
This partition has an addable $1$-node in row $2k$ and removable $1$-nodes in all other rows from $1$ to $2l$. But the fact that $k<l$ means that one of these removable nodes (namely, the one in row $2k+1$) is not normal, so $\ee1^{2l-1}\sid\mu=0$. But the spin branching rules give $\ee1^{2l-1}\mspn\la\neq0$, so $\mspn\la\neq\sid\mu$.
\end{enumerate}

\begin{eg}
Take $\la=(31,15,6)$. Then $\rs\la0=\la$ and
\[
\rs\la1=(29,13,5)=(9,5,1)+4(5,2,1),
\]
so that $k=l=m=3$. We have
\[
\la\dbl=(16,15,8,7,3,3),
\]
giving
\[
\mu=(16,15,8,7,4,2)=(6,5,4,3,2)+2(5,5,2,2,1,1).
\]
Now we examine the residues of the nodes of $\mu$ and the spin-residues of the nodes of $\la$:
\begin{align*}
\la&=\young(0110011001100110011001100110011,011001100110011,011001)
\\
\mu&=\yngres(2,16,15,8,7,4,2).
\end{align*}
The branching rules give $\ffd05\mspn\la=\mspn{33,17,6,1}$ and $\ffd05\sid\mu=\sid{17,16,9,8,4,2,1}$. And since $\ee1\mspn{33,17,6,1}\neq0=\ee1\sid{17,16,9,8,4,2,1}$ we deduce that $\mspn\la\neq\sid\mu$, so that $\mspn\la$ is reducible.
\end{eg}

\nextcase{$\rs\la0=\la$, $\rs\la1\in\onit$}

In this case, write
\[
\rs\la1=(4l+1,4l-3,\dots,5,2,1)+4\nu
\]
with $\nu$ a $2$-Carter partition of length $m\ls l$. $\la$ is obtained from $\rs\la1$ by adding $1$-\spams, and sufficiently many such nodes must be added that $\la$ has no $0$-\sprms. In particular, at least one node must be added in each of rows $1,\dots,l$, and exactly one node must be added in each of rows $l+1$ and $l+2$. But this means in particular that $\la_{l+2}=2$. Since by assumption $\la$ has at most one even part, all of $\la_1,\dots,\la_l$ must be odd, and therefore we must have
\[
\la=(4l+3,4l-1,\dots,7,3,2)+4\nu\in\thrt,
\]
contrary to assumption.

\nextcase{$\rs\la0\in\thr$, $\rs\la1=\la$}

Now write
\[
\rs\la0=(4l-1,4l-5,\dots,3)+4\nu
\]
with $\nu$ a $2$-Carter partition of length $m\ls l$. $\la$ is obtained from $\rs\la0$ by adding nodes of \spr $0$, and at least one such node must be added in each of rows $1,\dots,l$ to ensure that $\rs\la1=\la$. $\la$ can have at most one even part, so there can be at most one of these rows where we add a single $0$. We then also have the possibility of adding a $0$ in row $l+1$. If we add two nodes to each of rows $1,\dots,l$, then (regardless of whether we also add a node in row $l+1$) $\la$ is \qs, contrary to assumption. So there is $1\ls k\ls l$ such that we add a single node in row $k$ and two nodes in all other rows from $1$ to $l$.

Again we consider several subcases. First suppose $\la$ does have a node in row $l+1$; that is,
\[
\la=(4l+1,4l-3,\dots,4l-4k+9,\quad4l-4k+4,\quad4l-4k+1,4l-4k-3,\dots,1)+4\nu.
\]
Then
\[
\la\dbl=(2l+1,2l,\dots,2l-2k+4,\quad(2l-2k+2)^2,\quad2l-2k+1,2l-2k,\dots,1)+2(\dup\nu).
\]
\begin{enumerate}[label=(\arabic{case}\alph{enumi})]
\item\label{subc}
If $k\ls m$, then
\[
\mu=(2l+1,2l,\dots,2l-2k+4,\quad2l-2k+3,2l-2k+1,\quad2l-2k+1,2l-2k,\dots,1)+2(\dup\nu).
\]
and the modular branching rules give $\ee0^{2l}\sid\mu=0$. But $\ee0^{2l}\mspn\la\neq0$, so $\mspn\la\neq\sid\mu$.
\item
If $k>m$, then $\la$ is \qs, contrary to assumption.
\end{enumerate}
Now suppose $\la$ does not have a node in row $l+1$, i.e.
\[
\la=(4l+1,4l-3,\dots,4l-4k+9,\quad4l-4k+4,\quad4l-4k+1,4l-4k-3,\dots,5)+4\nu.
\]
\begin{enumerate}[resume,label=(\arabic{case}\alph{enumi})]
\item
If $m\gs k<l$, then (similarly to subcase \ref{subc}) $\ee0^{2l-1}\sid\mu=0$, while $\ee0^{2l-1}\mspn\la\neq0$, so $\mspn\la\neq\sid\mu$.
\item
If $l=1$, then $\la\in\frst$, contrary to assumption.
\item
If $1<k=l$, then (similarly to subcase \ref{subef}) $\ee0\ff1^{2l-2}$ kills $\sid\mu$ but not $\mspn\la$ so we cannot have $\mspn\la=\sid\mu$.
\item
Finally suppose $m<k<l$, and define
\[
\xi=(4l+1,4l-3,\dots,9,4)+4\nu.
\]
Then by applying \cref{dimrowrem} $k-1$ times followed by \cref{dimen1}, we obtain $\la\dblreg=\xi\dblreg$ and $\deg\spn\la>\deg\spn\xi$, so that by \cref{samerdoub} $\mspn\la$ is reducible.
\end{enumerate}

\nextcase{$\rs\la0\in\thrt$, $\rs\la1=\la$}

Now write
\[
\rs\la0=(4l-1,4l-5,\dots,3,2)+4\nu
\]
with $\nu$ a $2$-Carter partition. We obtain $\la$ from $\rs\la0$ by adding sufficiently many $0$-\spams to eliminate any $1$-\sprms. We certainly have to add the node $(l+2,1)$ (to eliminate the \sprm $(l+1,2)$). We also have to add at least one node in each of rows $1,\dots,l-1$; since $\la$ has at most one even part, we must add two nodes in each of these rows. If in addition we add two nodes in row $l$, then $\la=(4l+1,4l-3,\dots,5,2,1)+4\nu\in\onit$, contrary to assumption. So we must add no nodes in row $l$, and this means that $\nu_l=0$, since otherwise $\la$ has a $1$-\sprm.

So we have
\[
\la=(4l+1,4l-3,\dots,9,3,2,1)+4\nu.
\]
In particular, this implies $l>1$, since by assumption $\la\neq(3,2,1)$. So
\begin{align*}
\la\dbl&=(2l+1,2l,\dots,5,4,2,1^4)+2(\dup\nu),\\
\intertext{which implies that}
\mu&=(2l+1,2l,\dots,5,4,3,1)+2(\dup\nu\sqcup(1)).
\end{align*}
Now the modular branching rules (together with the fact that $l>1$) give $\ee0\ff1^{2l-2}\sid\mu=0$, while $\ee0\ff1^{2l-2}\mspn\la\neq0$.

\nextcase{$\rs\la0\in\thr\cup\thrt$, $\rs\la1\in\oni\cup\onit$}

In this case, let $b$ denote the total number of \sprms of $\la$. We claim first that the length of $\mu$ is at most $b$. This follows by considering the various possibilities for $\rs\la0$ and $\rs\la1$. Note that for each $i$ we have $\la_i=\max\left\{(\rs\la0)_i,(\rs\la1)_i\right\}$; since $\la$ can have at most two \sprms in row $i$, this means that $(\rs\la0)_i$ and $\rs\la1)_i$ differ by at most $2$. Moreover, $\len{\rs\la1}-\len{\rs\la0}$ is either $0$ or $1$. This leaves six distinct possibilities. We tabulate these below, giving the number $b$ for each, together with the length of $\la\dbl$.
\[
\begin{array}{ccccc}\hline
\rs\la0&\rs\la1&b&\len{\la\dbl}\\\hline
(4l-1,4l-5,\dots,3)+4\nu&(4l-3,4l-7,\dots,1)+4\xi&2l&2l\\
(4l-1,4l-5,\dots,3)+4\nu&(4l+1,4l-3,\dots,1)+4\xi&2l+1&2l+1\\
(4l-1,4l-5,\dots,3)+4\nu&(4l-3,4l-7,\dots,5,2,1)+4\xi&2l&2l+1\\
(4l-1,4l-5,\dots,3,2)+4\nu&(4l+1,4l-3,\dots,1)+4\xi&2l+1&2l+2\\
(4l-1,4l-5,\dots,3,2)+4\nu&(4l-3,4l-7,\dots,5,2,1)+4\xi&2l&2l+2\\
(4l-1,4l-5,\dots,3,2)+4\nu&(4l+1,4l-3,\dots,5,2,1)+4\xi&2l+1&2l+3\\\hline
\end{array}
\]
Now note that in the third and fourth cases the last two parts of $\la\dbl$ are both equal to $1$, so that the length of $\mu$ is less than the length of $\la\dbl$, and hence at most $b$. In the fifth and sixth cases, the last three parts of $\la\dbl$ are all equal to $1$, so the length of $\mu$ is at most $\len{\la\dbl}-2$, and hence less than or equal to $b$.

So our claim is proved. Now if $\mu$ has fewer than $b$ normal nodes, then for some $i\in\{0,1\}$ and some $r$ we have $\ee i^r\mspn\la\neq0$ while $\ee i^r\sid\mu=0$, so that $\mspn\la\neq\sid\mu$. So we can assume that $\mu$ has $b$ normal nodes, and since $\mu$ has length at most $b$ this means that $\mu$ must have a normal node at the end of every non-empty row. But the only way this can happen is if the nodes at the ends of the non-empty rows of $\mu$ all have the same residue (since if the nodes at the ends of rows $k$ and $k+1$ have different residues, then the one at the end of row $k+1$ is not normal). But if the removable nodes of $\mu$ all have residue $i$, then $\ee{1-i}\sid\mu=0$, while $\ee{1-i}\mspn\la\neq0$ (since by assumption $\rs\la{1-i}\neq\la$).

\begin{eg}
Take $\la=(21,11,5,2,1)$. Then
\begin{align*}
\rs\la0&=(19,11,3,2)=(11,7,3,2)+4(2,1)\in\thrt,\\
\rs\la1&=(21,9,5,2,1)=(13,9,5,2,1)+4(1)\in\onit.
\end{align*}
So we are in the final case in the table above, with $l=3$. $\la$ has seven \sprms. We have
\[
\la\dbl=(11,10,6,5,3,2,1^3),\qquad\mu=(11,10,7,5,4,2,1)
\]
so that $\mu$ can have at most seven normal nodes. But in fact $\mu$ has only three normal nodes (all of residue $0$). Hence $\ee1\mspn\la\neq0=\ee1\sid\mu$, so that $\mspn\la$ is reducible.
\end{eg}

This completes the proof of \cref{main}.

\section{Future directions}\label{futuresec}

\subsection{The alternating group}\label{ansec}

Let $\aaa n$ denote the alternating group, and define $\taaa n$ to be the inverse image of $\aaa n$ under the natural homomorphism $\tsss n\to\sss n$. Then $\taaa n$ is a double cover of $\aaa n$, and (except when $n\ls3$ or $n=6,7$) is the unique Schur cover of $\aaa n$. We can ask our main question for characters of $\taaa n$ as well, and the answer is closely related to that for $\tsss n$. In this section we make a few observations about this problem, which we hope to solve in a future paper. Throughout this section we work over fields large enough to be splitting fields for $\aaa n$.

As with $\tsss n$, the case of linear characters follows from the corresponding result for $\aaa n$, which is given in \cite[Theorem 3.1]{mfonaltred}. So here we need only be concerned with spin characters of $\taaa n$. We say that two characters of $\taaa n$ are \emph{conjugate} if we can transform one into the other by conjugating all elements of $\taaa n$ by an odd element of $\tsss n$. The ordinary irreducible spin characters of $\taaa n$ are easily constructed from those for $\tsss n$: for every $\la\in\cald^+(n)$ the restriction $\spn\la\res_{\taaa n}$ is a sum of two conjugate irreducible characters $\spn\la_+,\spn\la_-$, while for $\la\in\cald^-(n)$ the restriction $\spn\la_+\res_{\taaa n}=\spn\la_-\res_{\taaa n}$ is a self-conjugate irreducible character which we denote $\spn\la$. The set
\[
\lset{\spn\la_+,\spn\la_-}{\la\in\cald^+(n)}\cup\lset{\spn\la}{\la\in\cald^-(n)}
\]
is a complete set (without repeats) of ordinary irreducible spin characters of $\taaa n$. To help us to determine which of these characters remain irreducible in characteristic $2$, we consider restricting from $\tsss n$, using the fact that restriction to subgroups commutes with modular reduction. \emph{We now no longer allow ourselves to write $\spn\la$ to mean ``either $\spn\la_+$ or $\spn\la_-$'' when $\la\in\cald^-(n)$.} So $\spn\la$ for $\la\in\cald^-(n)$ will unambiguously mean the spin character for $\taaa n$ labelled by $\la$.

As with $\tsss n$, there are no irreducible spin characters of $\taaa n$ in characteristic $2$, so the irreducible $2$-modular characters of $\taaa n$ are the same as those for $\aaa n$. These are also constructed by restriction from $\sss n$, but to explain the situation here we need a definition: say that $\la\in\cald$ is an \emph{S-partition} if for every odd $l$ the $l$th ladder of $\la$ contains an even number of nodes. Then $\sid\la\res_{\aaa n}$ is a sum of two conjugate irreducible characters $\sid\la_+,\sid\la_-$ if $\la$ is an S-partition, and otherwise $\sid\la\res_{\aaa n}$ is a self-conjugate irreducible character. The irreducible characters arising in this way are (without repeats) all the irreducible $2$-modular characters of $\aaa n$. This result is due to Benson \cite[Theorem 1.1]{ben}, though as far as we are aware the characterisation of S-partitions we have given is new.

Now we consider what happens when we take an ordinary irreducible spin character of $\tsss n$, restrict to $\taaa n$ and reduce modulo $2$.

Take $\la\in\cald^-(n)$. Then the fact that restriction commutes with modular reduction gives
\[
\mspn\la=\modr{\spn\la_+}\res_{\taaa n}.
\]
So in order for $\mspn\la$ to be irreducible, $\modr{\spn\la_+}$ must be irreducible, and must remain irreducible on restriction to $\taaa n$. By \cref{spinreg}, this means that $\modr{\spn\la_+}=\sid{\la\dblreg}$, and that $\la\dblreg$ is not an S-partition.

In fact, it is easy to check that the second condition: an easy exercise shows that an S-partition must have even $2$-weight, whereas all the partitions in \cref{main} that lie in $\cald^-$ have odd \fbw, with the exception of the partitions $(4b)$ for $b\gs1$. If $\la=(4b)$ then $\la\dblreg=(2b+1,2b-1)$ is an S-partition. So for $\la\in\cald^-$ we conclude that $\mspn\la$ is irreducible \iff $\modr{\spn\la_+}$ is irreducible and $\la\neq(4b)$ for $b\gs1$.

Now consider $\la\in\cald^+(n)$. In this case restricting and reducing modulo $2$ gives
\[
\modr{\spn\la_+}+\modr{\spn\la_-}=\mspn\la\res_{\taaa n}.
\]
Since $\spn\la_+$ and $\spn\la_-$ are conjugate characters, $\modr{\spn\la_+}$ is irreducible \iff $\modr{\spn\la_-}$ is, and this happens \iff $\mspn\la\res_{\taaa n}$ has exactly two irreducible constituents. Hence either
\begin{itemize}
\item
$\mspn\la$ is irreducible, or
\item
$\mspn\la$ has exactly two irreducible constituents, which restrict to conjugate irreducible $2$-modular characters of $\aaa n$.
\end{itemize}
The main theorem of the present paper tells us exactly when the first situation occurs. For the second situation to occur, the two irreducible constituents of $\mspn\la$ would have to be equal in order for their restrictions to $\taaa n$ to be conjugate and irreducible, so by \cref{spinreg} we would have to have $\mspn\la=2\sid{\la\dblreg}$. In particular, this implies that $\ev\la$ would have to be $2$.

So to solve the main problem for spin characters of $\taaa n$, it remains to classify partitions $\la\in\cald$ with exactly two non-zero even parts, such that $\mspn\la=2\sid{\la\dblreg}$. Such partitions do occur; the first example is $(4,2,1)$, and we will see a family of examples below.

\subsection{Decomposition numbers for Rouquier blocks}

\begin{figure}[p]
\[
\begin{array}{r|c@{\,\,}c@{\,\,}c@{\,\,}c@{\,\,}c@{\,\,}c@{\,\,}c@{\,\,}c@{\,\,}c@{\,\,}c@{\,\,}c@{\,\,}c@{\,\,}c@{\,\,}c@{\,\,}c}
&\rt{(7)}
&\rt{(6,1)}
&\rt{(5,2)}
&\rt{(5,1^2)}
&\rt{(4,3)}
&\rt{(4,2,1)}
&\rt{(4,1^3)}
&\rt{(3^2,1)}
&\rt{(3,2^2)}
&\rt{(3,2,1^2)}
&\rt{(3,1^4)}
&\rt{(2^3,1)}
&\rt{(2^2,1^3)}
&\rt{(2,1^5)}
&\rt{(1^7)}
\\\hline
(3),(1)&\cdot&\cdot&\cdot&\cdot&\cdot&\cdot&\cdot&1&\cdot&\cdot
&\cdot&\cdot&\cdot&\cdot&\cdot\\
(2,1),(1)&\cdot&\cdot&\cdot&\cdot&\cdot&\cdot&\cdot&\cdot&\cdot&\cdot
&\cdot&\cdot&1&\cdot&\cdot\\
(1^3),(1)&\cdot&\cdot&\cdot&\cdot&\cdot&\cdot&\cdot&1&\cdot&\cdot
&\cdot&\cdot&\cdot&\cdot&1\\\hline
(2),(3)&\cdot&\cdot&\cdot&\cdot&\cdot&\cdot&\cdot&1&2&\cdot
&\cdot&1&1&\cdot&\cdot\\
(1^2),(3)&\cdot&\cdot&\cdot&\cdot&\cdot&\cdot&\cdot&1&2&\cdot
&2&1&1&1&1\\
(2),(2,1)&\cdot&\cdot&\cdot&\cdot&\cdot&\cdot&\cdot&\cdot&\cdot&\cdot
&\cdot&2&\cdot&\cdot&\cdot\\
(1^2),(2,1)&\cdot&\cdot&\cdot&\cdot&\cdot&\cdot&\cdot&\cdot&\cdot&\cdot
&\cdot&2&\cdot&2&\cdot\\\hline
(1),(5)&\cdot&\cdot&\cdot&4&\cdot&\cdot&2&1&4&1
&2&2&1&1&1\\
(1),(4,1)&\cdot&\cdot&\cdot&8&\cdot&\cdot&4&2&8&6
&4&4&2&2&\cdot\\
(1),(3,2)&\cdot&\cdot&\cdot&8&\cdot&\cdot&\cdot&2&8&6
&4&\cdot&2&\cdot&\cdot\\\hline
\varnothing,(7)&4&4&2&4&1&\cdot&2&1&2&1
&2&1&\cdot&1&1\\
\varnothing,(6,1)&8&8&12&16&2&4&4&2&8&6
&4&2&2&2&\cdot\\
\varnothing,(5,2)&16&8&16&24&6&12&4&8&12&8
&4&2&2&\cdot&\cdot\\
\varnothing,(4,3)&8&8&4&8&6&8&4&6&4&2
&\cdot&2&\cdot&\cdot&\cdot\\
\varnothing,(4,2,1)&8&\cdot&4&8&\cdot&8&\cdot&6&4&2
&\cdot&\cdot&\cdot&\cdot&\cdot
\end{array}
\]

\[
\begin{array}{r|c@{\,\,}c@{\,\,}c@{\,\,}c@{\,\,}c@{\,\,}c@{\,\,}c@{\,\,}c@{\,\,}c@{\,\,}c@{\,\,}c@{\,\,}c@{\,\,}c@{\,\,}c@{\,\,}c}
&\rt{(7)}
&\rt{(6,1)}
&\rt{(5,2)}
&\rt{(5,1^2)}
&\rt{(4,3)}
&\rt{(4,2,1)}
&\rt{(4,1^3)}
&\rt{(3^2,1)}
&\rt{(3,2^2)}
&\rt{(3,2,1^2)}
&\rt{(3,1^4)}
&\rt{(2^3,1)}
&\rt{(2^2,1^3)}
&\rt{(2,1^5)}
&\rt{(1^7)}
\\\hline
(3),(1)&\cdot&\cdot&\cdot&\cdot&\cdot&\cdot&\cdot&1&\cdot&\cdot
&\cdot&\cdot&\cdot&\cdot&\cdot\\
(2,1),(1)&\cdot&\cdot&\cdot&\cdot&\cdot&\cdot&\cdot&\cdot&\cdot&\cdot
&\cdot&\cdot&1&\cdot&\cdot\\
(1^3),(1)&\cdot&\cdot&\cdot&\cdot&\cdot&\cdot&\cdot&\cdot&\cdot&\cdot
&\cdot&\cdot&\cdot&\cdot&1\\\hline
(2),(3)&\cdot&\cdot&\cdot&\cdot&\cdot&\cdot&\cdot&1&2&\cdot
&\cdot&1&1&\cdot&\cdot\\
(1^2),(3)&\cdot&\cdot&\cdot&\cdot&\cdot&\cdot&\cdot&\cdot&\cdot&\cdot
&2&\cdot&1&1&1\\
(2),(2,1)&\cdot&\cdot&\cdot&\cdot&\cdot&\cdot&\cdot&\cdot&\cdot&\cdot
&\cdot&2&\cdot&\cdot&\cdot\\
(1^2),(2,1)&\cdot&\cdot&\cdot&\cdot&\cdot&\cdot&\cdot&\cdot&\cdot&\cdot
&\cdot&\cdot&\cdot&2&\cdot\\\hline
(1),(5)&\cdot&\cdot&\cdot&3&\cdot&\cdot&2&\cdot&2&1
&2&1&1&1&1\\
(1),(4,1)&\cdot&\cdot&\cdot&2&\cdot&\cdot&4&2&4&6
&4&2&2&2&\cdot\\
(1),(3,2)&\cdot&\cdot&\cdot&2&\cdot&\cdot&\cdot&2&4&6
&4&\cdot&2&\cdot&\cdot\\\hline
\varnothing,(7)&4&3&2&3&1&\cdot&2&\cdot&\cdot&1
&2&\cdot&\cdot&1&1\\
\varnothing,(6,1)&4&6&12&10&2&4&4&2&4&6
&4&\cdot&2&2&\cdot\\
\varnothing,(5,2)&4&2&16&16&6&12&4&8&8&8
&4&2&2&\cdot&\cdot\\
\varnothing,(4,3)&\cdot&2&4&6&6&8&4&6&4&2
&\cdot&2&\cdot&\cdot&\cdot\\
\varnothing,(4,2,1)&\cdot&\cdot&4&6&\cdot&8&\cdot&6&4&2
&\cdot&\cdot&\cdot&\cdot&\cdot\\
\end{array}
\]

\caption
{The matrices $E$ and $EA\v$ for a Rouquier block of weight $7$}\label{eav}
\end{figure}
One of the main results in this paper (\cref{allodd,allodd2}) is the calculation of some of the rows of the spin part of the decomposition matrix for a Rouquier block of $\tsss n$; specifically, the rows corresponding to \trps with no even parts greater than $2$. It would be interesting to extend this to determine the whole of the decomposition matrix for a Rouquier block. We expect this would be amenable to the same techniques, once suitable combinatorial expressions for the entries are found.

By way of example, we consider the spin part of the decomposition matrix for a Rouquier block of weight $7$. Let $\sigma$ be a $2$-core of length at least $6$ and $\tau$ the corresponding \fbc, and let $B$ be the block with $2$-core $\sigma$ and weight $7$. Then we write $E$ for the spin part of the decomposition matrix of $B$; we label the columns of $E$ by partitions of $7$, and the rows by pairs $(\alpha,\beta)$ of partitions with $\beta$ $2$-regular and $2|\alpha|+|\beta|=7$, and set
\[
E_{(\alpha,\beta)\mu}=\decs{(\tau+4\alpha\sqcup2\beta)}{(\sigma+2\mu)}.
\]
We also define $A$ to be the adjustment matrix for $\schur7$. The matrices $E$ and $EA\v$ are given in \cref{eav}.

Note that the submatrix consisting of the first three rows of $EA\v$ is the matrix $J$ from \cref{somedecoddsec}, in agreement with \cref{allodd2}. We leave the reader to try to work out the pattern in the rest of this matrix. We cannot even see why the entries of $EA\v$ are necessarily non-negative!

Another point to observe in the matrix $E$ is the sixth row, which shows a spin character whose $2$-modular reduction just has two equal composition factors. This is relevant to the discussion of the double cover of the alternating group in \cref{ansec}: it shows that if $\la$ is one of
\[
(19,7,4,3,2),(21,9,5,4,2,1),(23,11,7,4,3,2),(25,13,9,5,4,2,1),\dots
\]
then the ordinary spin characters of $\taaa n$ labelled by $\la$ remain irreducible in characteristic $2$.

\section{Index of notation}\label{indexnotnsec}

For the reader's convenience we conclude with an index of the notation we use in this paper. We provide references to the relevant subsections.

\newlength\colwi
\newlength\colwii
\newlength\colwiii
\setlength\colwi{2.5cm}
\setlength\colwiii{1cm}
\setlength\colwii\textwidth
\addtolength\colwii{-\colwi}
\addtolength\colwii{-\colwiii}
\addtolength\colwii{-1em}
\subsubsection*{Basic objects}
\vspace{-\topsep}
\begin{longtable}{@{}p{\colwi}p{\colwii}p{\colwiii}@{}}
$\bbf$&a field of characteristic $2$&\\
$\sss n$&the symmetric group of degree $n$&\ref{snsec}\\
$\aaa n$&the alternating group of degree $n$&\ref{ansec}\\
$\hhh n$&the Iwahori--Hecke algebra of $\sss n$ over $\bbc$, with $q=-1$&\ref{snsec}\\
$\schur n$&the Schur algebra of degree $n$ over $\bbf$&\ref{snsec}\\
$\qschur n$&the $q$-Schur algebra of degree $n$ over $\bbc$, with $q=-1$&\ref{snsec}\\
$\tsss n$&a double cover of $\sss n$&\ref{doublecoversec}\\
\end{longtable}

\subsubsection*{Partitions}
\vspace{-\topsep}
\begin{longtable}{@{}p{\colwi}p{\colwii}p{\colwiii}@{}}
$\calp$&the set of all partitions&\ref{partnsec}\\
$\calp(n)$&the set of all partitions of $n$&\ref{partnsec}\\
$\cald$&the set of all $2$-regular partitions&\ref{partnsec}\\
$\cald(n)$&the set of all $2$-regular partitions of $n$&\ref{partnsec}\\
$\dom$&the dominance order on $\calp(n)$&\ref{partnsec}\\
$a\la$&the partition $(a\la_1,a\la_2,\dots,)$&\ref{partnsec}\\
$\la+\mu$&the partition $(\la_1+\mu_1,\la_2+\mu_2,\dots)$&\ref{partnsec}\\
$\la\sqcup\mu$&the partition obtained by arranging all the parts of $\la$ and $\mu$ together in decreasing order&\ref{partnsec}\\
$\rs\la i$&the partition obtained by removing all the $i$-\sprms of $\la$&\ref{partnsec}\\
$(\la^{(0)},\la^{(1)})$&the $2$-quotient of $\la$&\ref{2coresec}\\
$\epsilon(\mu)$&the $2$-sign of $\mu$&\ref{2coresec}\\
$\ev\la$&the number of positive even parts of $\la$&\ref{irrdecsec}\\
$\cald^+(n)$&$\lset{\la\in\cald(n)}{\ev\la\text{ is even}}$&\ref{irrdecsec}\\
$\cald^-(n)$&$\lset{\la\in\cald(n)}{\ev\la\text{ is odd}}$&\ref{irrdecsec}\\
$a^\alpha_{\beta\gamma}$&the Littlewood--Richardson coefficient corresponding to $\alpha,\beta,\gamma\in\calp$ with $|\alpha|=|\beta|+|\gamma|$&\ref{inversesec}\\
$\kappa(\alpha,\mu)$&$a^\alpha_{\mu^{(0)}\mu^{(1)}}$ if $\mu$ has empty $2$-core, $0$ otherwise&\ref{inversesec}\\
$\la\reg$&the $2$-regularisation of $\la\in\calp(n)$&\ref{regnsec}\\
$\la\dbl$&the double of $\la\in\cald(n)$&\ref{regnsec}\\
$\la\dblreg$&$(\la\dbl)\reg$&\ref{regnsec}\\
$\mu\adr ir\la$&$\la$ is obtained from $\mu$ by adding $r$ nodes of residue $i$&\ref{brnchsec}\\
$\mu\ads ir\la$&$\la$ is obtained from $\mu$ by adding $r$ $i$-\spams&\ref{brnchsec}\\
$\thr$&the set of \trps $\tau+4\alpha$, where $\tau=(4l-1,4l-5,\dots,3)$ for $l\gs0$ and $\alpha$ is a $2$-Carter partition with $\len\alpha\ls l$&\ref{mainproofsec}\\
$\oni$&the set of \trps $\tau+4\alpha$, where $\tau=(4l-3,4l-7,\dots,1)$ for $l\gs1$ and $\alpha$ is a $2$-Carter partition with $\len\alpha\ls l$&\ref{mainproofsec}\\
$\thrt$&the set of \trps $\tau+4\alpha\sqcup(2)$, where $\tau=(4l-1,4l-5,\dots,3)$ for $l\gs0$ and $\alpha$ is a $2$-Carter partition with $\len\alpha\ls l$&\ref{mainproofsec}\\
$\onit$&the set of \trps $\tau+4\alpha\sqcup(2)$, where $\tau=(4l-3,4l-7,\dots,1)$ for $l\gs1$ and $\alpha$ is a $2$-Carter partition with $\len\alpha\ls l-1$&\ref{mainproofsec}\\
$\frst$&the set of \trps $(2b)$ or $(4b-2,1)$ for some $b\gs2$&\ref{mainproofsec}\\
$\ird$&$\thr\cup\oni\cup\thrt\cup\onit\cup\frst\cup\{(3,2,1)\}$&\ref{mainproofsec}\\
\end{longtable}

\subsubsection*{Modules and characters}
\vspace{-\topsep}
\begin{longtable}{@{}p{\colwi}p{\colwii}p{\colwiii}@{}}
$\spe\la$&the Specht module for $\sss n$ or $\hhh n$ corresponding to $\la\in\calp(n)$&\ref{irrdecsec}\\
$\jms\la$&the James module for $\sss n$ or $\hhh n$ corresponding to $\la\in\cald(n)$&\ref{irrdecsec}\\
$\weyl\la$&the Weyl module for $\schur n$ or $\qschur n$ corresponding to $\la\in\calp(n)$&\ref{irrdecsec}\\
$\schir\la$&the irreducible module for $\schur n$ or $\qschur n$ corresponding to $\la\in\calp(n)$&\ref{irrdecsec}\\
$D_{\la\mu}$&the decomposition number $\cm{\weyl\la}{\schir\mu}$ for $\schur n$ (equals the decomposition number $\cm{\spe\la}{\jms\mu}$ for $\bbf\sss n$ if $\mu\in\cald(n)$)&\ref{irrdecsec}\\
$\odc_{\la\mu}$&the decomposition number $\cm{\weyl\la}{\schir\mu}$ for $\qschur n$ (equals the decomposition number $\cm{\spe\la}{\jms\mu}$ for $\hhh n$ if $\mu\in\cald(n)$)&\ref{irrdecsec}\\
$A_{\la\mu}$&the $(\la,\mu)$ entry of the adjustment matrix for $\schur n$&\ref{irrdecsec}\\
$\ord\la$&the character of $\spe\la$ over $\bbc$&\ref{irrdecsec}\\
$\sid\la$&the Brauer character of $\jms\la$&\ref{irrdecsec}\\
$\spn\la$&the irreducible spin character labelled by $\la\in\cald^+(n)$&\ref{irrdecsec}\\
$\spn\la_+,\spn\la_-$&the irreducible spin characters labelled by $\la\in\cald^-(n)$&\ref{irrdecsec}\\
$\spn\la$&either $\spn\la_+$ or $\spn\la_-$, when $\la\in\cald^-(n)$&\ref{irrdecsec}\\
$\chm{\ \ }{\ \ }$&the standard inner product on characters&\ref{irrdecsec}\\
$\modr\chi$&the $2$-modular reduction of a character $\chi$&\ref{irrdecsec}\\
$\decs\la\mu$&the decomposition number $\cm{\mspn\la}{\sid\mu}$&\ref{irrdecsec}\\
$\prj\mu$&the character of the projective cover of $\jms\mu$&\ref{irrdecsec}\\
$\deg(\chi)$&the degree of a character $\chi$&\ref{degsec}\\
\end{longtable}

\subsubsection*{Branching rules}
\vspace{-\topsep}
\begin{longtable}{@{}p{\colwi}p{\colwii}p{\colwiii}@{}}
$\chi\res_{\tsss{n-1}}$&the restriction of $\chi$ to $\tsss{n-1}$&\ref{brnchsec}\\
$\chi\ind^{\tsss{n+1}}$&the character obtained by inducing $\chi$ to $\tsss{n+1}$&\ref{brnchsec}\\
$\ee i$&Robinson's $i$-restriction functor&\ref{brnchsec}\\
$\ff i$&Robinson's $i$-induction functor&\ref{brnchsec}\\
$\eed ir$&$\ee i^r/r!$&\ref{brnchsec}\\
$\ffd ir$&$\ff i^r/r!$&\ref{brnchsec}\\
$\epsilon_i\chi$&$\max\lset{r\gs0}{\ee i^r\chi\neq0}$&\ref{brnchsec}\\
$\varphi_i\chi$&$\max\lset{r\gs0}{\ff i^r\chi\neq0}$&\ref{brnchsec}\\
$\emx i\chi$&$\eed i{\epsilon_i\chi}\chi$&\ref{brnchsec}\\
$\fmx i\chi$&$\ffd i{\varphi_i\chi}\chi$&\ref{brnchsec}\\
\end{longtable}

\subsubsection*{Symmetric functions}
\vspace{-\topsep}
\begin{longtable}{@{}p{\colwi}p{\colwii}p{\colwiii}@{}}
$\sym$&the space of symmetric functions&\ref{symfnsec}\\
$s_\la$&the Schur function corresponding to $\la\in\calp$&\ref{symfnsec}\\
$h_r$&the $r$th complete homogeneous function&\ref{symfnsec}\\
$h_\la$&$h_{\la_1}h_{\la_2}\dots$&\ref{symfnsec}\\
$e_r$&the $r$th elementary symmetric function&\ref{symfnsec}\\
$e_\la$&$e_{\la_1}e_{\la_2}\dots$&\ref{symfnsec}\\
$\ip{\ \ }{\ \ }$&standard inner product on $\sym$&\ref{symfnsec}\\
$\ke\la\mu$&the coefficient of $e_\mu$ in $h_\la$&\ref{symfnsec}\\
\end{longtable}

\end{document}